\documentclass[11pt]{article}

\usepackage[a4paper,margin=1in]{geometry}
\usepackage{amsmath,amssymb,amsthm,mathtools}
\usepackage{mathrsfs}
\usepackage{enumitem}
\usepackage{xcolor}
\usepackage[colorlinks=true,linkcolor=black,citecolor=black,urlcolor=black]{hyperref}

\numberwithin{equation}{section}

\newtheorem{theorem}{Theorem}[section]
\newtheorem{proposition}[theorem]{Proposition}
\newtheorem{lemma}[theorem]{Lemma}
\newtheorem{corollary}[theorem]{Corollary}
\newtheorem{assumption}{Assumption}
\theoremstyle{definition}
\newtheorem{definition}[theorem]{Definition}
\newtheorem{example}[theorem]{Example}
\theoremstyle{remark}
\newtheorem{remark}[theorem]{Remark}

\newcommand{\R}{\mathbb R}

\newcommand{\F}{\mathcal F}
\newcommand{\cA}{\mathcal A}
\newcommand{\cP}{\mathcal P}
\newcommand{\cQ}{\mathcal Q}
\newcommand{\bS}{\mathsf S}
\newcommand{\bU}{\mathsf U}

\newcommand{\abs}[1]{\lvert#1\rvert}
\newcommand{\norm}[1]{\lVert#1\rVert}
\newcommand{\tr}{\operatorname{tr}}
\newcommand{\osc}{\operatorname{osc}}

\newcommand{\BUC}{\operatorname{BUC}}

\newcommand{\AuthorAddress}{}
\newcommand{\AuthorEmail}{}
\newcommand{\address}[1]{\renewcommand{\AuthorAddress}{#1}}
\newcommand{\email}[1]{\renewcommand{\AuthorEmail}{#1}}

\address{Department of Mathematical Sciences, Durham University,
Durham DH1 3LE, United Kingdom; \\

and School of Mathematics, Shandong University, Jinan 
250100, China
}
\email{huaizhong.zhao@durham.ac.uk
}

\title{Controlled nonstationary \texorpdfstring{$G$}{G}-diffusions:\\
dynamic programming, weighted Harris contraction \\
and entrance risk}
\author{Huaizhong Zhao\\[0.5em]
{\small\parbox{0.85\textwidth}{\centering\AuthorAddress}}\\[0.25em]
{\small\texttt{\AuthorEmail}}}
\date{}

\begin{document}
\maketitle

\begin{abstract}
We study controlled nonstationary \(G\)-diffusions in a canonical-law
formulation of covariance uncertainty.  The principal structural result is
an abstract weighted Harris--oscillation contraction theorem for
translation-covariant nonlinear evolution families.  A Lyapunov estimate for
a dominating terminal-perturbation bound and a local space--time
oscillation estimate for differences imply exponential contraction in a
weighted quotient norm, without positive homogeneity, subadditivity or
convexity.  For elementary sample-and-hold feedback controls and a
rectangular family of martingale laws, we construct a measurable
controlled state map and prove deterministic- and stopped-time dynamic
programming principles for the robust lower value.  Under uniform
ellipticity, differences of Bellman values belong to a common parabolic
Pucci class; a local Harnack argument and an all-policy Lyapunov estimate
verify the abstract hypotheses.  For zero running cost, the contraction
yields a unique tempered entrance nonlinear expectation and exponential
attraction.  With running cost, it yields a unique normalized
nonstationary corrector and an additive cocycle, together with exponential
convergence of relative values.  Subtracting a scalar primitive of the
cocycle gives an all-time viscosity solution of the backward
Bellman--Isaacs equation.  Autonomous and periodic coefficients yield
additive eigenpairs, stationary and periodic ergodic equations, and
convergence of long-time average costs.  The results provide a long-time
theory for optimally controlled dynamic risk evaluation under
nonstationary covariance uncertainty.
\end{abstract}



\tableofcontents

\section{Introduction}

Nonlinear expectations provide a systematic way to model probability
uncertainty when no single reference law is regarded as reliable; see, for
example, \cite{GilboaSchmeidler1989,Peng2019}.  In the \(G\)-expectation
framework, uncertainty in quadratic variation produces a canonical nonlinear
diffusion theory and a fully nonlinear parabolic equation; see
\cite{Peng2007,DenisHuPeng2011,SonerTouziZhang2011}.  Closely related
nondominated formulations arise in the theory of second-order BSDEs and
stochastic control under nonlinear kernels; see
\cite{CheriditoSonerTouziVictoir2007,SonerTouziZhang2012,
PossamaiTanZhou2018}.  On canonical path space, the same uncertainty is
represented by a nondominated family of local-martingale laws, whose stability
under conditioning and pasting is the probabilistic form of time consistency
\cite{Nutz2012,NutzVanHandel2013}.  This viewpoint is particularly useful for the control problem
considered here: the controller chooses an elementary sample-and-hold
feedback rule \(\alpha\in\cA\), while nature selects an admissible
covariance scenario after each action has been fixed from
left-endpoint information.  For a terminal cost \(\varphi\), the
robust lower value is
\[
 \bS_{s,t}\varphi(x)
 :=
 \inf_{\alpha\in\cA}\sup_{P\in\cP_{\cQ}}
 E^P\!\left[\varphi\!\left(X_t^{s,x,\alpha}\right)\right],
 \qquad s\leq t,
\]
where \(X^{s,x,\alpha}\) is the controlled nonstationary
\(G\)-diffusion started from \(x\) at calendar time \(s\), and
\(\cP_{\cQ}\) is the family of admissible martingale laws with
quadratic-variation density in the fixed covariance set \(\cQ\).
We also study the corresponding problem with accumulated running costs.

The control problem also admits a risk-measure interpretation.
Control enters through the terminal loss
\(\varphi(X_t^{s,x,\alpha})\), which depends on the chosen feedback
rule \(\alpha\).  Under the loss convention, a supremum of
expectations over a family of probability laws is a coherent risk
functional; under the payoff convention one changes sign, writing
\(\rho(\xi)=\sup_P E^P[-\xi]\)
\cite{ArtznerDelbaenEberHeath1999}.
Thus the inner supremum evaluates the coherent risk of the loss
generated by a fixed feedback rule, while the outer infimum selects
the smallest worst-case risk attainable by control.
The optimized value remains monotone and cash translation covariant,
but need not be subadditive and hence need not be a coherent risk
measure.  Stability under conditioning and pasting yields recursive
time consistency of the conditional risk evaluation; closely related
rectangularity conditions arise in recursive multiple-priors models
\cite{Riedel2004,Delbaen2006,EpsteinSchneider2003}.
The BSDE approach to recursive valuation provides another closely
related nonlinear-evaluation framework
\cite{ElKarouiPengQuenez1997}.
The problem addressed here is the long-time stability of this
optimally controlled risk evaluation under nonstationary covariance
uncertainty.

The canonical-law formulation makes the controller--nature information
structure explicit.  The controller chooses an action from the
information available at the left endpoint of each holding interval,
while nature selects an admissible covariance scenario for that
interval.  This leads to a Bellman--Isaacs operator with an outer
infimum over control actions and an inner supremum over covariance
matrices.

Related problems are studied by
Cosso and Pham~\cite{CossoPham2019}, who establish dynamic programming
and viscosity characterizations for generalized McKean--Vlasov games,
and by Achdou, Lauri\`ere and Lions~\cite{AchdouLauriereLions2021},
who derive coupled HJB--Fokker--Planck optimality systems for feedback
control of diffusions conditioned on survival.
In the present setting, the technical task is to establish
deterministic- and stopped-time dynamic programming for the prescribed
elementary canonical feedback class under nondominated covariance
laws.  The construction provides law-independent representatives of
state-path stopping rules and realises continuous approximate
covariance responses within the admissible-law family.  These results
identify the robust value with the corresponding viscosity solution
and, together with uniform parabolicity, provide the difference
estimates used in the weighted long-time contraction argument.

Dynamic programming and viscosity characterizations for finite-horizon
stochastic recursive control under \(G\)-expectation were developed
in \cite{HuJi2017}.  That paper concerns a prescribed finite
terminal horizon, for which the Lipschitz and continuity assumptions
provide the required estimates without a dissipativity condition;
infinite-horizon behaviour is not considered there.
Infinite-horizon recursive control was studied in
\cite{HuWangInfinite2018}, while ergodic \(G\)-BSDEs and their
applications to large-time asymptotics and ergodic control were
studied in \cite{HuWangErgodic2018}.  These two works impose global strong
dissipativity conditions of pairwise type, together with quantitative
compatibility conditions on the diffusion and BSDE coefficients.
For its long-time analysis, the present paper does not impose such
global pairwise dissipativity.  Instead, it combines the all-policy
Lyapunov drift condition \eqref{eq:controlled-dissipativity} with
uniform ellipticity and local Harnack estimates to obtain weighted
contraction for differences of values.  This provides a nonstationary
framework for entrance nonlinear expectations and normalized
corrector--cocycle families, together with exponential relative-value
convergence.

The covariance formulation also provides a natural interface with statistical
learning.  In the canonical representation of \(G\)-Brownian motion, the model
uncertainty is carried by the quadratic-variation density, a pathwise
characteristic naturally targeted by high-frequency covariance estimators.  In
the companion paper \cite{ZhaoLearningAmbiguity2026}, the fixed set \(\cQ\) is
replaced by covariance correspondences constructed from offline data or updated
during online operation.  Statistical coverage of these correspondences
yields high-confidence upper bounds on the population-model costs or risks of
the policies evaluated using them, while their convergence to a population
covariance set transfers to finite-horizon values and long-time control
objects.  The fixed-set dynamic programming, regularity and contraction
theory developed here provides the theoretical foundation for that data-driven
analysis.

Nonstationarity is not merely a formal generalization of an autonomous model.
Volatility, drift and running cost may vary with calendar time, as in seasonal
systems, changing market conditions and regime-dependent environments.
Time-inhomogeneous stochastic systems and their long-time statistical
structures have been studied through evolution systems of measures, periodic
and quasi-periodic measures, and entrance measures; see, for example,
\cite{DaPratoRockner2007,KunzeLorenziLunardi2010,Lorenzi2013,
FengZhao2020,FengZhaoZhong2023,FengQuZhao2026,LuZhao2026}.
Consequently, the value operators form a two-parameter evolution rather than a
homogeneous semigroup.  
For the controlled covariance-uncertain evolution studied here, we construct
entrance nonlinear expectations and, in the presence of running costs,
normalized time-dependent correctors accompanied by an additive scalar
cocycle.  The corrector family \((h_t)\) and the cocycle \(\Lambda\) satisfy
\[
 \begin{aligned}
  \bU_{s,t}^f h_t&=h_s+\Lambda(s,t),
       &&s\leq t,\\
  \Lambda(s,t)&=\Lambda(s,r)+\Lambda(r,t),
       &&s\leq r\leq t.
 \end{aligned}
\]
In the autonomous case the cocycle
reduces to \(\Lambda(s,t)=\rho(t-s)\), giving the usual additive eigenpair
\((\rho,h)\), with \(h\) solving the stationary Bellman--Isaacs equation.  In
the genuinely nonstationary case a constant eigenvalue is generally too rigid;
instead one obtains a normalized corrector family \((h_t)\) together with
\(\Lambda\).  After separating the scalar cocycle, the corrector is
characterized by an all-time infinite-horizon Bellman--Isaacs equation.  Under
periodic coefficients the cocycle further reduces, after a periodic scalar
gauge change, to a constant average cost \(\rho\) and a periodic corrector
solving a periodic ergodic Bellman--Isaacs equation.  For controlled diffusions
driven by a single Brownian reference law, a related periodic ergodic equation
is obtained in \cite{MaZhao2026} by a different method.

{
Nonstationarity and the possible loss of sublinearity under control create
two difficulties for the long-time analysis.  First, time dependence
replaces a single invariant expectation by an entrance family indexed by
calendar time.  Second, the outer infimum over feedback controls need not
preserve sublinearity, so arguments relying on the algebra of sublinear
Markov semigroups cannot be applied directly.}
The structural result of the paper is therefore formulated before the
controlled model is analysed.

We prove an abstract weighted Harris--oscillation theorem for
translation-covariant nonlinear evolution families: a Lyapunov estimate for a
dominating terminal-perturbation bound and a local oscillation reduction for
differences imply strict contraction in a weighted quotient by constants, and
hence exponential contraction in continuous time.  The theorem does not
require positive homogeneity, subadditivity or convexity.  
Its local hypothesis is formulated as a space--time oscillation
estimate for differences, rather than a common-measure minorisation
of controlled transition laws.  A common-measure minorisation is
less natural for the nonlinear value evolution considered here,
whereas the oscillation condition directly addresses differences
of value functions.  
In the present \(G\)-diffusion setting,
parabolic regularity and Harnack estimates provide the required
local oscillation reduction without constructing a common
minorising measure.
The mechanism is inspired by the Lyapunov--contraction structure of the
classical Harris theorem \cite{HairerMattingly2011}. Related abstract nonlinear
Harris--oscillation ideas and their application to uncontrolled
\(G\)-diffusions are developed in \cite{ZhaoZhao2026}.
The present result builds on that one-function Harris--oscillation mechanism.
Translation covariance supplies the centring argument, while
terminal-perturbation domination and local oscillation control of differences
provide the two-function formulation needed for controlled Bellman evolutions. Two-function weighted-seminorm contraction for nonlinear
risk maps was established by Shen et al.~\cite[Theorem~3.11]{ShenStannatObermayer2013}
under Lyapunov and generalized Doeblin assumptions, without requiring convexity
or positive homogeneity. 
The contribution here is an abstract contraction criterion based on local
space--time oscillation for nonstationary two-parameter evolutions, together
with its verification by parabolic Harnack chains for differences of
controlled covariance-uncertain diffusion values.  This yields a common
quantitative framework for entrance nonlinear expectations, nonstationary
correctors and additive cocycles.

The main contributions are as follows.
\begin{enumerate}[label=\textup{(\roman*)}]
\item We give a canonical controlled formulation with full-information
sample-and-hold feedback rules and make the information order
\(\inf_\alpha\sup_P\) explicit.  The class permits arbitrarily fine decision
grids but excludes instantaneous action--covariance loops.  A measurable
canonical state construction, analytic selection, stability under
concatenation and rectangularity yield deterministic- and stopped-time dynamic
programming principles.  Continuous approximate covariance selectors are
realized interval by interval in the viscosity proof.

\item We first prove, at an abstract level, a weighted Harris--oscillation
contraction theorem for translation-covariant nonlinear evolutions.  The
hypotheses are a two-sided terminal-perturbation domination, a Lyapunov estimate
for the dominating perturbation function, and a local oscillation contraction
for differences.  The conclusion is an
explicit strict block contraction, and hence exponential contraction, in a
weighted quotient norm.  The result applies equally to zero-cost nonlinear
expectations and to Bellman evolutions with running costs.

Under uniform ellipticity, positive-time Krylov--Safonov regularity and the
viscosity theorem of sums place differences of controlled values in a common
Pucci class.  A local Harnack-chain argument verifies the abstract oscillation
hypothesis, while an all-policy weak dissipativity estimate verifies the
Lyapunov hypothesis.

\item In the zero-running-cost problem, the contraction constructs a unique
tempered entrance nonlinear expectation and yields exponential loss of the
initial state.  With running cost, it gives a relative-value theory and
constructs a unique normalized time-dependent corrector \(h_t\) together with
an additive cocycle \(\Lambda(s,t)\), with exponential convergence of relative
values.

In the autonomous case this produces a unique additive eigenvalue \(\rho\) and
corrector \(h\), modulo constants, characterized both by
\[
 \bU_t^f h=h+\rho t
\]
and by the stationary ergodic Bellman--Isaacs equation
\[
 \rho
 =
 \inf_{a\in A}\sup_{Q\in\cQ}
 \left\{
 f(x,a)+b_{a,Q}(x)\cdot Dh
 +\frac12\tr\!\left(c_{a,Q}(x)D^2h\right)
 \right\}.
\]
{For genuinely nonstationary coefficients, writing}
\(\Lambda(s,t)=A(t)-A(s)\) and \(W(s,x)=h_s(x)-A(s)\), the long-time corrector
satisfies the all-time infinite-horizon equation
\[
 -\partial_s W-H_f(s,x,D_xW,D_x^2W)=0
\]
in the viscosity sense.  Under periodic coefficients there are a constant
average cost \(\rho\) and a periodic corrector \(\widetilde h(s,x)\) satisfying
\[
 -\partial_s\widetilde h
 -H_f(s,x,D_x\widetilde h,D_x^2\widetilde h)
 +\rho=0,
 \qquad
 \widetilde h(s+T,x)=\widetilde h(s,x).
\]
These are infinite-horizon corrector equations: unlike the finite-horizon
Bellman--Isaacs equations, they are selected by additive relations,
normalization and, in the periodic case, temporal periodicity rather than by
prescribed terminal data.
\end{enumerate}

The scope is deliberately structural.  The covariance ambiguity set is fixed
and rectangular, while the coefficients and running cost may depend on calendar
time.  State- or action-dependent covariance correspondences lie beyond this
fixed-set framework and require additional controlled martingale-problem and
regularity arguments.  The paper is organized around the abstract contraction
principle.  Section~2 states and proves the weighted Harris--oscillation
theorem.  The subsequent sections construct the canonical controlled evolution,
prove the deterministic- and stopped-time DPPs, and develop positive-time
regularity, the Bellman--Isaacs equation and the local Harnack--oscillation
mechanism.  The abstract hypotheses are then verified for the controlled
nonstationary \(G\)-diffusion.  The remaining sections construct entrance
expectations and develop the running-cost theory: the autonomous problem yields
the additive eigenpair and stationary ergodic Bellman--Isaacs equation, while
the genuinely nonstationary and periodic problems lead respectively to
additive cocycles and the corresponding all-time and periodic infinite-horizon
corrector equations.

\section[An abstract weighted Harris--oscillation principle for translation-covariant evolution families]{An abstract weighted Harris--oscillation principle\\ for translation-covariant evolution families}
\label{sec:abstract-Harris}

This section isolates the contraction mechanism used throughout the paper.
It is independent of the controlled diffusion, the canonical probability
space and the Bellman--Isaacs equation.

Let \(E\) be a nonempty set and let
\(V:E\to[0,\infty)\) satisfy \(V(x)<\infty\) for every \(x\in E\).
Define
\[
 \mathscr B_V(E)
 :=
 \left\{
   \varphi:E\to\R:
   \norm{\varphi}_V
   :=\sup_{x\in E}
      \frac{\abs{\varphi(x)}}{1+V(x)}
   <\infty
 \right\}.
\]
Let \(\mathscr X_V\subset\mathscr B_V(E)\) be a linear subspace
containing the constant functions, and suppose that the nonlinear
evolution considered below maps \(\mathscr X_V\) into itself. For \(\beta>0\), define
\begin{align}
 \norm{\varphi}_{\beta,V}
 &:=\sup_{x\in E}
   \frac{\abs{\varphi(x)}}{1+\beta V(x)},\\
 [\varphi]_{\beta,V}
 &:=\sup_{x\neq y}
 \frac{\abs{\varphi(x)-\varphi(y)}}
 {2+\beta V(x)+\beta V(y)}.
\end{align}
Whenever the right-hand side is finite, the elementary quotient identity
\begin{equation}\label{eq:controlled-quotient-identity}
 [\varphi]_{\beta,V}
 =\inf_{c\in\R}\norm{\varphi-c}_{\beta,V}
\end{equation}
identifies \([\cdot]_{\beta,V}\) with the quotient norm modulo constants.

For a nonempty set \(A\) and a real-valued function \(g\) on \(A\), define
\[
 \osc_A g
 :=
 \sup_{z,z'\in A}\abs{g(z)-g(z')}
 =
 \sup_A g-\inf_A g,
\]
with values in \([0,\infty]\).
When \(A\subset\R\times E\), this is understood as the
space--time oscillation:
\[
 \osc_A g
 =
 \sup_{(r,x),(r',y)\in A}
 \abs{g(r,x)-g(r',y)}.
\]

Let \((\mathsf U_{s,t})_{s\leq t}\) be an evolution family on
\(\mathscr X_V\):
\begin{equation}\label{eq:abstract-evolution-property}
 \mathsf U_{s,s}=\mathrm{Id},\qquad
 \mathsf U_{s,t}=\mathsf U_{s,r}\mathsf U_{r,t},
 \qquad s\leq r\leq t.
\end{equation}
Assume throughout this section that it is translation covariant,
\begin{equation}\label{eq:abstract-translation-covariance}
 \mathsf U_{s,t}(\varphi+c)=\mathsf U_{s,t}\varphi+c,
 \qquad c\in\R.
\end{equation}
No positive homogeneity, subadditivity, convexity or concavity is assumed.

{\renewcommand{\theassumption}{A.1}
\begin{assumption}
\label{ass:abstract-terminal-domination}
There are nonnegative functions
\(m_{s,t}:E\to[0,\infty)\), \(s\leq t\), such that for every
\(\varphi,\psi\in\mathscr X_V\) and every \(c,d\geq0\),
\begin{equation}\label{eq:abstract-terminal-domination}
 \abs{\varphi-\psi}\leq c+dV
 \quad\Longrightarrow\quad
 \abs{\mathsf U_{s,t}\varphi-\mathsf U_{s,t}\psi}
 \leq c+d\,m_{s,t}.
\end{equation}
{
The family \((m_{s,t})_{s\leq t}\) need not itself form an evolution family.
}
\end{assumption}
}

{\renewcommand{\theassumption}{A.2}
\begin{assumption}
\label{ass:controlled-Lyapunov}
Fix a block length \(a>0\).
There are \(\gamma\in(0,1)\), \(K\geq0\), \(\Gamma_a<\infty\) and
\(K_a<\infty\) such that, uniformly in \(s\in\R\),
\begin{align}
 m_{s,s+a}&\leq\gamma V+K,
 \label{eq:controlled-block-Lyapunov}\\
 m_{s,t}&\leq\Gamma_aV+K_a,
 \qquad 0\leq t-s\leq a.
 \label{eq:controlled-short-Lyapunov}
\end{align}
\end{assumption}
}
Choose \(R\) so that
\begin{equation}\label{eq:controlled-R-choice}
 R>\frac{2K}{1-\gamma},
 \qquad C_R:=\{V\leq R\}\neq\emptyset.
\end{equation}
{
For terminal data \(\varphi,\psi\in\mathscr X_V\), denote the difference of the evolved terminal functions by
}
{\begin{equation*}
 D_{r}^{s,a;\varphi,\psi}(x)
 :=\mathsf U_{r,s+a}\varphi(x)-\mathsf U_{r,s+a}\psi(x),
 \qquad s\leq r\leq s+a.
\end{equation*}}

{\renewcommand{\theassumption}{A.3}
\begin{assumption}
\label{ass:controlled-LOC}
For every \(s\in\R\) there is a set
\(\mathfrak Q_{s,R}\subset[s,s+a]\times E\) containing
\(\{s\}\times C_R\), and there is \(\vartheta_R\in[0,1)\), independent
of \(s,\varphi,\psi\), such that
\begin{equation}\label{eq:abstract-controlled-LOC}
 \osc_{C_R}\!\left( D_{s}^{s,a;\varphi,\psi}(\cdot)\right)
 \leq\vartheta_R
 \osc_{\mathfrak Q_{s,R}}\left (D_{\cdot}^{s,a;\varphi,\psi}(\cdot)\right ).
\end{equation}
Moreover,
\begin{equation}\label{eq:controlled-local-moment-bound}
 L_R
 :=\sup_{s\in\R}\sup_{(r,x)\in\mathfrak Q_{s,R}}
 m_{r,s+a}(x)<\infty.
\end{equation}
\end{assumption}
}
\begin{theorem}
\label{thm:controlled-weighted-Harris}
Suppose that \((\mathsf U_{s,t})_{s\leq t}\) satisfies
\eqref{eq:abstract-evolution-property} and
\eqref{eq:abstract-translation-covariance}, and that
Assumptions~\ref{ass:abstract-terminal-domination},
\ref{ass:controlled-Lyapunov} and \ref{ass:controlled-LOC} hold.
{Then there are \(\beta>0\), \(q\in(0,1)\), \(\lambda>0\) and
\(C_0<\infty\) such that}
\begin{equation}\label{eq:controlled-block-contraction}
 [\mathsf U_{s,s+a}\varphi-\mathsf U_{s,s+a}\psi]_{\beta,V}
 \leq q[\varphi-\psi]_{\beta,V},
 \qquad s\in\R,
\end{equation}
and
\begin{equation}\label{eq:controlled-continuous-contraction}
 [\mathsf U_{s,t}\varphi-\mathsf U_{s,t}\psi]_{\beta,V}
 \leq C_0e^{-\lambda(t-s)}[\varphi-\psi]_{\beta,V},
 \qquad s\leq t.
\end{equation}
\end{theorem}

\begin{proof}
Choose \(\beta>0\) so that
{\begin{equation*}
 q_{\rm in}:=\vartheta_R(1+\beta L_R)<1,
\end{equation*}}
and set
\begin{align}
 \gamma_0&:=\gamma+\frac{2K}{R}<1,\\
 q_{\rm out}
 &:=\frac{2+\beta\gamma_0R}{2+\beta R}<1,\\
 q&:=q_{\rm in}\vee q_{\rm out}.
\end{align}
Now fix \(\varphi,\psi\in\mathscr X_V\) and put
\(d_0:=[\varphi-\psi]_{\beta,V}<\infty\); finiteness follows from
\(\mathscr X_V\subset\mathscr B_V(E)\).
By
\eqref{eq:controlled-quotient-identity}, for every \(\varepsilon>0\) there
is \(c\in\R\) such that, with \(\widetilde\psi:=\psi+c\),
{\begin{equation*}
 \abs{\varphi-\widetilde\psi}
 \leq(d_0+\varepsilon)(1+\beta V).
\end{equation*}}
Translation covariance gives
\[
 \mathsf U_{r,t}\widetilde\psi
 =\mathsf U_{r,t}\psi+c,
\]
so subtracting \(c\) does not change the weighted seminorm of the output
difference.  By Assumption~\ref{ass:abstract-terminal-domination},
\begin{equation}\label{eq:abstract-weighted-difference-bound}
 \abs{\mathsf U_{r,t}\varphi(x)
       -\mathsf U_{r,t}\widetilde\psi(x)}
 \leq(d_0+\varepsilon)\bigl(1+\beta m_{r,t}(x)\bigr).
\end{equation}

Fix \(s\) and \(x,y\in E\), and write \(S:=V(x)+V(y)\).  If
\(S\leq R\), then \(x,y\in C_R\).  The local oscillation hypothesis~\ref{ass:controlled-LOC} and 
\eqref{eq:abstract-weighted-difference-bound} give
\begin{align*}
 &\abs{\bigl(\mathsf U_{s,s+a}\varphi
              -\mathsf U_{s,s+a}\widetilde\psi\bigr)(x)
          -\bigl(\mathsf U_{s,s+a}\varphi
              -\mathsf U_{s,s+a}\widetilde\psi\bigr)(y)}\\
 &\qquad\leq
 2\vartheta_R(d_0+\varepsilon)(1+\beta L_R).
\end{align*}
Since the denominator in \([\cdot]_{\beta,V}\) is at least \(2\), the
corresponding quotient is bounded by
\(q_{\rm in}(d_0+\varepsilon)\).

If \(S\geq R\), the block Lyapunov bound gives
\begin{align*}
 &\abs{\bigl(\mathsf U_{s,s+a}\varphi
              -\mathsf U_{s,s+a}\widetilde\psi\bigr)(x)
          -\bigl(\mathsf U_{s,s+a}\varphi
              -\mathsf U_{s,s+a}\widetilde\psi\bigr)(y)}\\
 &\qquad\leq(d_0+\varepsilon)
 \bigl(2+\beta\gamma S+2\beta K\bigr)\\
 &\qquad\leq(d_0+\varepsilon)
 \bigl(2+\beta\gamma_0S\bigr).
\end{align*}
The function
\[
 S\longmapsto\frac{2+\beta\gamma_0S}{2+\beta S}
\]
is decreasing, so its supremum over \(S\geq R\) equals
\(q_{\rm out}\).  Taking the supremum over \(x\neq y\), then letting
\(\varepsilon\downarrow0\), proves
\eqref{eq:controlled-block-contraction}.

For \(0\leq t-s\leq a\), the short-time Lyapunov estimate and
\eqref{eq:abstract-weighted-difference-bound} yield
\begin{equation}\label{eq:controlled-short-seminorm}
 [\mathsf U_{s,t}\varphi-\mathsf U_{s,t}\psi]_{\beta,V}
 \leq M_a[\varphi-\psi]_{\beta,V},
 \qquad
 M_a:=\max\{\Gamma_a,1+\beta K_a\}.
\end{equation}
Write \(t-s=na+r\), with \(n\geq0\) and \(0\leq r<a\).  The evolution
property, \eqref{eq:controlled-block-contraction}, and
\eqref{eq:controlled-short-seminorm} give
\[
 [\mathsf U_{s,t}\varphi-\mathsf U_{s,t}\psi]_{\beta,V}
 \leq q^nM_a[\varphi-\psi]_{\beta,V}.
\]
{Set
\begin{equation}\label{eq:controlled-contraction-constants}
 \lambda=-a^{-1}\log q,
 \qquad
 C_0=q^{-1}M_a.
\end{equation}}
Since \(n>(t-s)/a-1\), this proves
\eqref{eq:controlled-continuous-contraction}.
\end{proof}

\begin{remark}
The theorem separates the global and local mechanisms.  The functions
\(m_{s,t}\) control the propagation of terminal perturbations and supply the
Lyapunov drift outside \(C_R\), whereas
\eqref{eq:abstract-controlled-LOC} supplies the strict reduction inside
\(C_R\).  In the controlled diffusion below, \(m_{s,t}\) is the all-policy
moment bound and the local oscillation estimate is obtained from a
parabolic Harnack chain.  For Bellman evolutions with a common running cost,
the same bound controls terminal perturbations because the accumulated
running cost cancels when two terminal data are compared.
\end{remark}

\begin{remark}
\label{rem:abstract-two-functions}
The two-function formulation of
Theorem~\ref{thm:controlled-weighted-Harris} is deliberate.
For a nonlinear evolution, a contraction estimate for a single function
does not in general imply a contraction estimate for differences, since
\[
 \mathsf U_{s,t}\varphi-\mathsf U_{s,t}\psi
 \neq
 \mathsf U_{s,t}(\varphi-\psi).
\]
The passage to two terminal functions is possible here because the evolution
is translation covariant and the weighted seminorm is the quotient norm
modulo constants.

Indeed, if
\[
 d=[\varphi-\psi]_{\beta,V},
\]
then, by \eqref{eq:controlled-quotient-identity}, one may choose
\(c\in\R\) so that
\[
 \abs{\varphi-(\psi+c)}
 \leq (d+\varepsilon)(1+\beta V).
\]
Translation covariance gives
\[
 \mathsf U_{s,t}(\psi+c)
 =\mathsf U_{s,t}\psi+c,
\]
and hence by \eqref{eq:controlled-quotient-identity} again
\[
 \bigl[
   \mathsf U_{s,t}\varphi
   -\mathsf U_{s,t}(\psi+c)
 \bigr]_{\beta,V}
 =
 \bigl[
   \mathsf U_{s,t}\varphi
   -\mathsf U_{s,t}\psi
 \bigr]_{\beta,V}.
\]
Thus the terminal difference may be centred by an arbitrary constant
without changing the quantity to be contracted.  The domination assumption
then propagates the centred terminal difference, while the local oscillation
assumption contracts its state-dependent part.

This is useful beyond the zero-running-cost case.  In particular, Bellman
evolutions with a common running cost are generally neither constant
preserving nor positively homogeneous, but they remain translation
covariant.  Their differences therefore fall naturally within the present
two-function framework.

The formulation should also be compared with the weighted
Harris--oscillation result of Zhao--Zhao
\cite{ZhaoZhao2026}.  There the abstract theory is developed for a
nonlinear Markov semigroup in the sublinear-expectation setting, and the
contraction is formulated for the evolution of a single terminal function,
with the eventual aim of constructing an invariant sublinear expectation.
The present theorem uses the same Lyapunov--oscillation mechanism but
extends it in two directions relevant here: the evolution may be genuinely
nonstationary, and the contraction is between two terminal values.  The
latter extension does not rely on sublinearity or positive homogeneity;
its essential algebraic ingredient is translation covariance, which permits
one terminal function to be shifted by a constant without changing the
quotient distance of the two evolved values.  Taking \(\psi=0\) in a
zero-preserving evolution recovers the corresponding one-function
contraction.

The two-function viewpoint also has a direct antecedent in the nonlinear
risk-map theorem of Shen et al.~\cite[Theorem~3.11]{ShenStannatObermayer2013}.
Their result contracts the weighted seminorm of the difference of two
risk evaluations under a Lyapunov bound and a generalized Doeblin
condition.  Their risk maps need not be convex or positively homogeneous.
The present theorem uses terminal-perturbation domination and local
space--time oscillation reduction for a two-parameter evolution; the latter
condition will be verified by parabolic Harnack chains.  No implication
between the two sets of local hypotheses is asserted here.  Translation
covariance supplies the common centring mechanism, but the domination and
local oscillation hypotheses are also essential to the present proof.
\end{remark}


\section{The problem of controlled nonstationary \texorpdfstring{$G$}{G}-diffusion}

\subsection{Canonical uncertainty and the information structure}

\paragraph{Canonical increment space and calendar-time realization.}
Let
\[
 \Omega=C_0(\R_+;\R^m),
 \qquad
 B_u(\omega)=\omega_u,
\]
where \(B=(B_u)_{u\geq0}\) is the canonical coordinate process
{
indexed by elapsed time}, and let
\[
 \F_u:=\sigma(B_v:0\leq v\leq u),
 \qquad
 \mathbb F=(\F_u)_{u\geq0},
\]
be its canonical filtration, without completion or right-continuous augmentation.  This one-sided space is the canonical
space for future increments and will continue to be used for conditioning,
analytic selection and pasting.

{
Equip \(\Omega\) with the topology of uniform convergence on compact time
intervals, under which it is a Polish space, and denote its Borel
\(\sigma\)-field by \(\mathcal B(\Omega)\).  We write
\(\mathfrak P(\Omega)\) for the space of Borel probability measures on
\((\Omega,\mathcal B(\Omega))\), endowed with the topology of weak
convergence and its Borel \(\sigma\)-field.
}

Write \(\mathbb S_+^m\) for the set of nonnegative symmetric
\(m\times m\) matrices, and fix a nonempty compact convex set
\[
 \cQ\subset\mathbb S_+^m.
\]
{
Throughout this paper, \(\cQ\) is fixed: it does not depend on calendar
time, the state, or the controller's action.
}

Fix a progressively measurable pathwise quadratic-variation functional
\(\langle B\rangle\) which agrees with the stochastic quadratic variation
under every continuous local-martingale law.  Let
\begin{equation}\label{eq:absolutely-continuous-martingale-laws}
 \mathfrak M_a
 :=\left\{
 P\in\mathfrak P(\Omega):
 \begin{array}{l}
 B\text{ is a continuous local martingale under }P,\\
 \langle B\rangle\text{ is absolutely continuous with respect to}\ du,
 \ P\text{-a.s.}
 \end{array}
 \right\}.
\end{equation}
Following the construction around equation~(4.2) of
\cite{NutzVanHandel2013}, take the componentwise limit superior
\[
 D_u(\omega)
 :=\limsup_{n\to\infty}
 n\left(
 \langle B\rangle_u(\omega)
 -\langle B\rangle_{(u-1/n)^+}(\omega)
 \right),
 \qquad u>0,
\]
and define the progressively measurable \(\mathbb S^m\)-valued path
functional
\begin{equation}\label{eq:pathwise-covariance-density}
 \widehat Q_u(\omega)
 :=
 \begin{cases}
  D_u(\omega),
  &{
  D_u^{ij}(\omega)\in\R
   \text{ for every }1\leq i,j\leq m},\\
  0,&\text{otherwise},
 \end{cases}
 \qquad \widehat Q_0:=0.
\end{equation}
Then, for every \(P\in\mathfrak M_a\),
\begin{equation}\label{eq:pathwise-covariance-identity}
 \langle B\rangle_t
 =\int_0^t\widehat Q_u\,du
 \quad P\text{-a.s. for every }t\geq0,
 \qquad
 \widehat Q_u=\frac{d\langle B\rangle_u}{du}
 \quad du\times P\text{-a.e.}
\end{equation}
The class \(\mathfrak M_a\) is a Borel subset of
\(\mathfrak P(\Omega)\); see Lemma~4.4 of
\cite{NutzVanHandel2013}.

Define \(\cP_{\cQ}\) by
\begin{equation}\label{eq:qv}
 \cP_{\cQ}
 :=\left\{
  P\in\mathfrak M_a:
  \widehat Q_u\in\cQ
  \quad du\times P\text{-a.e.}
 \right\}.
\end{equation}
Thus each \(P\in\cP_{\cQ}\) represents one admissible absolutely continuous
covariance scenario, and \(\cP_{\cQ}\) represents the entire covariance
ambiguity family.  For every \(P\in\mathfrak M_a\), the second branch in
\eqref{eq:pathwise-covariance-density} is used only on a
\(du\times P\)-null set, by
\eqref{eq:pathwise-covariance-identity}.  Hence the particular fallback value \(0\)
does not affect either \eqref{eq:pathwise-covariance-identity} or the
\(du\times P\)-a.e. covariance constraint in \eqref{eq:qv}. No
derivative identity is asserted outside \(\mathfrak M_a\); in particular,
continuous local-martingale laws with a singular quadratic-variation
component are excluded from all covariance ambiguity families used below.

{
The state variable and the evolution operators are indexed by calendar
 time, whereas the canonical noise \(B\) and its quadratic variation are
 indexed by elapsed time.  For a problem starting at calendar time
 \(s\in\R\), the two clocks and their filtrations are related by
\begin{equation}\label{eq:calendar-noise-reindexing}
 r=s+u,\qquad u=r-s\geq0,\qquad
 \F_r^s:=\F_{r-s},\quad r\geq s.
\end{equation}
This re-indexes the information, not the driving process.  Every initial
calendar time uses the same canonical coordinate \(B\) on \(\Omega\) and
 the same admissible law family \(\cP_{\cQ}\); the starting time enters
 the coefficients through \(s+u\).  In particular, under every
 \(P\in\mathfrak M_a\), the quadratic-variation identity remains
\[
 d\langle B\rangle_u=\widehat Q_u\,du,
 \qquad du\times P\text{-a.e.}
\]
When a trajectory is restarted after elapsed time \(u\), its continuation
is driven by the future increments
\(v\mapsto B_{u+v}-B_u\), as encoded by the shift and pasting operations
below.  All conditioning, selection and pasting arguments remain on the
one-sided increment space \(\Omega\); no two-sided noise is introduced.
}

\paragraph{Universal measurability.}
For \(P\in\mathfrak P(\Omega)\), let \(\F_\tau^P\) denote the
\(P\)-completion of \(\F_\tau\).  The universal completion of
\(\F_\tau\) is
\[
 \F_\tau^*
 :=
 \bigcap_{P\in\mathfrak P(\Omega)}\F_\tau^P.
\]
A map from \(\Omega\) into a standard Borel space is called
\(\F_\tau^*\)-universally measurable if the inverse image of every Borel set
belongs to \(\F_\tau^*\). 
Recall that
\[
 \mathcal B(\mathfrak P(\Omega))
 =
 \sigma\bigl\{
  P\mapsto P(A):A\in\mathcal B(\Omega)
 \bigr\};
\]
see Bertsekas and Shreve~\cite[Proposition~7.25]{BertsekasShreve1996}.
Consequently, a map
\[
 \omega\longmapsto\nu_\omega\in\mathfrak P(\Omega)
\]
is
\(\F_\tau^*/\mathcal B(\mathfrak P(\Omega))\)-measurable if and only if
\(\omega\mapsto\nu_\omega(A)\) is \(\F_\tau^*\)-measurable for every
\(A\in\mathcal B(\Omega)\).

Universal measurability is needed because selectors obtained by projection
and analytic measurable-selection arguments need not be
\(\F_\tau\)-measurable.  For every fixed \(P\), however, an
\(\F_\tau^*\)-universally measurable kernel admits an
\(\F_\tau\)-measurable \(P\)-version.  Such a version is sufficient for
conditioning and pasting under \(P\).

\paragraph{Conditioning and pasting of probability laws.}
For a finite \(\mathbb F\)-stopping time \(\tau\), define the future
increment path by
\[
 (\theta_\tau\omega)_v
 :=
 \omega_{\tau(\omega)+v}-\omega_{\tau(\omega)},
 \qquad v\geq0.
\]
If \(P\in\mathfrak M_a\), let \(P^{\tau,\omega}\) denote a regular
conditional law of \(\theta_\tau B\) given \(\F_\tau\).  Since \(\Omega\) is
Polish, such regular conditional laws exist.  Lemma~4.7 of
\cite{NutzVanHandel2013} gives
\begin{equation}\label{eq:Ma-conditioning-and-density-shift}
 P^{\tau,\omega}\in\mathfrak M_a
 \quad\text{for \(P\)-almost every \(\omega\)},
 \qquad
 \widehat Q_v(\theta_\tau\omega)
 =\widehat Q_{\tau(\omega)+v}(\omega)
 \quad dv\times P(d\omega)\text{-a.e.}
\end{equation}
In particular, if \(P\in\cP_{\cQ}\), the defining covariance constraint is
local in time, and therefore
\[
 P^{\tau,\omega}\in\cP_{\cQ}
 \qquad\text{for \(P\)-almost every }\omega.
\]

Conversely, suppose that \(P\in\cP_{\cQ}\) and that
\[
 \omega\longmapsto\nu_\omega\in\cP_{\cQ}
\]
is an \(\F_\tau^*\)-universally measurable continuation kernel.  After
choosing an \(\F_\tau\)-measurable \(P\)-version, its pasting with \(P\) at
\(\tau\) is the probability measure \(P\otimes_\tau\nu\) defined by
\[
 (P\otimes_\tau\nu)(A)
 :=
 \int_\Omega\int_\Omega
 \mathbf 1_A
 \bigl(\omega\otimes_{\tau(\omega)}\widetilde\omega\bigr)
 \,\nu_\omega(d\widetilde\omega)\,P(d\omega),
 \qquad A\in\mathcal B(\Omega),
\]
where, for \(u\geq0\),
\[
 (\omega\otimes_u\widetilde\omega)_v
 :=
 \begin{cases}
  \omega_v,
     &v\leq u,\\
  \omega_u+\widetilde\omega_{v-u},
     &v>u.
 \end{cases}
\]
The pasted law agrees with \(P\) on the information available at time
\(\tau\), while its conditional future-increment law after observing the
past \(\omega\) is \(\nu_\omega\).  More generally, pasting a law in
\(\mathfrak M_a\) with a kernel taking values in \(\mathfrak M_a\) again
produces a law in \(\mathfrak M_a\); its quadratic-variation density is the
concatenation of the two densities.  This is the absolute-continuity part of
the pasting result in Lemma~4.9 of \cite{NutzVanHandel2013}.  Since here both
densities take values in \(\cQ\), their concatenation does also.  Consequently,
\[
 P\otimes_\tau\nu\in\cP_{\cQ}.
\]

\paragraph{Rectangularity.}
The preceding two stability properties are referred to as
\emph{rectangularity} of the ambiguity family \(\cP_{\cQ}\).  More
precisely, rectangularity means that:

\begin{enumerate}
\item the conditional future-increment law of every
      \(P\in\cP_{\cQ}\), at every finite stopping time, again belongs to
      \(\cP_{\cQ}\) almost surely;

\item pasting an admissible initial law with a measurably selected family of
      admissible continuation laws again produces a law in
      \(\cP_{\cQ}\).
\end{enumerate}

Thus nature may choose its future covariance scenario after observing the
past, without leaving the original ambiguity family.  Rectangularity is
exactly the time-consistency property that allows the worst-case expectation
to satisfy deterministic- and stopped-time dynamic programming principles.
The family \(\cP_{\cQ}\) is a Borel subset of
\(\mathfrak P(\Omega)\) and has the conditioning and pasting properties
described above; see, for example, \cite{NutzVanHandel2013}.

\paragraph{Elementary sample-and-hold feedback controls.}
Let \(A\) be a compact metric action space.  An
\emph{elementary canonical feedback rule} is specified by \(\mathbb F\)-stopping times
\[
 0=\tau_0<\tau_1<\cdots,
 \qquad \tau_k\uparrow\infty,
\]
with only finitely many switching times on every bounded time interval, and
by \(A\)-valued, \(\F_{\tau_k}\)-measurable decision variables \(\xi_k\).
The associated control is
\begin{equation}\label{eq:elementary-feedback-rule}
 \alpha_u(\omega)=\xi_k(\omega)
 \quad\text{when}\quad
 \tau_k(\omega)<u\leq\tau_{k+1}(\omega),
 \qquad \alpha_0=\xi_0.
\end{equation}
Elementary feedback strategies of this type have been used in
stochastic differential games and robust control; see
\cite{Sirbu2014Elementary,Sirbu2014Robust}.
In those formulations, an action is selected from the available
state history at a stopping rule and held until the next decision.
Related approximation results using time-discretized Markov
strategies are developed in \cite{Sirbu2015Markov}.
The present formulation uses the same sample-and-hold principle,
but the decision variables are measurable with respect to the
canonical noise history.

{Let \(\cA\) be the smallest class containing all such elementary feedback
rules and closed under countable eventwise pasting.  More precisely, at a
finite stopping time \(\tau\), the controller may select one of countably
many continuation rules \((\beta^n)_{n\geq1}\) according to an
\(\F_\tau\)-measurable partition \((C_n)_{n\geq1}\) of \(\Omega\): on
\(C_n\), the rule \(\beta^n\) is used after \(\tau\).  This operation is
defined explicitly in \eqref{eq:control-pasting}.  Every member of \(\cA\)
is \(\mathbb F\)-progressively measurable, but the converse is deliberately not
imposed.
}

{
Thus, at a decision time \(\tau_k\), the controller observes the available
history and selects \(\xi_k\).  This action is applied immediately after
\(\tau_k\) and remains in force up to and including the next decision time
\(\tau_{k+1}\).
}
Deterministic decision grids are included, as are state-feedback
rules for which \(\xi_k\) is a Borel function of the state history observed
at \(\tau_k\).  Allowing arbitrarily fine locally finite grids gives the
continuous-decision idealisation, while each individual admissible rule has
a well-defined left-endpoint information order on every holding interval.

{
For a problem started at calendar time \(s\), the calendar-time realization
of an elapsed-time rule \(\alpha\in\cA\) is
\[
 \alpha_r^{[s]}(\omega):=\alpha_{r-s}(\omega),
 \qquad r\geq s.
\]
We write
\[
 \cA_s:=\{\alpha^{[s]}:\alpha\in\cA\}.
\]
The brackets distinguish calendar-time translation from the shifted
continuation rule \(\alpha^{u,\omega}\) defined below, and we retain this
distinction throughout.  Thus \(\alpha_u\) always denotes an elapsed-time
rule, whereas \(\alpha_r^{[s]}\) denotes its calendar-time realization.
}

The controller commits to the elementary feedback rule before nature
chooses \(P\in\cP_{\cQ}\).  Operationally, at each decision time
\(\tau_k\), \(k=0,1,\ldots\), the action to be used on the next holding
interval is determined from the information available at \(\tau_k\);
only thereafter may nature select the admissible covariance scenario
over that interval. This gives the sequential information order
\[
 \inf_{\alpha\in\cA}\sup_{P\in\cP_{\cQ}},
\]
and excludes an instantaneous algebraic loop in which the action at time
\(u\) is defined as a function of the quadratic-variation density at that
same time.  It is not the reversed upper value
\(\sup_P\inf_\alpha\).

\paragraph{Shifting and pasting feedback rules.}
For \(u\geq0\), \(\omega,\widetilde\omega\in\Omega\), and
\(\alpha\in\cA\), define the shifted feedback rule by
\begin{equation}\label{eq:shifted-control}
 \alpha^{u,\omega}_v(\widetilde\omega)
 :=
 \alpha_{u+v}
 \bigl(\omega\otimes_u\widetilde\omega\bigr),
 \qquad v\geq0.
\end{equation}
It is the continuation rule obtained after the past path \(\omega\) has
been observed up to time \(u\).  If \(\tau\) is a finite stopping time, the
same definition with \(u=\tau(\omega)\) gives
\[
 \alpha^{\tau,\omega}_v(\widetilde\omega)
 :=
 \alpha_{\tau(\omega)+v}
 \bigl(\omega\otimes_{\tau(\omega)}\widetilde\omega\bigr),
 \qquad v\geq0.
\]
For every fixed past \(\omega\), this is again an admissible elementary
feedback rule: if \(u\) lies inside a holding interval, the current action
is retained until the next shifted switching time.

We also require the ability to choose different continuation controls on
different events observed at the stopping time.  Let
\(C_n\in\F_\tau\), \(n\geq1\), form a countable partition of \(\Omega\), and
let \(\beta^n\in\cA\).  Let \(\eta\in\Omega\) denote the canonical path on
which the pasted feedback rule is evaluated.  The feedback rule obtained by
using \(\alpha\) before \(\tau\) and \(\beta^n\) after \(\tau\) on \(C_n\)
is
\begin{equation}\label{eq:control-pasting}
 \overline\alpha_v(\eta)
 :=
 \begin{cases}
  \alpha_v(\eta),
     & v\leq\tau(\eta),\\[1mm]
  \beta^n_{v-\tau(\eta)}(\theta_\tau\eta),
     & v>\tau(\eta)\ \text{and }\eta\in C_n.
 \end{cases}
\end{equation}
{
Here both \(v\) and \(\tau\) are measured in elapsed time from the
initial time of the feedback rule.  If the controlled system starts at
calendar time \(s\), the corresponding calendar times are \(s+v\) and
\(s+\tau\).
}
Since \((C_n)_{n\geq1}\) is a partition,
exactly one continuation branch applies.  The stopping time \(\tau\) becomes
an additional switching time.  On each path exactly one locally finite
continuation grid is used, and standard progressive-measurability criteria
show that
\(\overline\alpha\in\cA\).

This closure of \(\cA\) under feedback pasting is the controller's
counterpart of the rectangularity of \(\cP_{\cQ}\).  In the proofs of the
dynamic programming principles, continuation controls are pasted according
to the observed state, while continuation probability laws are pasted
according to the covariance scenario selected by nature.

\begin{remark}
\label{rem:feedback}
The controller observes the canonical history generated by
\(\mathbb F\) at its decision times; the model is therefore a
full-information rather than a partial-observation formulation.
Restricting observation to the controlled state would give a
different problem.  The sample-and-hold restriction is natural
in applications: premium, reinsurance, reserve and portfolio
decisions are implemented at discrete operational times, and
covariance information is itself estimated over a nonzero
observation window.  The action is selected using current
information at each decision time; it is merely held until the
next update.  The class remains stable under shifts, countable
eventwise pastings and stopping-time concatenations, which are
the operations required by the dynamic programming principles
below.

The elementary restriction also has a mathematical role.
Under a fixed admissible law, the coefficient assumptions permit
the usual strong-solution construction for a prescribed
progressively measurable control.  Here, however, we additionally
require an admissible law realising a covariance response to the
chosen canonical feedback, as in
Lemma~\ref{lem:elementary-covariance-realisation}.
Progressive measurability alone does not ensure this property:
the pathwise covariance density \(\widehat Q_u\) is itself
progressively measurable, so an unrestricted canonical control
may depend on that density at the same time at which nature's
response is to be realised.

For example, in one noise dimension, take \(A=[0,1]\) and
\(\cQ=[1,2]\), and consider
\[
 \alpha_u(\omega)
 :=\mathbf 1_{\{\widehat Q_u(\omega)\leq3/2\}},
 \qquad
 \widehat q(q,x,a):=1+a.
\]
The rule \(\alpha\) is progressively measurable, and
\(\widehat q\) is a continuous covariance response taking values
in \(\cQ\).  If a law \(P\in\cP_{\cQ}\) realised this response
on \([0,T]\), with \(T>0\), then
\[
 \widehat Q_u
 =1+\mathbf 1_{\{\widehat Q_u\leq3/2\}}
 \quad du\times P\text{-a.e. on }[0,T].
\]
This is impossible: if \(\widehat Q_u\leq3/2\), the right-hand
side equals \(2\), whereas if \(\widehat Q_u>3/2\), it equals \(1\).
Thus the covariance-response conclusion need not hold for an
arbitrary progressively measurable canonical control.

For an elementary feedback, the action on each holding interval
is fixed from the history at its left endpoint, allowing the
response law to be constructed with that action frozen and
then concatenated over successive intervals.
The same structure is used in the canonical state construction
of Lemma~\ref{lem:controlled-canonical-state-map}.
Thus the restriction supports these constructions and the
specified controller--nature information structure; it is not
a requirement of the abstract contraction theorem itself.
\end{remark}

\begin{remark}
\label{rem:fixed-operational-grids}
If an application prescribes a deterministic grid
\(\pi=\{0=t_0<\cdots<t_N\}\), one may restrict \(\cA\) to rules whose
actions are updated only at the times in \(\pi\).  The resulting value is no
smaller than the value used in this paper, since the controller has fewer
decision opportunities.  The present class allows stopping-time grids so
that it is closed under the stopped DPP and the exit-time localization used
in the Harnack argument.  A quantitative comparison with a particular fixed
grid is a separate implementation-error question.
\end{remark}

\subsection{The controlled nonstationary \texorpdfstring{$G$}{G}-diffusion}

Let
\[
 b:\R\times\R^d\times A\to\R^d,
 \quad
 h_{ij}:\R\times\R^d\times A\to\R^d,
 \quad
 \sigma_j:\R\times\R^d\times A\to\R^d
\]
be Borel functions, with $h_{ij}=h_{ji},~i,j=1,\cdots,m$.  Put
$\sigma=(\sigma_1,\ldots,\sigma_m)$.

{\renewcommand{\theassumption}{B}
\begin{assumption}
\label{ass:coeff}
For every compact interval $I\subset\R$ there are $L_I<\infty$ and a modulus
$\rho_I$, with $\rho_I(r)\downarrow0$ as $r\downarrow0$, such that, for
$c=b,h_{ij},\sigma_j$,
\begin{align}
 \abs{c(q,x,a)-c(q,y,a)}
 &\leq L_I\abs{x-y},\label{eq:Lip}\\
 \abs{c(q,x,a)}
 &\leq L_I(1+\abs{x}),\label{eq:growth}\\
 \abs{c(q,x,a)-c(q',x,a)}
 &\leq \rho_I(\abs{q-q'})(1+\abs{x})\label{eq:time-modulus}
\end{align}
for $q,q'\in I$, $x,y\in\R^d$ and $a\in A$.
\end{assumption}
}

{
The controlled equation is motivated by an SDE driven by
\(G\)-Brownian motion.  The representation theorem of
Denis--Hu--Peng~\cite{DenisHuPeng2011} connects this viewpoint
with upper expectations over covariance-uncertain martingale
laws; the canonical formulation using \(\cP_{\cQ}\) is described
in \cite[Section~3]{NutzVanHandel2013}.
We adopt the latter viewpoint as the starting point of the
construction: for each fixed feedback rule, we construct the
state under every admissible law and evaluate its terminal
cost by taking the supremum over \(\cP_{\cQ}\).
This retains the \(G\)-expectation interpretation on its usual
domain, while allowing the pathwise-defined canonical feedback rules used
in the present paper.
}

{
For \(s\in\R\), \(x\in\R^d\) and an elapsed-time feedback rule
\(\alpha\in\cA\), retain calendar-time indexing for the state
\(X^{s,x,\alpha}=(X_r^{s,x,\alpha})_{r\geq s}\), but write its
 equation on the elapsed-time clock \(u\geq0\), using the single
 canonical process \(B\):
\begin{equation}\label{eq:controlled-SDE}
\begin{split}
 X_{s+u}^{s,x,\alpha}={}&x
 +\int_0^u
 b(s+v,X_{s+v}^{s,x,\alpha},\alpha_v)\,dv\\
 &+\sum_{i,j=1}^m\int_0^u
 h_{ij}(s+v,X_{s+v}^{s,x,\alpha},\alpha_v)
 \,d\langle B^i,B^j\rangle_v\\
 &+\sum_{j=1}^m\int_0^u
 \sigma_j(s+v,X_{s+v}^{s,x,\alpha},\alpha_v)\,dB_v^j,
 \qquad u\geq0.
\end{split}
\end{equation}
Here \(s+v\) is calendar time, while \(v\) is elapsed time for the
noise and the feedback rule.  In particular, \(X_s^{s,x,\alpha}=x\).
For \(r\geq s\), the value \(X_r^{s,x,\alpha}\) is obtained by taking
\(u=r-s\), as in \eqref{eq:calendar-noise-reindexing}.
The calendar-time control notation remains
\(\alpha_{s+v}^{[s]}=\alpha_v\).
Thus the state and the two-parameter evolution \(\bS_{s,t}\) retain their
calendar-time indices, without a start-dependent notation for the driver.
}

The next lemma supplies the canonical version of the controlled state that
is needed below.  The restriction to elementary feedbacks is essential in
its proof: deterministic-grid Euler schemes are used only while the action
is frozen, and the resulting solution maps are then concatenated at the
feedback switching times.

\begin{lemma}
\label{lem:controlled-canonical-state-map}
Let Assumption~\ref{ass:coeff} hold.  For every \(\alpha\in\cA\), there
exists a map
\begin{equation}\label{eq:canonical-controlled-state-map}
 (s,x,\omega,r)\longmapsto X_r^{s,x,\alpha}(\omega),
 \qquad
 (s,x,\omega,r)\in
 \left\{
  (s,x,\omega,r)\in\R\times\R^d\times\Omega\times\R:
  s\leq r
 \right\},
\end{equation}
which is jointly Borel.  For every fixed \((s,x)\), the map
\[
 (u,\omega)\longmapsto X_{s+u}^{s,x,\alpha}(\omega),
 \qquad u\geq0,
\]
is progressively measurable with respect to \(\mathbb F\).  Moreover, the following three properties hold:

First, under every \(P\in\cP_{\cQ}\), the process
\[
 \bigl(X_{s+u}^{s,x,\alpha}\bigr)_{u\geq0}
\]
is the unique continuous strong solution of
\eqref{eq:controlled-SDE}, up to indistinguishability.
The same conclusion holds under any \(P\in\mathfrak M_a\) for which
\(\widehat Q\) is essentially bounded on every finite horizon.
Continuity here is asserted under each such law, up to
indistinguishability, not for every canonical path \(\omega\).

Second, let \(\gamma\in\cA\), let \(\tau\) be a bounded canonical
stopping time, and let \(y:\Omega\to\R^d\) be an
\(\F_\tau\)-measurable Borel map.  For fixed \(s,t\in\R\) such that
\(0\leq\tau\leq t-s\), the continuation map
\begin{equation}\label{eq:measurable-controlled-continuation}
 (\omega,\widetilde\omega)\longmapsto
 X_t^{s+\tau(\omega),y(\omega),\gamma^{\tau,\omega}}
      (\widetilde\omega)
\end{equation}
is \(\F_\tau\otimes\mathcal B(\Omega)\)-measurable, and in
particular Borel.  Consequently, if
\(\phi:\R^d\to\overline\R\) is bounded above and upper semianalytic,
then the map
\begin{equation}\label{eq:usa-controlled-continuation}
 (\omega,P')\longmapsto
 \int_\Omega
  \phi\!\left(
   X_t^{s+\tau(\omega),y(\omega),\gamma^{\tau,\omega}}
      (\widetilde\omega)
  \right)
 \,P'(d\widetilde\omega)
\end{equation}
is upper semianalytic on
\(\Omega\times\mathfrak P(\Omega)\), and is universally measurable
with respect to \(\F_\tau\) in its first argument.  It is
\(\F_\tau\otimes\mathcal B(\mathfrak P(\Omega))\)-measurable when
\(\phi\) is bounded Borel.

Third, with \(s\), \(\gamma\), and \(\tau\) as above, fix
\(x\in\R^d\), let \(P\in\cP_{\cQ}\), and let
\(\omega\mapsto\nu_\omega\in\cP_{\cQ}\) be an
\(\F_\tau^*\)-universally measurable kernel, replaced by an
\(\F_\tau\)-measurable \(P\)-version.  Then
\begin{align}
 &X_{s+\tau(\omega)+v}^{s,x,\gamma}
   (\omega\otimes_{\tau(\omega)}\widetilde\omega)
 \notag\\
 &\qquad=
 X_{s+\tau(\omega)+v}^{
    s+\tau(\omega),X_{s+\tau(\omega)}^{s,x,\gamma}(\omega),
    \gamma^{\tau,\omega}}
   (\widetilde\omega)
 \label{eq:controlled-shifted-flow}
\end{align}
for \(P(d\omega)\nu_\omega(d\widetilde\omega)\)-almost every pair
and all \(v\geq0\).  In particular, if \(\overline\alpha\) is
obtained by using \(\alpha\) up to \(\tau\) and \(\beta^n\)
afterwards on the events of a countable
\(\F_\tau\)-measurable partition \((C_n)_{n\geq1}\), as in
\eqref{eq:control-pasting}, then
\begin{equation}\label{eq:controlled-stopped-flow}
 X_{s+\tau(\omega)+v}^{s,x,\overline\alpha}
   (\omega\otimes_{\tau(\omega)}\widetilde\omega)
 =
 X_{s+\tau(\omega)+v}^{
    s+\tau(\omega),X_{s+\tau(\omega)}^{s,x,\alpha}(\omega),
    \beta^n}
   (\widetilde\omega)
\end{equation}
on \(\{\omega\in C_n\}\), for
\(P(d\omega)\nu_\omega(d\widetilde\omega)\)-almost every pair
and all \(v\geq0\).
Both flow identities remain valid, with the same proof, for any
rectangular family of laws in \(\mathfrak M_a\) whose covariance
densities take values in a fixed compact subset of \(\mathbb S_+^m\).
\end{lemma}

\begin{proof}
We adapt the universal Euler construction for the uncontrolled equation in
\cite{ZhaoZhao2026}.  The additional point is to perform that construction
with a frozen action and then concatenate the resulting maps along the
elementary feedback rule.

\emph{Step 1: a universal frozen-action map.}
Fix \(a\in A\).  For \(\Delta_n=2^{-n}\), put
\(\eta_n(v)=\Delta_n\lfloor v/\Delta_n\rfloor\).  Starting from calendar
time \(q\) and state \(z\), define \(\Phi^{n;q,z,a}\) recursively by
\begin{align}
 \Phi_u^{n;q,z,a}
 ={}&z+
 \int_0^u b\bigl(q+\eta_n(v),
       \Phi_{\eta_n(v)}^{n;q,z,a},a\bigr)\,dv
 \notag\\
 &+\sum_{i,j=1}^m\int_0^u
 h_{ij}\bigl(q+\eta_n(v),
       \Phi_{\eta_n(v)}^{n;q,z,a},a\bigr)
 \,d\langle B^i,B^j\rangle_v
 \notag\\
 &+\sum_{j=1}^m\int_0^u
 \sigma_j\bigl(q+\eta_n(v),
       \Phi_{\eta_n(v)}^{n;q,z,a},a\bigr)\,dB_v^j.
 \label{eq:frozen-action-Euler}
\end{align}
Every integral in \eqref{eq:frozen-action-Euler} is a finite sum on compact
time intervals.  Since the coefficients are Borel, induction over the
deterministic grid shows that
\[
 (q,z,a,\omega,u)\longmapsto
 \Phi_u^{n;q,z,a}(\omega)
\]
is jointly Borel and progressively measurable in \((u,\omega)\).

{
Let \(\Phi^{P;q,z,a}\) be the usual strong solution of \eqref{eq:controlled-SDE} with the action frozen
at \(a\), under \(P\in\cP_{\cQ}\).  We give the estimate needed below.
Set
\[
 K_{\cQ}:=\sup_{Q\in\cQ}\|Q\|_{\mathrm F}<\infty
\]
and, for \(Q\in\cQ\), define
\[
 \overline b(r,x,a,Q)
 :=b(r,x,a)+\sum_{i,j=1}^m h_{ij}(r,x,a)Q_{ij}.
\]
For \(P\in\cP_{\cQ}\), write
\[
 Q_v^P:=\widehat Q_v,
 \qquad
 d\langle B\rangle_v=Q_v^P\,dv,
 \quad dv\times P\text{-a.e.}
\]
Then \(\Phi^{P;q,z,a}\) and \(\Phi^{n;q,z,a}\) satisfy, under \(P\),
the equations
\begin{align*}
 \Phi_u^{P;q,z,a}
 ={}&z+\int_0^u
 \overline b(q+v,\Phi_v^{P;q,z,a},a,Q_v^P)\,dv\\
 &+\sum_{j=1}^m\int_0^u
 \sigma_j(q+v,\Phi_v^{P;q,z,a},a)\,dB_v^j
\end{align*}
and
\begin{align*}
 \Phi_u^{n;q,z,a}
 ={}&z+\int_0^u
 \overline b\bigl(q+\eta_n(v),
        \Phi_{\eta_n(v)}^{n;q,z,a},a,Q_v^P\bigr)\,dv\\
 &+\sum_{j=1}^m\int_0^u
 \sigma_j\bigl(q+\eta_n(v),
        \Phi_{\eta_n(v)}^{n;q,z,a},a\bigr)\,dB_v^j.
\end{align*}

Let \(p\geq2\), let \(I\) be a compact calendar interval, let
\(K\subset\R^d\) be compact, and let \(H<\infty\).  Choose a compact
interval \(I_H\) containing
\(\{q+u:q\in I,\ 0\leq u\leq H\}\).  Compactness of
\(\cQ\), Assumption~\ref{ass:coeff}, the BDG inequality and Gronwall's
lemma first give the uniform moment estimate
\[
 \sup_{\substack{
       n\geq1,\ P\in\cP_{\cQ},\ q\in I,\\
       z\in K,\ a\in A}}
 E^P\!\left[
  \sup_{0\leq u\leq H}
  \left(
   |\Phi_u^{n;q,z,a}|^p+|\Phi_u^{P;q,z,a}|^p
  \right)
 \right]
 \leq C_{p,I,K,H}.
\]
On each Euler interval, the finite-variation terms have increments of order
\(\Delta_n\), whereas the martingale term has increments of order
\(\Delta_n^{1/2}\).  Consequently,
\[
 \sup_{\substack{
       P\in\cP_{\cQ},\ q\in I,\ z\in K,\ a\in A,\\
       0\leq v\leq H}}
 E^P\!\left[
  \left|
   \Phi_v^{n;q,z,a}-\Phi_{\eta_n(v)}^{n;q,z,a}
  \right|^p
 \right]
 \leq C_{p,I,K,H}\Delta_n^{p/2}.
\]

Put
\[
 e_u^n:=\Phi_u^{n;q,z,a}-\Phi_u^{P;q,z,a},
 \qquad
 F_n(u):=E^P\!\left[\sup_{0\leq r\leq u}|e_r^n|^p\right].
\]
The Lipschitz and time-continuity assumptions imply, for \(0\leq v\leq H\),
\begin{align*}
 &\left|
  \overline b\bigl(q+\eta_n(v),
      \Phi_{\eta_n(v)}^{n;q,z,a},a,Q_v^P\bigr)
  -\overline b(q+v,\Phi_v^{P;q,z,a},a,Q_v^P)
 \right|^p\\
 &\qquad\leq C\Bigl(
   |e_v^n|^p
   +\left|\Phi_v^{n;q,z,a}
      -\Phi_{\eta_n(v)}^{n;q,z,a}\right|^p
   +\rho_{I_H}(\Delta_n)^p
      \bigl(1+|\Phi_{\eta_n(v)}^{n;q,z,a}|^p\bigr)
 \Bigr),
\end{align*}
and the same bound holds for the corresponding difference of the diffusion
coefficients.  Applying H\"older's inequality to the finite-variation term
and the BDG inequality to the martingale term therefore yields
\[
 F_n(u)
 \leq C_{p,I,K,H}\int_0^u F_n(v)\,dv
 +C_{p,I,K,H}
  \left(\Delta_n^{p/2}+\rho_{I_H}(\Delta_n)^p\right).
\]
Gronwall's lemma now gives
\begin{align}
 &\sup_{\substack{P\in\cP_{\cQ},\ q\in I,\ z\in K,\ a\in A}}
 E^P\!\left[
  \sup_{0\leq u\leq H}
  \left|\Phi_u^{n;q,z,a}-\Phi_u^{P;q,z,a}\right|^p
 \right]
 \notag\\
 &\qquad\leq
 C_{p,I,K,H}
 \left(\Delta_n^{p/2}
       +\rho_{I_H}(\Delta_n)^p\right),
 \label{eq:frozen-action-Euler-error}
\end{align}
The estimate is uniform in \(a\) because the constants in
Assumption~\ref{ass:coeff} are uniform in the action.
}
Taking \(p=2\), choose, by a diagonal argument, a deterministic
subsequence for
which the right-hand side of
\eqref{eq:frozen-action-Euler-error} is summable on every compact calendar
interval and finite horizon.  Relabel this subsequence by \(n\), and define
{\begin{equation*}
 \Phi_u(q,z,a,\omega)
 :=\lim_{n\to\infty}\Phi_u^{n;q,z,a}(\omega)
\end{equation*}}
whenever the coordinatewise limit exists and is finite, and set it equal to
\(z\) otherwise.  The resulting map is jointly Borel and progressive.
Markov's inequality and Borel--Cantelli, applied to
\eqref{eq:frozen-action-Euler-error}, show that for every fixed
\((q,z,a)\) and every \(P\in\cP_{\cQ}\), the convergence is locally uniform
in \(u\), \(P\)-almost surely.  Hence the same map
\(\Phi(q,z,a,\cdot)\) is a version of the unique frozen-action strong
solution under each admissible law; the exceptional set may depend on the
law.  The identical estimate, with a constant depending on a local bound
for \(\widehat Q\), proves the same assertion for every
\(P\in\mathfrak M_a\) whose density is essentially bounded on finite
horizons.

\emph{Step 2: concatenation along an elementary rule.}
Suppose first that \(\alpha\) is represented by
\((\tau_k,\xi_k)_{k\geq0}\) as in
\eqref{eq:elementary-feedback-rule}.  Put \(Z_0^{s,x}=x\) and recursively
define
\begin{equation}\label{eq:elementary-state-recursion}
 Z_{k+1}^{s,x}(\omega)
 :=\Phi_{\tau_{k+1}(\omega)-\tau_k(\omega)}
 \left(s+\tau_k(\omega),Z_k^{s,x}(\omega),\xi_k(\omega),
       \theta_{\tau_k}\omega\right).
\end{equation}
For \(u>0\), on
\(\{\tau_k<u\leq\tau_{k+1}\}\), set
\begin{equation}\label{eq:elementary-state-map}
 {
 X_{s+u}^{s,x,\alpha}(\omega)}
 :=\Phi_{u-\tau_k(\omega)}
 \left(s+\tau_k(\omega),Z_k^{s,x}(\omega),\xi_k(\omega),
       \theta_{\tau_k}\omega\right),
 \qquad
 {%
 X_s^{s,x,\alpha}(\omega):=x.}
\end{equation}
The shift map \((u,\omega)\mapsto\theta_u\omega\) is Borel.  Thus induction
in \eqref{eq:elementary-state-recursion}, followed by the countable
partition in \eqref{eq:elementary-state-map}, proves the joint Borel and
progressive measurability asserted in
\eqref{eq:canonical-controlled-state-map}.  Local finiteness of the
switching times ensures that only finitely many recursions occur on each
bounded time interval.  The same induction, using nonanticipativity of
\(\Phi\), shows that \(Z_k^{s,x}\) is
\(\F_{\tau_k}\)-measurable.

Fix \(P\in\cP_{\cQ}\).  At \(\tau_k\), condition on \(\F_{\tau_k}\).
By rectangularity, the conditional future-increment law belongs to
\(\cP_{\cQ}\), while \(s+\tau_k\), \(Z_k^{s,x}\), and \(\xi_k\) are fixed
under the conditional problem.  Step~1 therefore identifies
\eqref{eq:elementary-state-map}, up to the next switching time, with the
unique strong solution having the frozen action \(\xi_k\).  Induction over
the finitely many switches on each finite horizon shows that
{
\(\bigl(X_{s+u}^{s,x,\alpha}\bigr)_{u\geq0}\) solves
\eqref{eq:controlled-SDE} under \(P\).}
Pathwise uniqueness follows interval by interval from
\eqref{eq:Lip}.  {
If \(\alpha\) is obtained by countable
eventwise pasting, apply the preceding construction on each member of the
measurable partition and paste the resulting maps; they agree before the
pasting time, so this}
preserves progressive measurability as well as the solution property.  This
proves the first assertion for every
\(\alpha\in\cA\).

\emph{Step 3: measurable continuation.}
The concatenation map
\((u,\omega,\widetilde\omega)\mapsto
\omega\otimes_u\widetilde\omega\) is Borel.  For a fixed elementary
representation of \(\gamma\), write
\(\eta=\omega\otimes_{\tau(\omega)}\widetilde\omega\), and let
\(k_0(\omega,\widetilde\omega)\) be the unique integer satisfying
\[
 \tau_{k_0}(\eta)\leq\tau(\omega)<\tau_{k_0+1}(\eta).
\]
This integer is Borel because each of its level sets is Borel.  Up to the
irrelevant value of the control at elapsed time zero, the shifted rule has
initial action \(\xi_{k_0}(\eta)\), and its subsequent switching times and
decisions are
\[
 \tau_{k_0+\ell}(\eta)-\tau(\omega),
 \qquad \xi_{k_0+\ell}(\eta),
 \qquad \ell\geq1.
\]
These are Borel functions of \((\omega,\widetilde\omega)\), and their
dependence on \(\omega\) is only through the stopped path
\(\omega_{\cdot\wedge\tau(\omega)}\).  Hence
the recursion \eqref{eq:elementary-state-recursion}--\eqref{eq:elementary-state-map}, with starting data
\((s+\tau(\omega),y(\omega))\), proves that
\eqref{eq:measurable-controlled-continuation} is
\(\F_\tau\otimes\mathcal B(\Omega)\)-measurable.  {
Countable
eventwise pastings are handled branchwise.}  The standard integration theorem for upper
semianalytic functions now gives
\eqref{eq:usa-controlled-continuation}, including its universal
\(\F_\tau\)-measurability; for a bounded Borel integrand, integration
against the {
variable law \(P'\)} preserves the stated product measurability.
Joint measurability in the continuation time variable follows from the same
recursion, so integrating a bounded Borel running cost in time preserves
these conclusions.

\emph{Step 4: stopped flow.}
{
Rectangularity gives
\(\overline P:=P\otimes_\tau\nu\in\cP_{\cQ}\).  Under the product measure
\(P(d\omega)\nu_\omega(d\widetilde\omega)\), disintegrate the state under
\(\overline P\) and split the elapsed-time integrals in
\eqref{eq:controlled-SDE} at \(\tau(\omega)\), corresponding to calendar
 time \(s+\tau(\omega)\).  On the pasted path
\(\eta=\omega\otimes_{\tau(\omega)}\widetilde\omega\),
\[
 B_{\tau(\omega)+v}(\eta)-B_{\tau(\omega)}(\eta)
 =B_v(\widetilde\omega),\qquad v\geq0.
\]
The continuation coefficients are evaluated at calendar time
\(s+\tau(\omega)+v\).
The quadratic variation after \(\tau(\omega)\) is likewise the
quadratic variation of this future-increment path, by the
absolute-continuity pasting property established above, and}
\[
 \gamma_{\tau(\omega)+v}
   (\omega\otimes_{\tau(\omega)}\widetilde\omega)
 =\gamma_v^{\tau,\omega}(\widetilde\omega).
\]
Thus both sides of \eqref{eq:controlled-shifted-flow} solve the same
continuation equation, with the same initial state and driver.  Pathwise
uniqueness gives equality for rational \(v\), and continuity gives equality
simultaneously for all \(v\geq0\).  On \(C_n\), the shifted continuation of
\(\overline\alpha\) is \(\beta^n\), so
\eqref{eq:controlled-stopped-flow} follows.
\end{proof}

\begin{lemma}
\label{lem:controlled-flow}
Let Assumption~\ref{ass:coeff} hold.  For every \(p\geq2\), compact calendar interval \(I\subset\R\),
and \(0<H<\infty\), there is \(C=C(p,I,H)<\infty\) such that,
for every \(s\in I\) and \(x,y\in\R^d\),
\begin{align}
 \sup_{\alpha\in\cA}\sup_{P\in\cP_{\cQ}}
 E^P\!\left[\sup_{s\leq r\leq s+H}
  \abs{X_r^{s,x,\alpha}}^p\right]
 &\leq C(1+\abs{x}^p),\label{eq:moment}\\
 \sup_{\alpha\in\cA}\sup_{P\in\cP_{\cQ}}
 E^P\!\left[\sup_{s\leq r\leq s+H}
  \abs{X_r^{s,x,\alpha}-X_r^{s,y,\alpha}}^p\right]
 &\leq C\abs{x-y}^p.\label{eq:space-estimate}
\end{align}
There is also a modulus \(\widehat\rho_{I,H,R}\) such that, whenever
\(s,s'\in I\) and \(\abs{x}\leq R\),
\begin{align}
&\sup_{\alpha\in\cA}\sup_{P\in\cP_{\cQ}}
 E^P\!\left[
  \abs{X_t^{s,x,\alpha}
       -X_t^{s',x,\alpha}}^p
 \right]^{1/p}
 \leq
 \widehat\rho_{I,H,R}(\abs{s-s'})
 \label{eq:time-estimate}
\end{align}
whenever \(0\leq t-s,t-s'\leq H\).  Consequently, by
\eqref{eq:space-estimate} and the triangle inequality, after possibly
enlarging the constant, whenever \(\abs{x},\abs{y}\leq R\),
\[
 \sup_{\alpha\in\cA}\sup_{P\in\cP_{\cQ}}
 E^P\!\left[
  \abs{X_t^{s,x,\alpha}-X_t^{s',y,\alpha}}^p
 \right]^{1/p}
 \leq
 C\abs{x-y}+\widehat\rho_{I,H,R}(\abs{s-s'}).
\]
\end{lemma}

\begin{proof} Set \(I_H:=[\inf I,\sup I+H]\), so that
\([s,s+H]\subset I_H\) for every \(s\in I\), and apply the
coefficient bounds of Assumption~\ref{ass:coeff} on \(I_H\).
Compactness of $\cQ$ bounds the quadratic-variation density uniformly.
The Burkholder--Davis--Gundy inequality, \eqref{eq:Lip}--\eqref{eq:growth}
and Gronwall's lemma give \eqref{eq:moment} and
\eqref{eq:space-estimate}, uniformly in $\alpha$ and $P$.
{
For \eqref{eq:time-estimate}, compare
\(u\mapsto X_{s+u}^{s,x,\alpha}\) and
\(u\mapsto X_{s'+u}^{s',x,\alpha}\) on their common elapsed-time
interval \(0\leq u\leq(t-s)\wedge(t-s')\), under the same law \(P\),
with the same canonical driver \(B_u\) and feedback rule \(\alpha_u\).
Only the calendar-time arguments \(s+u\) and \(s'+u\) of the
coefficients differ.  Use \eqref{eq:time-modulus} on a compact interval
\(I_H\) containing \(\{s+u:s\in I,\ 0\leq u\leq H\}\), and apply
Gronwall's lemma.  The terminal values at calendar time \(t\) correspond
to elapsed horizons \(t-s\) and \(t-s'\); their unmatched interval has
length \(|s-s'|\) and is controlled by BDG and the moment bound.
This gives \eqref{eq:time-estimate}, with a modulus built from
\(r^{1/2}+\rho_{I_H}(r)\).
}
\end{proof}

{
For \(s\leq t\), \(\varphi\in\BUC_b(\R^d)\) and \(\alpha\in\cA\), set
}
\begin{align}
 J_{s,t}^{\alpha}\varphi(x)
 &:=\sup_{P\in\cP_{\cQ}}
 E^P[\varphi(X_t^{s,x,\alpha})],\label{eq:policy-value}\\
\bS_{s,t}\varphi(x)
&:=\inf_{\alpha\in\cA}J_{s,t}^{\alpha}\varphi(x).
 \label{eq:controlled-evolution}
\end{align}
{
Only the restriction of \(\alpha\) to \([0,t-s]\) matters.}
The operator
$\bS_{s,t}$ is the lower value for a controller who moves first and covariance
nature who responds within the fixed rectangular family $\cP_{\cQ}$.

\begin{remark}\label{rem:Where-the-G-Brownian-motion-enters}
Define, on the natural canonical domain,
\[
 \widehat{\mathbb E}[\xi]
 :=
 \sup_{P\in\cP_{\cQ}}E^P[\xi],
 \qquad
 G(A):=\frac12\sup_{Q\in\cQ}\tr(QA),
 \quad A\in\mathbb S^m.
\]
By the canonical representation of \(G\)-expectation
\cite{DenisHuPeng2011}, the restriction of
\(\widehat{\mathbb E}\) to the usual \(G\)-expectation domain agrees with
the \(G\)-expectation generated by the displayed function \(G\).  In this
representation, the canonical process is a \(G\)-Brownian motion with
covariance uncertainty set \(\cQ\).  Hence, for every fixed feedback rule
\(\alpha\) and every \(\varphi\in\BUC_b(\R^d)\),
\[
 J_{s,t}^{\alpha}\varphi(x)
 =
 \widehat{\mathbb E}\!\left[
   \varphi\bigl(X_t^{s,x,\alpha}\bigr)
 \right],
\]
where the right-hand side is interpreted through the canonical family
\(\cP_{\cQ}\).  Thus \eqref{eq:controlled-SDE} is a controlled \(G\)-SDE.

For the uncontrolled \(G\)-SDE semigroup, the same canonical-law
representation extends to bounded Borel terminal functions; see
\cite{ZhaoZhao2026}.  That extension is important for the associated
capacity and ergodicity theory.  We do not require it here, since the
controlled evolution is studied on \(\BUC_b(\R^d)\) and its weighted
continuous extension.

The outer infimum over \(\alpha\) gives \(\bS_{s,t}\).  It is this
additional optimisation, rather than the \(G\)-expectation itself, that
destroys sublinearity in general.  Accordingly, the invariant-capacity
ergodicity theory for sublinear expectations does not directly apply to
the controlled evolution.  The relevant long-time objects here are instead
entrance nonlinear expectations, normalized correctors and additive
eigenvalues.

In applications, the same structure is often described without using the
term ``\(G\)-Brownian motion'': each \(P\in\cP_{\cQ}\) represents one
admissible volatility or covariance scenario, the supremum is worst-case
evaluation, and rectangularity permits the scenario to be continued
nonanticipatively after each observation time.  Thus uncertain-volatility
models and robust stochastic control provide concrete interpretations of the canonical
\(G\)-Brownian formulation.
\end{remark}

\begin{proposition}
\label{prop:controlled-continuity}
Under Assumption~\ref{ass:coeff}, for every
$\varphi\in\BUC_b(\R^d)$,
\[
 (s,t,x)\longmapsto\bS_{s,t}\varphi(x),\qquad s\leq t,
\]
is continuous locally uniformly on bounded time-state cylinders; for fixed
$s\leq t$, the function $\bS_{s,t}\varphi$ belongs to
$\BUC_b(\R^d)$.  The same assertions hold for
$J_{s,t}^{\alpha}\varphi$, with moduli that can be chosen independently of
$\alpha$.  Moreover, $\bS_{s,t}$ is monotone, translation invariant,
positively homogeneous and nonexpansive:
\[
 \bS_{s,t}(\varphi+c)=\bS_{s,t}\varphi+c,
 \qquad
 \bS_{s,t}(a\varphi)=a\bS_{s,t}\varphi,
 \quad c\in\R,\ a\geq0,
\]
and
\begin{equation}\label{eq:nonexpansive}
 \norm{\bS_{s,t}\varphi-\bS_{s,t}\psi}_\infty
 \leq\norm{\varphi-\psi}_\infty.
\end{equation}
In general it is neither sublinear nor superlinear.
\end{proposition}

\begin{proof}
Apply the modulus of continuity of $\varphi$ to
{\eqref{eq:space-estimate} and \eqref{eq:time-estimate}}, first on the event on
which the two terminal states are close and then on its complement.  The
moment estimate and Markov's inequality control the complement uniformly in
$\alpha$ and $P$.  Taking first the supremum over $P$ and then the infimum
over $\alpha$ preserves the common modulus.  Monotonicity, translation
invariance and positive homogeneity are immediate from
\eqref{eq:policy-value}--\eqref{eq:controlled-evolution}.  For every $\alpha$,
\[
 \abs{J_{s,t}^{\alpha}\varphi(x)
       -J_{s,t}^{\alpha}\psi(x)}
 \leq\norm{\varphi-\psi}_\infty;
\]
taking infima proves \eqref{eq:nonexpansive}.  The failure of additivity or
subadditivity is caused by the outer infimum over controls.
\end{proof}

\begin{remark}[Dynamic risk-measure interpretation]
\label{rem:dynamic-risk-measure-interpretation}
For a bounded terminal loss \(L\), the functional
\[
 \mathcal R(L):=\sup_{P\in\cP_{\cQ}}E^P[L]
\]
is coherent: it is monotone, cash additive, positively homogeneous and
subadditive.  Equivalently, for a terminal payoff \(\xi\),
\(\rho(\xi):=\mathcal R(-\xi)\) is a coherent monetary risk measure.  The
rectangularity of \(\cP_{\cQ}\) gives time-consistent conditional versions on
canonical path space.  Thus \(J_{s,t}^{\alpha}\varphi(x)\) is the robust risk
of the terminal loss generated by the fixed feedback rule \(\alpha\), whereas
\(\bS_{s,t}\varphi(x)\) is the minimal such risk attainable by admissible
feedback control.  Since an infimum of coherent risk functionals need not be
subadditive, Proposition~\ref{prop:controlled-continuity} deliberately treats
\(\bS_{s,t}\) as a monetary nonlinear evaluation rather than as a coherent
risk measure.  The deterministic- and stopped-time DPPs below express the
recursive time consistency of this controlled risk evaluation.
\end{remark}

\section{Dynamic programming}

\subsection{Deterministic-time dynamic programming}

For $s\leq r\leq t$, define
{\begin{equation*}
 \mathcal R_{s,r,t}\varphi(x)
 :=\inf_{\alpha\in\cA}\sup_{P\in\cP_{\cQ}}
 E^P\!\left[
  \bS_{r,t}\varphi(X_r^{s,x,\alpha})
 \right].
\end{equation*}}
{%
Again, only the part of the elapsed-time rule \(\alpha\) on
\([0,r-s]\) is used in this expression.}

\begin{theorem}
\label{thm:controlled-deterministic-DPP}
Under Assumption~\ref{ass:coeff}, for every
$\varphi\in\BUC_b(\R^d)$ and $s\leq r\leq t$,
\begin{equation}\label{eq:controlled-deterministic-DPP}
 \bS_{s,t}\varphi(x)
 =\inf_{\alpha\in\cA}\sup_{P\in\cP_{\cQ}}
 E^P\!\left[
  \bS_{r,t}\varphi(X_r^{s,x,\alpha})
 \right].
\end{equation}
Equivalently,
\begin{equation}\label{eq:controlled-evolution-property}
 \bS_{s,t}=\bS_{s,r}\bS_{r,t}
 \quad\text{on }\BUC_b(\R^d).
\end{equation}
\end{theorem}

\begin{proof}
Write $u=r-s$, $v=t-r$ and
$V(y)=\bS_{r,t}\varphi(y)$.  Proposition~\ref{prop:controlled-continuity}
shows that $V$ is bounded and continuous.

\emph{Step 1: $\bS_{s,t}\varphi\geq\mathcal R_{s,r,t}\varphi$.}
Fix a full-horizon feedback rule $\gamma\in\cA$ and a first-stage law
$P_0\in\cP_{\cQ}$.  For each past path $\omega$, let
$\gamma^{u,\omega}$ be the shifted rule \eqref{eq:shifted-control} and put
{\begin{equation*}
 H(\omega,P')
 :=\int_\Omega
  \varphi\!\left(
   X_t^{r,X_r^{s,x,\gamma}(\omega),\gamma^{u,\omega}}
   (\widetilde\omega)
  \right)
 \,P'(d\widetilde\omega),
 \qquad P'\in\mathfrak P(\Omega).
\end{equation*}}
Lemma~\ref{lem:controlled-canonical-state-map}, applied with
\(\tau\equiv u\) and
\(y(\omega)=X_r^{s,x,\gamma}(\omega)\), shows that
$(\omega,\widetilde\omega)\mapsto
X_t^{r,X_r^{s,x,\gamma}(\omega),\gamma^{u,\omega}}
(\widetilde\omega)$ is Borel.  Hence
\((\omega,P')\mapsto H(\omega,P')\) is Borel, is
\(\F_u\)-measurable in its first argument, and in
particular is upper semianalytic as a function of \((\omega,P')\).  Since
\(\cP_{\cQ}\) is Borel, analytic selection gives an
\(\F_u^*\)-universally measurable kernel
\(\omega\mapsto\nu_\omega^\varepsilon\in\cP_{\cQ}\) such that
\begin{equation}\label{eq:nature-selector-det}
H(\omega,\nu_\omega^\varepsilon)
 \geq\sup_{P'\in\cP_{\cQ}}H(\omega,P')-\varepsilon.
\end{equation}
Moreover, for every $\omega$,
\begin{equation}\label{eq:continuation-dominates-value}
\sup_{P'\in\cP_{\cQ}}H(\omega,P')
 =J_{r,t}^{\gamma^{u,\omega}}\varphi
   (X_r^{s,x,\gamma}(\omega))
 \geq V(X_r^{s,x,\gamma}(\omega)).
\end{equation}
On the other hand, after replacing the universally measurable selector by an
\(\F_u\)-measurable \(P_0\)-version, define
\[
 \overline P^\varepsilon:=P_0\otimes_u\nu^\varepsilon.
\]
Rectangularity gives \(\overline P^\varepsilon\in\cP_{\cQ}\).  By the
definition of the pasted law,
\begin{align*}
 &E^{\overline P^\varepsilon}
   \!\left[\varphi(X_t^{s,x,\gamma})\right]\\
 &\quad=
 \int_\Omega\int_\Omega
 \varphi\!\left(
  X_t^{s,x,\gamma}
  (\omega\otimes_u\widetilde\omega)
 \right)
 \nu_\omega^\varepsilon(d\widetilde\omega)\,P_0(d\omega).
\end{align*}
Since \(u=r-s\), the shifted-flow identity
\eqref{eq:controlled-shifted-flow}, with continuation length \(t-r\), gives
\[
 X_t^{s,x,\gamma}(\omega\otimes_u\widetilde\omega)
 =X_t^{r,X_r^{s,x,\gamma}(\omega),\gamma^{u,\omega}}
  (\widetilde\omega)
\]
for
\(P_0(d\omega)\nu_\omega^\varepsilon(d\widetilde\omega)\)-almost every
pair.  Consequently,
\[
 E^{\overline P^\varepsilon}
   \!\left[\varphi(X_t^{s,x,\gamma})\right]
 =E^{P_0}\!\left[H(\omega,\nu_\omega^\varepsilon)\right].
\]
It follows from \eqref{eq:nature-selector-det} and
\eqref{eq:continuation-dominates-value} that
\begin{align*}
 J_{s,t}^{\gamma}\varphi(x)
 &\geq E^{\overline P^\varepsilon}
       [\varphi(X_t^{s,x,\gamma})]\\
 &=E^{P_0}[H(\omega,\nu_\omega^\varepsilon)]\\
 &\geq E^{P_0}[V(X_r^{s,x,\gamma})]-\varepsilon.
\end{align*}
Take the supremum over $P_0$.  The restriction of $\gamma$ to $[0,u]$ is
an admissible first-stage feedback rule, so
\[
 J_{s,t}^{\gamma}\varphi(x)
 \geq\mathcal R_{s,r,t}\varphi(x)-\varepsilon.
\]
Taking the infimum over $\gamma$ and then letting
$\varepsilon\downarrow0$ proves Step 1.

\emph{Step 2: $\bS_{s,t}\varphi\leq\mathcal R_{s,r,t}\varphi$.}
For every $y\in\R^d$, choose $\beta^y\in\cA$ such that
\begin{equation}\label{eq:controller-pointwise-selector}
 J_{r,t}^{\beta^y}\varphi(y)\leq V(y)+\varepsilon.
\end{equation}
No measurable selection from the control space is being asserted here.
Both sides of \eqref{eq:controller-pointwise-selector} are continuous in the
initial state.  Therefore there is an open neighbourhood $U_y$ of $y$ such
that
{\begin{equation*}
 J_{r,t}^{\beta^y}\varphi(z)\leq V(z)+2\varepsilon,
 \qquad z\in U_y.
\end{equation*}}
Second countability gives a countable subcover $(U_{y_n})_{n\geq1}$.  Define
the Borel partition
\[
 C_1=U_{y_1},\qquad
 C_n=U_{y_n}\setminus\bigcup_{k<n}U_{y_k},\quad n\geq2,
\]
discarding empty sets, and put $\beta^n=\beta^{y_n}$.

Choose $\alpha^0\in\cA$ such that
{\begin{equation*}
 \sup_{P\in\cP_{\cQ}}E^P[V(X_r^{s,x,\alpha^0})]
 \leq\mathcal R_{s,r,t}\varphi(x)+\varepsilon.
\end{equation*}}
At calendar time $r$ (elapsed time $u$), concatenate $\alpha^0$ with
$\beta^n$ on the event
$\{X_r^{s,x,\alpha^0}\in C_n\}$.  The resulting rule
$\overline\alpha$ belongs to $\cA$.
Let \(P\in\cP_{\cQ}\), and let \(P^{u,\omega}\) be a regular conditional
law of the future increment path after elapsed time \(u\).  By
rectangularity,
\[
 P^{u,\omega}\in\cP_{\cQ}
 \qquad\text{for \(P\)-almost every \(\omega\)}.
\]
On the event \(\{X_r^{s,x,\alpha^0}\in C_n\}\), the shifted continuation
of \(\overline\alpha\) is \(\beta^n\).  Therefore, conditioning at \(u\)
and applying the controlled flow identity give
\begin{align*}
 &E^P\!\left[
   \varphi(X_t^{s,x,\overline\alpha})
   \,\middle|\,\F_u
  \right](\omega)\\
 &\quad=
 E^{P^{u,\omega}}\!\left[
  \varphi\!\left(
   X_t^{r,X_r^{s,x,\alpha^0}(\omega),\beta^n}
  \right)
 \right]\\
 &\quad\leq
 \sup_{P'\in\cP_{\cQ}}
 E^{P'}\!\left[
  \varphi\!\left(
   X_t^{r,X_r^{s,x,\alpha^0}(\omega),\beta^n}
  \right)
 \right]\\
 &\quad=
 J_{r,t}^{\beta^n}\varphi
  \bigl(X_r^{s,x,\alpha^0}(\omega)\bigr)\\
 &\quad\leq
 V\bigl(X_r^{s,x,\alpha^0}(\omega)\bigr)+2\varepsilon
\end{align*}
on \(\{X_r^{s,x,\alpha^0}\in C_n\}\).  Since these events, \(n\geq1\),
form an \(\F_u\)-measurable partition, the resulting conditional inequality
holds \(P\)-almost surely on all of \(\Omega\).  Taking expectations and
using the tower property yields
\[
 E^P\!\left[\varphi(X_t^{s,x,\overline\alpha})\right]
 \leq
 E^P\!\left[V(X_r^{s,x,\alpha^0})\right]+2\varepsilon.
\]
Since this holds for every \(P\in\cP_{\cQ}\), taking the supremum gives
\[
 J_{s,t}^{\overline\alpha}\varphi(x)
 \leq
 \sup_{P\in\cP_{\cQ}}
 E^P\!\left[V(X_r^{s,x,\alpha^0})\right]+2\varepsilon.
\]
Consequently,
\begin{align*}
 \bS_{s,t}\varphi(x)
 &\leq J_{s,t}^{\overline\alpha}\varphi(x)\\
 &\leq\sup_{P\in\cP_{\cQ}}
 E^P[V(X_r^{s,x,\alpha^0})]+2\varepsilon\\
 &\leq\mathcal R_{s,r,t}\varphi(x)+3\varepsilon.
\end{align*}
Letting \(\varepsilon\downarrow0\) proves Step 2 and
\eqref{eq:controlled-deterministic-DPP}.  Since
\(V=\bS_{r,t}\varphi\), the right-hand side of that identity is, by the
definition of \(\bS_{s,r}\),
\[
 \bS_{s,r}V(x)
 =\bS_{s,r}\bigl(\bS_{r,t}\varphi\bigr)(x).
\]
Hence \(\bS_{s,t}=\bS_{s,r}\bS_{r,t}\).
\end{proof}

\subsection{Stopped dynamic programming}

Let $s<t$ and let $\tau$ be a fixed $\mathbb F$-stopping time satisfying
$0\leq\tau\leq t-s$.  Thus \(s+\tau\) is the associated calendar-time
stopping time.  For $\alpha\in\cA$, the state
$X_{s+\tau}^{s,x,\alpha}$ is $\F_\tau$-measurable.  Define
\begin{equation}\label{eq:stopped-Bellman-RHS}
 {\mathcal R_{s,s+\tau,t}}\varphi(x)
 :=\inf_{\alpha\in\cA}\sup_{P\in\cP_{\cQ}}
 E^P\!\left[
  \bS_{s+\tau,t}\varphi(X_{s+\tau}^{s,x,\alpha})
 \right].
\end{equation}
The integrand is measurable because
$(a,y)\mapsto\bS_{a,t}\varphi(y)$ is continuous.

\begin{theorem}
\label{thm:controlled-stopped-DPP}
Under Assumption~\ref{ass:coeff}, for every
$\varphi\in\BUC_b(\R^d)$,
\begin{equation}\label{eq:controlled-stopped-DPP}
 \bS_{s,t}\varphi(x)
 =\inf_{\alpha\in\cA}\sup_{P\in\cP_{\cQ}}
 E^P\!\left[
  \bS_{s+\tau,t}\varphi(X_{s+\tau}^{s,x,\alpha})
 \right].
\end{equation}
For $\tau\equiv r-s$, this is
Theorem~\ref{thm:controlled-deterministic-DPP}.
\end{theorem}

\begin{proof}
Put
\[
 V(a,y):=\bS_{a,t}\varphi(y),\qquad a\leq t.
\]
This function is bounded and continuous by
Proposition~\ref{prop:controlled-continuity}.

\emph{Step 1: $\bS_{s,t}\varphi\geq
{\mathcal R_{s,s+\tau,t}}\varphi$.}
Fix $\gamma\in\cA$ and $P_0\in\cP_{\cQ}$.  For a past path $\omega$, set
\[
 a(\omega)=s+\tau(\omega),\qquad
 y(\omega)=X_{s+\tau}^{s,x,\gamma}(\omega)
\]
and let $\gamma^{\tau,\omega}$ be the shifted continuation rule. Define
{\begin{equation*}
 H(\omega,P')
 :=\int_\Omega
   \varphi\!\left(
    X_t^{a(\omega),y(\omega),\gamma^{\tau,\omega}}
      (\widetilde\omega)
   \right)
 \,P'(d\widetilde\omega),
 \qquad P'\in\mathfrak P(\Omega).
\end{equation*}}
Lemma~\ref{lem:controlled-canonical-state-map}, applied with
\(y(\omega)=X_{s+\tau}^{s,x,\gamma}(\omega)\), shows that
\((\omega,P')\mapsto H(\omega,P')\) is Borel and is
\(\F_\tau\)-measurable in its first argument.  In particular, it is upper
semianalytic.  Since \(\cP_{\cQ}\) is Borel, analytic selection provides an
\(\F_\tau^*\)-universally measurable kernel
\(\omega\mapsto\nu_\omega^\varepsilon\in\cP_{\cQ}\) satisfying
{\begin{equation*}
 H(\omega,\nu_\omega^\varepsilon)
 \geq\sup_{P'\in\cP_{\cQ}}H(\omega,P')-\varepsilon.
\end{equation*}}
Moreover,
\[
 \sup_{P'\in\cP_{\cQ}}H(\omega,P')
 =J_{a(\omega),t}^{\gamma^{\tau,\omega}}\varphi(y(\omega))
 \geq V(a(\omega),y(\omega)).
\]
After replacing the selector by an \(\F_\tau\)-measurable
\(P_0\)-version, define
\[
 \overline P^\varepsilon:=P_0\otimes_\tau\nu^\varepsilon.
\]
Rectangularity gives \(\overline P^\varepsilon\in\cP_{\cQ}\).  By the
same disintegration-and-pasting calculation as in the proof of
Theorem~\ref{thm:controlled-deterministic-DPP}, now applied at \(\tau\)
and using the stopped flow identity,
\[
 E^{\overline P^\varepsilon}
  [\varphi(X_t^{s,x,\gamma})]
 =E^{P_0}[H(\omega,\nu_\omega^\varepsilon)].
\]
Consequently,
\begin{align*}
 J_{s,t}^{\gamma}\varphi(x)
 &\geq E^{\overline P^\varepsilon}
       [\varphi(X_t^{s,x,\gamma})]\\
 &=E^{P_0}[H(\omega,\nu_\omega^\varepsilon)]\\
 &\geq E^{P_0}[V(s+\tau,X_{s+\tau}^{s,x,\gamma})]
       -\varepsilon.
\end{align*}
Take the supremum over $P_0$, then the infimum over $\gamma$, and use the
fact that the stopped restriction of every $\gamma$ is admissible in
\eqref{eq:stopped-Bellman-RHS}.  Letting $\varepsilon\downarrow0$ proves
Step 1.

\emph{Step 2: $\bS_{s,t}\varphi\leq
{\mathcal R_{s,s+\tau,t}}\varphi$.}
For each $(a,y)\in(-\infty,t]\times\R^d$, choose
$\beta^{a,y}\in\cA$ such that
{\begin{equation*}
 J_{a,t}^{\beta^{a,y}}\varphi(y)
 \leq V(a,y)+\varepsilon.
\end{equation*}}
For fixed $\beta^{a,y}$, both sides are continuous at $(a,y)$.  Hence there
is an open neighbourhood $U_{a,y}$ such that
{\begin{equation*}
 J_{q,t}^{\beta^{a,y}}\varphi(z)
 \leq V(q,z)+2\varepsilon,
 \qquad (q,z)\in U_{a,y}.
\end{equation*}}
Choose a countable subcover $(U_n)_{n\geq1}$ of
$[s,t]\times\R^d$, choose the associated controls $\beta^n$, and form a
Borel partition $(D_n)_{n\geq1}$ subordinate to this cover.

Choose $\alpha^0\in\cA$ satisfying
{\begin{align*}
 &\sup_{P\in\cP_{\cQ}}
 E^P[V(s+\tau,X_{s+\tau}^{s,x,\alpha^0})]\\
 &\hspace{25mm}\leq
 {\mathcal R_{s,s+\tau,t}}\varphi(x)+\varepsilon.
\end{align*}}
The events
\[
 C_n:=\{(s+\tau,X_{s+\tau}^{s,x,\alpha^0})\in D_n\}
\]
form a countable $\F_\tau$-measurable partition.  Concatenate $\alpha^0$
before the calendar time $s+\tau$ (elapsed time $\tau$) with $\beta^n$
afterwards on $C_n$, and denote the resulting feedback rule by
$\overline\alpha$.  Fix \(P\in\cP_{\cQ}\), and let
\(P^{\tau,\omega}\) be a regular conditional law of the future increment
path after \(\tau\).  By rectangularity,
\[
 P^{\tau,\omega}\in\cP_{\cQ}
 \qquad\text{for \(P\)-almost every \(\omega\)}.
\]
On \(C_n\), the shifted continuation of \(\overline\alpha\) is
\(\beta^n\).  Conditioning at \(\tau\) and using the stopped flow identity
therefore give
\begin{align*}
 &E^P\!\left[
   \varphi(X_t^{s,x,\overline\alpha})
   \,\middle|\,\F_\tau
  \right](\omega)\\
 &\quad=
 E^{P^{\tau,\omega}}\!\left[
  \varphi\!\left(
   X_t^{s+\tau(\omega),
         X_{s+\tau}^{s,x,\alpha^0}(\omega),\beta^n}
  \right)
 \right]\\
 &\quad\leq
 J_{s+\tau(\omega),t}^{\beta^n}\varphi
  (X_{s+\tau}^{s,x,\alpha^0}(\omega))\\
 &\quad\leq
 V(s+\tau(\omega),X_{s+\tau}^{s,x,\alpha^0}(\omega))
 +2\varepsilon.
\end{align*}
Because \((C_n)_{n\geq1}\) is an \(\F_\tau\)-measurable partition, this
conditional inequality holds \(P\)-almost surely on all of \(\Omega\).
Taking expectations and using the tower property yields
\[
 E^P[\varphi(X_t^{s,x,\overline\alpha})]
 \leq
 E^P[V(s+\tau,X_{s+\tau}^{s,x,\alpha^0})]+2\varepsilon.
\]
Since this holds for every \(P\in\cP_{\cQ}\), taking the supremum gives
\begin{align*}
 \bS_{s,t}\varphi(x)
 &\leq J_{s,t}^{\overline\alpha}\varphi(x)\\
 &\leq\sup_{P\in\cP_{\cQ}}
 E^P[V(s+\tau,X_{s+\tau}^{s,x,\alpha^0})]+2\varepsilon\\
 &\leq{\mathcal R_{s,s+\tau,t}}\varphi(x)+3\varepsilon.
\end{align*}
Let \(\varepsilon\downarrow0\) to finish the proof.
\end{proof}

\begin{remark}
\label{rem:BUC}
For the covariance supremum, upper-semi-analytic selection is sufficient.
Although the controls are elementary in time, their stopping grids and
left-endpoint decision variables range over a large class for which no
convenient Polish coding is imposed here.  The continuity argument in the two
DPP proofs avoids using such a topology: one chooses a near-optimal
continuation at each time-state point, enlarges it to an open neighbourhood,
and then uses a countable Borel partition.  Extending the controlled DPP to
arbitrary bounded Borel or upper-semi-analytic terminal data requires either
a standard-Borel coding of elementary or relaxed feedback policies, or an
independent positive-time smoothing argument.  Such an extension should be treated as a separate theorem, but it is
not needed for the results established in this paper.
\end{remark}


\subsection{Feedback-compatible state-path stopping rules}

We now construct the \(\mathbb F\)-stopping-time representatives needed for
state-dependent localisation.  Lemma~\ref{lem:controlled-canonical-state-map}
gives \(\mathbb F\)-progressively measurable canonical state coordinates and
continuous trajectories under each admissible law.  It does not assert
continuity for every canonical path.  Accordingly, we distinguish a stopping
rule on continuous state paths from its law-independent stopping-time
representative on the canonical noise space.  The construction below produces
one representative valid under every admissible law, without completing or
right-continuously augmenting \(\mathbb F\).

To formulate stopping rules in terms of the controlled state, we
introduce a separate canonical state-path space.  Unlike the noise
space in Section~3.1, this space contains \(\R^d\)-valued paths
with arbitrary initial states.  A stopping rule on this space can
then be applied to the state trajectories generated by different
feedback controls.

Let
\[
 \Omega^{\mathrm{state}}:=C([0,\infty);\R^d).
\]
For each \(u\geq0\), define the coordinate projection
\[
 \pi_u:\Omega^{\mathrm{state}}\longrightarrow\R^d,
 \qquad
 \pi_u(\mathbf x):=\mathbf x_u,
\]
and set
\[
 \mathcal G_u^{\mathrm{state}}
 :=\sigma(\pi_r:0\leq r\leq u).
\]
For \(\mathbf x\in\Omega^{\mathrm{state}}\), write
\(\mathbf x^u_r:=\mathbf x_{r\wedge u}\).  If
\(\mathbf y_0=\mathbf x_u\), define
{\begin{equation*}
 (\mathbf x\otimes_u^{\mathrm{state}}\mathbf y)_r
 :=
 \begin{cases}
  \mathbf x_r,&0\leq r\leq u,\\
  \mathbf y_{r-u},&r>u.
 \end{cases}
\end{equation*}}
For \(z_0=(\ell_0,x_0)\) and \(r>0\), write
\begin{equation}\label{eq:canonical-parabolic-cylinder}
 \mathcal C_r(z_0):=(\ell_0-r^2,\ell_0]\times B_r(x_0).
\end{equation}

\begin{definition}
\label{def:feedback-compatible-stopping}
Fix \(0<H<\infty\). A map
\[
 \vartheta:\Omega^{\mathrm{state}}\longrightarrow[0,H]
\]
is a \emph{state-path stopping rule} if
\(\{\vartheta\leq u\}\in\mathcal G_u^{\mathrm{state}}\) for every
\(0\leq u\leq H\).
For fixed \((s,x)\), write
{\begin{equation*}
 \mathbf X^{s,x,\alpha}
 :=\bigl(X_{s+u}^{s,x,\alpha}\bigr)_{u\geq0},
\end{equation*}}
viewed as a continuous-path process under each admissible law.
A family \((\tau^\alpha)_{\alpha\in\cA}\) of
\(\mathbb F\)-stopping times taking values in \([0,H]\) is called
\emph{feedback compatible} if there is a state-path stopping rule
\(\vartheta\) such that
\begin{equation}\label{eq:generated-feedback-stopping-time}
 \tau^\alpha
 =\vartheta\bigl(\mathbf X^{s,x,\alpha}\bigr)
 \quad P\text{-a.s. for every }P\in\cP_{\cQ},
 \qquad \alpha\in\cA.
\end{equation}
On the right, the continuous strong-solution version of
\eqref{eq:controlled-SDE} under \(P\) is used; it is
indistinguishable from the canonical coordinates of
Lemma~\ref{lem:controlled-canonical-state-map}.
The same map
\(\tau^\alpha:\Omega\to[0,H]\) must work for every
\(P\in\cP_{\cQ}\), since it is defined pathwise on the canonical
space and hence independently of the choice of probability law.
Only the almost-sure identity
\(\tau^\alpha
 =\vartheta(\mathbf X^{s,x,\alpha})\)
is interpreted separately under each \(P\).
The parameters \((s,x)\) and those specifying the state-path rule are
held fixed when the family is considered.
\end{definition}

The next lemma proves that such representatives exist and have exactly the
pasting property needed by the stopped DPP.  Any such law-independent
representatives satisfying \eqref{eq:generated-feedback-stopping-time} may be used: all
pasting and flow identities required below are asserted under each
admissible law.

\begin{lemma}
\label{lem:feedback-compatible-exit-rules}
Suppose Assumption~\ref{ass:coeff} holds and let \(\vartheta\) be a
state-path stopping rule bounded by \(H\).  Then:
\begin{enumerate}[label=\textup{(\roman*)}]
\item For every fixed \((s,x)\) and every \(\alpha\in\cA\), there is a
single \(\mathbb F\)-stopping time \(\tau^\alpha\) satisfying
\eqref{eq:generated-feedback-stopping-time}.  The construction is explicit
and does not depend on an admissible law.
\item Suppose \(C_n\in\F_{\tau^\alpha}\) form a countable partition and
\(\beta^n\in\cA\).  If \(\overline\alpha\) is obtained by using
\(\alpha\) up to \(\tau^\alpha\) and \(\beta^n\) afterwards on
\(C_n\), as in \eqref{eq:control-pasting}, then
{\begin{equation*}
 \tau^{\overline\alpha}=\tau^\alpha
 \quad P\text{-a.s. for every }P\in\cP_{\cQ}.
\end{equation*}}
Under each such law the two controlled state paths agree up to this common
time, up to indistinguishability.
\item If \(u\leq\vartheta(\mathbf x)\) and
\(\mathbf y_0=\mathbf x_u\), then
{\begin{equation*}
 \vartheta^{u,\mathbf x}(\mathbf y)
 :=\vartheta(\mathbf x\otimes_u^{\mathrm{state}}\mathbf y)-u
\end{equation*}}
is a nonnegative state-path stopping rule for the future state path, bounded by
\(H-u\).  In particular,
\[
 \vartheta(\mathbf x\otimes_u^{\mathrm{state}}\mathbf y)
 =u+\vartheta^{u,\mathbf x}(\mathbf y).
\]
\end{enumerate}
\end{lemma}

\begin{proof}
\emph{Step 1: an identity for stopped paths and the available path regularity.}
For a path \(\mathbf x\in\Omega^{\mathrm{state}}\), recall that
\[
 \mathbf x^{\vartheta(\mathbf x)}_r
 :=
 \mathbf x_{r\wedge\vartheta(\mathbf x)}.
\]
Thus, if two paths agree up to the stopping time
\(\vartheta(\mathbf x)\), then they have the same stopping time.
Equivalently, the stopping decision is determined by the stopped path. To see this, we need to show that 
a state-path stopping rule has the standard non-anticipative
property
\begin{equation}\label{eq:state-path-stopping-property}
 \mathbf x^{\vartheta(\mathbf x)}
 =\mathbf y^{\vartheta(\mathbf x)}
 \quad\Longrightarrow\quad
 \vartheta(\mathbf y)=\vartheta(\mathbf x).
\end{equation}
Indeed, set \(v:=\vartheta(\mathbf x)\).  The hypothesis means that
\(\mathbf x\) and \(\mathbf y\) agree on \([0,v]\).
First let \(q<v\) be rational.  Since
\(\vartheta(\mathbf x)=v>q\),
\[
 \mathbf x\notin\{\vartheta\leq q\}.
\]
The event \(\{\vartheta\leq q\}\) belongs to
\(\mathcal G_q^{\mathrm{state}}\), and therefore its membership
depends only on the restriction of a path to \([0,q]\).  Since
\(\mathbf x\) and \(\mathbf y\) agree on this interval,
\[
 \mathbf y\notin\{\vartheta\leq q\}.
\]
Thus \(\vartheta(\mathbf y)>q\) for every rational \(q<v\), and hence
\[
 \vartheta(\mathbf y)\geq v.
\]
On the other hand,
\(\mathbf x\in\{\vartheta\leq v\}\), and
\(\{\vartheta\leq v\}\in\mathcal G_v^{\mathrm{state}}\).
Since the two paths agree on \([0,v]\), it follows that
\(\mathbf y\in\{\vartheta\leq v\}\).  Hence
\[
 \vartheta(\mathbf y)\leq v.
\]
Consequently \(\vartheta(\mathbf y)=v=\vartheta(\mathbf x)\).

Fix \((s,x,\alpha)\), and abbreviate the canonical state coordinates by
\(Y_u:=X_{s+u}^{s,x,\alpha}\), so that \(Y_0=x\) on all of \(\Omega\).
The moment estimate in Lemma~\ref{lem:controlled-flow}, the controlled
SDE and the Burkholder--Davis--Gundy inequality give
\begin{equation}\label{eq:canonical-lift-increment-moment}
 \sup_{P\in\cP_{\cQ}}E^P[|Y_v-Y_u|^8]
 \leq C_{s,H}(1+|x|^8)|v-u|^4,
 \qquad 0\leq u\leq v\leq H,
\end{equation}
with a constant that may also be chosen uniformly in \(\alpha\).
For the drift increment, H\"older's inequality gives a bound of order
\(|v-u|^8\); for the martingale increment, BDG gives a bound of order
\(|v-u|^4\).  The uniform eighth-moment bound on \([0,H]\) controls their
coefficients, proving \eqref{eq:canonical-lift-increment-moment}.
Thus, under every admissible \(P\), the continuous state trajectory is
\(\kappa\)-H\"older on \([0,H]\) for \(\kappa:=1/4\), almost surely,
by Kolmogorov's continuity theorem.  Its H\"older version and its continuous
strong-solution version are indistinguishable, since they agree at all
rational times on a common \(P\)-full set.  This concerns the sample paths
of the SDE and does not use the positive-time regularity of the value
function proved later.  We use it only under each law, not to modify the
canonical coordinates on exceptional paths.

\smallskip
\noindent\emph{Step 2: continuous reconstructions from the observed prefixes.}
For \(n\geq1\), let
\begin{equation}\label{eq:canonical-lift-prefix-time}
 \begin{split}
 \zeta_n:=H\wedge\inf\bigl\{q\in\mathbb Q\cap(0,H):
 &\ \text{there is }p\in\mathbb Q\cap[0,q)\text{ with}\\
 &|Y_q-Y_p|>n(q-p)^\kappa\bigr\},
 \end{split}
\end{equation}
where the infimum of the empty set is \(+\infty\).
The variables \(\zeta_n\) increase with \(n\).  For \(0<t\leq H\),
\[
 \{\zeta_n<t\}
 =\bigcup_{\substack{p,q\in\mathbb Q\cap[0,H)\\0\leq p<q<t}}
       \{|Y_q-Y_p|>n(q-p)^\kappa\}
 \in\F_t.
\]
We do \emph{not} claim that \(\zeta_n\) is an \(\mathbb F\)-stopping time:
only the strict-event property just displayed is needed.
It implies that \(t\wedge\zeta_n\) is \(\F_t\)-measurable for each
\(t\).  Indeed, for \(v\leq t\) its strict sublevel event is
\(\{t\wedge\zeta_n<v\}=\{\zeta_n<v\}\in\F_v\subset\F_t\), whereas
for \(v>t\) that event is the entire space.

For each \(\omega\), the rational coordinates
\(Y_q(\omega)\), \(0\leq q<\zeta_n(\omega)\), satisfy the
\(\kappa\)-H\"older bound with constant \(n\).  If \(\zeta_n>0\),
let \(Y^n\) be their unique continuous extension to \([0,\zeta_n]\),
frozen after \(\zeta_n\).  If \(\zeta_n=0\), set \(Y^n\equiv x\).
This defines a continuous path on \([0,\infty)\) for every \(\omega\).
The reconstruction is \(\mathbb F\)-adapted.  To see this explicitly, put
\[
 d_j(v):=2^{-j}\max\{\lceil2^jv\rceil-1,0\},\qquad v\geq0.
\]
Then \(d_j(v)<v\) for \(v>0\), \(d_j(v)\to v\), and
{\begin{equation*}
 Y_t^n=\lim_{j\to\infty}Y_{d_j(t\wedge\zeta_n)},
 \qquad t\geq0.
\end{equation*}}
The limit exists and is finite for every path by the H\"older bound.
For fixed \(t\), each variable before the limit is \(\F_t\)-measurable:
its time argument is \(\F_t\)-measurable, takes countably many
nonnegative rational values not exceeding \(t\), and \(Y_q\) is
\(\F_q\)-measurable.  Hence \(Y^n\) is \(\mathbb F\)-adapted; its all-path
continuity also makes it \(\mathbb F\)-progressively measurable and a Borel map into
\(\Omega^{\mathrm{state}}\).

These reconstructions agree on their common prefixes:
\begin{equation}\label{eq:canonical-lift-prefix-consistency}
 Y^m_t=Y^n_t\quad\text{for }0\leq t\leq\zeta_n,
 \qquad m\geq n,
\end{equation}
pathwise.  This follows first at rational times strictly below
\(\zeta_n\), where both equal \(Y\), and then at the endpoint by
continuity.  It is also immediate if \(\zeta_n=0\).

\smallskip
\noindent\emph{Step 3: a single law-independent stopping-time representative}.
Since \(Y^n\) is an all-path continuous \(\mathbb F\)-adapted process,
\(\sigma_n:=\vartheta(Y^n)\) is an \(\mathbb F\)-stopping time.
Accept its stopping decision only if it occurs within the reconstructed
prefix, by setting
{\begin{equation*}
 T_n:=
 \begin{cases}
  \sigma_n,&\sigma_n\leq\zeta_n,\\
  H,&\sigma_n>\zeta_n.
 \end{cases}
\end{equation*}}
This is again an \(\mathbb F\)-stopping time.  In fact,
\(\{\zeta_n<\sigma_n\}\in\F_{\sigma_n}\), because for every \(t\geq0\),
\[
 \{\zeta_n<\sigma_n\}\cap\{\sigma_n\leq t\}
 =\bigcup_{q\in\mathbb Q\cap(0,t)}
    \{\zeta_n<q\}\cap\{q<\sigma_n\leq t\}\in\F_t.
\]
Therefore the event on which \(\sigma_n\) is retained is known at
\(\sigma_n\); explicitly, for \(t<H\),
\(\{T_n\leq t\}=\{\sigma_n\leq t\}\cap\{\sigma_n\leq\zeta_n\}\in\F_t\).

If \(T_n(\omega)<H\), then \(\sigma_n\leq\zeta_n\leq\zeta_m\)
for \(m\geq n\), and \eqref{eq:canonical-lift-prefix-consistency} together
with \eqref{eq:state-path-stopping-property} implies
\[
 \sigma_m(\omega)=\sigma_n(\omega),\qquad
 T_m(\omega)=T_n(\omega),\qquad m\geq n.
\]
Consequently, for every path the sequence \((T_n)\) is either identically
\(H\), or becomes constant as soon as a value below \(H\) occurs.
Its infimum is therefore attained on every path.  Define
{\begin{equation*}
 \tau^\alpha:=\inf_{n\geq1}T_n.
\end{equation*}}
For every \(t\geq0\), attainment gives the exact identity
{\begin{equation*}
 \{\tau^\alpha\leq t\}
 =\bigcup_{n\geq1}\{T_n\leq t\}\in\F_t.
\end{equation*}}
Thus \(\tau^\alpha\) is an \(\mathbb F\)-stopping time.  Attainment is essential here: an arbitrary
decreasing limit of \(\mathbb F\)-stopping times would not justify this conclusion
without right-continuity of the filtration.

Fix now \(P\in\cP_{\cQ}\).  On a \(P\)-full set, the state path is
continuous and has a finite \(\kappa\)-H\"older constant on \([0,H]\).
For each such path and every sufficiently large integer \(n\), there is
no violation in \eqref{eq:canonical-lift-prefix-time}, so \(\zeta_n=H\) and
\(Y^n=Y\) on \([0,H]\).  Hence
\(T_n=\vartheta(\mathbf X^{s,x,\alpha})\).
For any smaller \(n\), either \(T_n=H\), or
\(\sigma_n\leq\zeta_n\) and \(Y^n\) agrees with the continuous state
path up to \(\sigma_n\).  In the latter case
\eqref{eq:state-path-stopping-property} again gives
\(T_n=\vartheta(\mathbf X^{s,x,\alpha})\).
Taking the infimum proves \eqref{eq:generated-feedback-stopping-time}.
The construction used only the canonical coordinates and \(\vartheta\), not
\(P\), so it is the same representative under every admissible law.
This proves (i).  For a fixed state-path rule it is also jointly Borel in
\((s,x,\omega)\), since all its operations are Borel and use countably
many state coordinates and limits.

\smallskip
\noindent\emph{Step 4: pasting invariance and residual rules.}
Let \(\overline\alpha\) be as in (ii), and fix an admissible \(P\).
Both \(\tau^\alpha\) and \(\tau^{\overline\alpha}\) are already
\(\mathbb F\)-stopping times.  The two controls agree up to \(\tau^\alpha\), so
pathwise uniqueness for the controlled SDE gives
\[
 X_{s+u}^{s,x,\overline\alpha}=X_{s+u}^{s,x,\alpha},
 \qquad 0\leq u\leq\tau^\alpha,\quad P\text{-a.s.}
\]
Intersect the \(P\)-full sets on which both state paths are continuous
and both identities in \eqref{eq:generated-feedback-stopping-time} hold.
On that intersection, \eqref{eq:state-path-stopping-property} gives
\[
 \tau^{\overline\alpha}
 =\vartheta(\mathbf X^{s,x,\overline\alpha})
 =\vartheta(\mathbf X^{s,x,\alpha})
 =\tau^\alpha.
\]
This proves (ii), also for any other choice of law-independent
representatives satisfying \eqref{eq:generated-feedback-stopping-time}.
The pasted control remains in the original class \(\cA\), since
the pasting time is an \(\mathbb F\)-stopping time and the
partition is \(\F_{\tau^\alpha}\)-measurable.

For (iii), the concatenated path agrees with \(\mathbf x\) up to \(u\).
The stopping-rule property prevents it from stopping before \(u\)
when \(u\leq\vartheta(\mathbf x)\).  Also, for \(0\leq v\leq H-u\),
\[
 \{\vartheta^{u,\mathbf x}\leq v\}
 =\{\mathbf y:
       \vartheta(\mathbf x\otimes_u^{\mathrm{state}}\mathbf y)
          \leq u+v\}
\]
belongs to \(\mathcal G_v^{\mathrm{state}}\), the canonical
state-path sigma-field up to time \(v\).
The stated shift identity follows from the definition.  The case \(u=H\)
has the identically zero residual rule.  Otherwise, applying (i)
to the residual rule gives a law-independent
\(\mathbb F\)-stopping-time representative for the restarted problem;
together with the stopped flow of
Lemma~\ref{lem:controlled-canonical-state-map}, these representatives
have the same residual-time interpretation under each admissible law.

\end{proof}

Let \(z_0=(\ell_0,x_0)\), \(r>0\), and
\((\ell,x)\in\mathcal C_r(z_0)\), with \(\ell>0\).  With
\(\inf\emptyset:=\ell\), define on continuous state paths
{\begin{equation*}
 \vartheta_{r,z_0}^{\ell}(\mathbf x)
 :=\inf\left\{0\leq u\leq\ell:
    \ell-u\leq\ell_0-r^2
    \ \text{or}\ |\mathbf x_u-x_0|\geq r\right\}.
\end{equation*}}
For a compact \(\Gamma\subset(0,\ell]\times\R^d\), define
{\begin{equation*}
 \sigma_\Gamma^\ell(\mathbf x)
 :=\inf\left\{0\leq u\leq\ell:
       (\ell-u,\mathbf x_u)\in\Gamma\right\}.
\end{equation*}}
{
On \(\Omega^{\mathrm{state}}\), the path
\(u\mapsto(\ell-u,\mathbf x_u)\) is continuous.  Each rule above is
therefore a capped first hitting time of a closed set, and the event that
it occurs by time \(u\) is Borel and determined by the state path up to
\(u\).  Hence both are state-path stopping rules, and so is their minimum.
By Lemma~\ref{lem:feedback-compatible-exit-rules}(i), they admit
law-independent \(\mathbb F\)-stopping-time representatives on the canonical
noise space; the minimum of two such representatives represents the minimum
of the state-path rules under every admissible law.  Part~(ii) gives pasting
invariance.  Thus these representatives generate feedback-compatible
families, including the compact-set hitting-before-exit rules used in the
Krylov--Safonov argument.
}

All subsequent controlled hitting and exit times are understood as
law-independent \(\mathbb F\)-stopping-time representatives supplied by this
lemma.  Their usual first-hit formulae hold under every admissible law on
the continuous state trajectories.  They are therefore legitimate times for
the canonical conditioning and pasting
arguments, as well as for the lawwise hitting estimates and stopped
It\^o formulae.  No all-path continuous modification of the canonical
state map has been assumed.

\section{Positive-time regularity under uniform parabolicity}

\subsection{Uniformly parabolic controlled framework}

For $a\in A$ and $Q\in\cQ$, define the effective drift and covariance
\begin{align}
 b_{a,Q}(q,x)
 &:={}
 b(q,x,a)+\sum_{i,j=1}^mQ_{ij}h_{ij}(q,x,a),
 \label{eq:controlled-effective-drift}\\
 {c_{a,Q}}(q,x)
 &:={}
 \sigma(q,x,a)Q\sigma(q,x,a)^{\mathsf T}.
 \label{eq:controlled-effective-covariance}
\end{align}
For a smooth function \(\phi\), write
\begin{equation}\label{eq:controlled-generator}
 \mathcal L_{a,Q}^{q}\phi(x)
 :=b_{a,Q}(q,x)\cdot D\phi(x)
   +\frac12\tr\!\left({c_{a,Q}}(q,x)D^2\phi(x)\right).
\end{equation}

{\renewcommand{\theassumption}{C}
\begin{assumption}
\label{ass:uniform-parabolic}
The coefficients are jointly continuous in $(q,x,a)$.  The Lipschitz and
linear-growth constants in Assumption~\ref{ass:coeff} can be chosen
independently of calendar time, and the time modulus can be chosen uniformly
in $(x,a)$ after division by $1+\abs{x}$.  There are constants
$0<\underline\lambda\leq\overline\lambda<\infty$ such that
\begin{equation}\label{eq:controlled-effective-ellipticity}
 \underline\lambda I_d
 \leq {c_{a,Q}}(q,x)
 \leq\overline\lambda I_d
\end{equation}
for every $(q,x,a,Q)\in\R\times\R^d\times A\times\cQ$.  Consequently, for
every bounded $O\subset\R^d$,
\begin{equation}\label{eq:controlled-local-drift-bound}
 B_O
 :=\sup_{\substack{q\in\R,\,x\in O\\a\in A,\,Q\in\cQ}}
 \abs{b_{a,Q}(q,x)}<\infty.
\end{equation}
\end{assumption}
}

The bounds are deliberately uniform in the action.  Without this
uniformity, taking the infimum over all feedback rules can destroy both the
regularity estimate and the Harnack constant.

\subsection{Positive-time regularity}

We first record the stopping-rule form of the DPP needed for exit-time
localisation.  Feedback compatibility is meant in the precise sense of
Definition~\ref{def:feedback-compatible-stopping}.  In particular,
Lemma~\ref{lem:feedback-compatible-exit-rules} proves both
the \(\mathbb F\)-stopping-time property and pasting invariance for the
cylinder exit and closed-set hitting rules
used below.

\begin{proposition}
\label{prop:compatible-stopped-DPP}
{Under Assumption~\ref{ass:coeff}, let \(s<t\), let
\(\varphi\in\BUC_b(\R^d)\), and let}
$\alpha\mapsto\tau^\alpha$ be feedback compatible with
$0\leq\tau^\alpha\leq t-s$.  Then
\begin{equation}\label{eq:compatible-stopped-DPP}
 \bS_{s,t}\varphi(x)
 =\inf_{\alpha\in\cA}\sup_{P\in\cP_{\cQ}}
 E^P\!\left[
  \bS_{s+\tau^\alpha,t}\varphi
  (X_{s+\tau^\alpha}^{s,x,\alpha})
 \right].
\end{equation}
\end{proposition}

\begin{proof}
Put \(V(a,y):=\bS_{a,t}\varphi(y)\).  This is bounded and continuous
by Proposition~\ref{prop:controlled-continuity}.
For each fixed feedback rule \(\alpha\),
Lemma~\ref{lem:feedback-compatible-exit-rules} ensures that
\(\tau^\alpha\) is an \(\mathbb F\)-stopping time.
The conditioning and analytic-selection argument used in the proof
of Theorem~\ref{thm:controlled-stopped-DPP} can therefore be applied
with \(\tau=\tau^\alpha\).
For each fixed feedback rule, this is a single law-independent
\(\mathbb F\)-stopping-time representative.  Its equality with the state-path
rule in \eqref{eq:generated-feedback-stopping-time} is used only under
an admissible law, including any law obtained by pasting.

For the lower bound, fix a full-horizon rule \(\gamma\).  Conditional on
the history at \(\tau^\gamma\), the shifted continuation
\(\gamma^{\tau^\gamma,\omega}\) is admissible.  The continuation-cost map
is Borel by Lemma~\ref{lem:controlled-canonical-state-map}, and analytic
selection supplies an almost worst-case continuation law.  Pasting this law
at \(\tau^\gamma\) gives
\[
 J_{s,t}^{\gamma}\varphi(x)
 \geq\sup_{P\in\cP_{\cQ}}
 E^P\!\left[
  V(s+\tau^\gamma,
    X_{s+\tau^\gamma}^{s,x,\gamma})
 \right]-\varepsilon.
\]
Taking the infimum over \(\gamma\) proves the lower inequality in
\eqref{eq:compatible-stopped-DPP}.

For the reverse inequality, for each
\((a,y)\in[s,t]\times\mathbb R^d\), choose
\(\beta^{a,y}\in\mathcal A\) with
\begin{equation}\label{eq:compatible-continuation-near-optimal}
J_{a,t}^{\beta^{a,y}}\varphi(y)
 \leq V(a,y)+\varepsilon.
\end{equation}
Continuity of the fixed-policy cost and of \(V\) enlarges this inequality,
with error \(2\varepsilon\), to an open neighbourhood of \((a,y)\).
Since \([s,t]\times\R^d\) is second countable, we can choose
countably many of these neighbourhoods that cover it.
Write \(U_n\) for their intersections with \([s,t]\times\R^d\),
and let \(\beta^n\) be the continuation rule associated with \(U_n\).
Define
\[
 D_1:=U_1,
 \qquad
 D_n:=U_n\setminus\bigcup_{k<n}U_k,
 \quad n\geq2.
\]
After discarding empty sets, the sets \(D_n\) form a countable
Borel partition of \([s,t]\times\R^d\).
Since \(D_n\subset U_n\), the estimate for \(\beta^n\) remains
valid at every point of \(D_n\).

Having selected the continuation rules \(\beta^n\), we now choose
a nearly optimal rule for the stopped first-stage problem.
By the definition of the infimum, there exists \(\alpha^0\in\cA\)
such that
\begin{align*}
 &\sup_{P\in\cP_{\cQ}}
 E^P\!\left[
   V\!\left(s+\tau^{\alpha^0},
     X_{s+\tau^{\alpha^0}}^{s,x,\alpha^0}\right)
 \right]\\
 &\qquad\leq
 \inf_{\alpha\in\cA}\sup_{P\in\cP_{\cQ}}
 E^P\!\left[
   V\!\left(s+\tau^\alpha,
     X_{s+\tau^\alpha}^{s,x,\alpha}\right)
 \right]+\varepsilon.
\end{align*}
Put
\[
 C_n:=\left\{\omega\in\Omega:
  \left(
   s+\tau^{\alpha^0}(\omega),
   X_{s+\tau^{\alpha^0}(\omega)}^{s,x,\alpha^0}(\omega)
  \right)\in D_n
 \right\}.
\]
The stopped time--state pair takes values in
\([s,t]\times\R^d\).  Since the sets \(D_n\) form a partition
of this space, the events \(C_n\) are pairwise disjoint and satisfy
\[
 \bigcup_{n\geq1}C_n=\Omega.
\]
Moreover, each \(C_n\) belongs to
\(\F_{\tau^{\alpha^0}}\), because the stopped time--state pair
is measurable with respect to this sigma-field.
Some \(C_n\) may be empty; after discarding these events and their
associated continuation rules, we obtain a countable, possibly
finite, \(\F_{\tau^{\alpha^0}}\)-measurable partition of \(\Omega\).
These events form an \(\mathcal F_{\tau^{\alpha^0}}\)-measurable partition.
Paste \(\beta^n\) after \(\tau^{\alpha^0}\) on \(C_n\), and denote the
result by \(\overline\alpha\).  By
Lemma~\ref{lem:feedback-compatible-exit-rules}, under every
\(P\in\mathcal P_{\mathcal Q}\), almost surely,
\[
 \tau^{\overline\alpha}=\tau^{\alpha^0},
 \qquad
 X_{s+u}^{s,x,\overline\alpha}
 =X_{s+u}^{s,x,\alpha^0},
 \quad 0\leq u\leq\tau^{\alpha^0}.
\]
Conditioning at this common stopping time and using the selected
continuation on each \(C_n\) gives
\begin{align*}
 J_{s,t}^{\overline\alpha}\varphi(x)
 &\leq
 \sup_{P\in\cP_{\cQ}}E^P\!\left[
  V(s+\tau^{\alpha^0},
    X_{s+\tau^{\alpha^0}}^{s,x,\alpha^0})
 \right]+2\varepsilon.
\end{align*}
The choice of \(\alpha^0\) bounds the last expression by the right-hand
side of \eqref{eq:compatible-stopped-DPP} plus \(3\varepsilon\).  Letting
\(\varepsilon\downarrow0\) proves the result.
\end{proof}

Recall \(\mathcal C_r(z_0)\) from
\eqref{eq:canonical-parabolic-cylinder}, and write
\[
 d_{\rm p}((\ell,x),(k,y))
 :=\abs{x-y}+\abs{\ell-k}^{1/2}.
\]

We first specify the local setting for the following hitting estimate.
Fix \(\zeta=(k_0,y_0)\), \(\rho>0\), and set
\[
 A_\rho(\zeta)
 :=(k_0-6\rho^2,k_0-5\rho^2)\times B_\rho(y_0).
\]
Assume that
\(\overline{\mathcal C_{3\rho}(\zeta)}
  \subset(0,\infty)\times\mathbb R^d\).
Suppose that, for every
\(z=(\ell,x)\in\mathcal C_\rho(\zeta)\), every admissible feedback
rule \(\alpha\), and every admissible law \(P\), the following
decomposition holds up to the first exit from
\(\mathcal C_{3\rho}(\zeta)\):
\[
 \mathsf X_q^{\ell,x,\alpha}
 =x+\int_0^q\beta_u^{\alpha,P}\,du+M_q^{\alpha,P},
 \qquad
 \abs{\beta_u^{\alpha,P}}\leq B,
\]
where \(M^{\alpha,P}\) is a continuous local \(P\)-martingale satisfying
\[
 \lambda I_d
 \leq\frac{d\langle M^{\alpha,P}\rangle_u}{du}
 \leq\Lambda I_d.
\]
Write
\[
 Z_q^{z,\alpha}
 :=(\ell-q,\mathsf X_q^{\ell,x,\alpha}).
\]
For a closed set \(\Gamma\subset A_\rho(\zeta)\), let
\(\sigma_\Gamma^{z,\alpha}\) and \(\tau_{3\rho}^{z,\alpha}\) denote the
law-independent \(\mathbb F\)-stopping-time representatives of the
corresponding state-path hitting and exit rules from
Lemma~\ref{lem:feedback-compatible-exit-rules}.  Under every admissible
law, they have the first-hit descriptions
\begin{align}
 \sigma_\Gamma^{z,\alpha}
 &=\inf\left\{0\leq q\leq\ell:
       Z_q^{z,\alpha}\in\Gamma\right\},
 \label{eq:controlled-Gamma-hitting-time}\\
 \tau_{3\rho}^{z,\alpha}
 &=\inf\left\{0\leq q\leq\ell:
       Z_q^{z,\alpha}\notin\mathcal C_{3\rho}(\zeta)\right\},
 \label{eq:controlled-cylinder-exit-time}
\end{align}
\(P\)-almost surely, with \(\inf\emptyset:=\ell\).

\begin{lemma}
\label{lem:controlled-uniform-hitting}
Under the preceding assumptions and notation, there exists
\(p_0=p_0(d,\lambda,\Lambda,\rho B)\in(0,1)\) such that, for every
closed set \(\Gamma\subset A_\rho(\zeta)\) satisfying
\[
 \abs{\Gamma}\geq\frac14\abs{A_\rho(\zeta)},
\]
one has
\begin{equation}\label{eq:controlled-uniform-hitting}
 \inf_{\substack{z\in\mathcal C_\rho(\zeta)\\
                 \alpha\in\cA,\ P\in\cP_{\cQ}}}
 P\!\left\{
   \sigma_\Gamma^{z,\alpha}<\tau_{3\rho}^{z,\alpha}
 \right\}
 \geq p_0.
\end{equation}
If \(0<\rho\leq r\), then \(p_0\) may be chosen to depend on
the drift bound only through \(rB\).
\end{lemma}

\begin{proof}
Fix \(z=(\ell,x)\), \(\alpha\), and \(P\), and make the parabolic
change of variables
\[
 \widehat Z_u
 :=\left(
    \frac{\ell-k_0}{\rho^2}-u,
    \frac{\mathsf X_{\rho^2u}^{\ell,x,\alpha}-y_0}{\rho}
   \right).
\]
This sends the initial cylinder \(\mathcal C_\rho(\zeta)\), the exit
cylinder \(\mathcal C_{3\rho}(\zeta)\), and the target slab
\(A_\rho(\zeta)\) respectively to
\[
 (-1,0]\times B_1,
 \qquad (-9,0]\times B_3,
 \qquad (-6,-5)\times B_1.
\]
The image
\[
 \widehat\Gamma
 :=\left\{
  \left(\frac{k-k_0}{\rho^2},\frac{y-y_0}{\rho}\right):
  (k,y)\in\Gamma
 \right\}
\]
satisfies
\[
 \abs{\widehat\Gamma}
 =\rho^{-(d+2)}\abs{\Gamma}
 \geq\frac14\abs{(-6,-5)\times B_1}.
\]
The spatial component of \(\widehat Z\) has drift bounded by \(\rho B\).
Its martingale part
\[
 \widehat M_u:=\rho^{-1}M_{\rho^2u}^{\alpha,P}
\]
satisfies
\[
 \lambda I_d
 \leq\frac{d\langle\widehat M\rangle_u}{du}
 \leq\Lambda I_d.
\]
Writing \(\widehat a_u=d\langle\widehat M\rangle_u/du\),
L\'evy's characterisation represents the martingale part, up to the
localisation time, as
\(d\widehat M_u=\widehat a_u^{1/2}\,dW_u\).
For a precise adapted-coefficient version of the classical hitting
estimates \cite{KrylovSafonov1979,Krylov1987}, we use
\cite[Corollary~3.6, p.~14, of the arXiv version]{Krylov2023LdPlusOne}.
The processes in Sections~2--3 of that reference have adapted drift
and diffusion coefficients; no Markov assumption is imposed there.

To match its forward-time cylinders with our remaining-time coordinates,
put
\[
 \theta:=\frac{\ell-k_0}{\rho^2}\in(-1,0],
 \qquad
 \widehat Y_u:=\frac{\mathsf X_{\rho^2u}^{\ell,x,\alpha}-y_0}{\rho},
\]
and define
\[
 \Gamma_\theta
 :=\{(u,y):(\theta-u,y)\in\widehat\Gamma\}.
\]
Then
\[
 \Gamma_\theta\subset(4,6)\times B_1
 \subset[3,6)\times B_{\sqrt3},
 \qquad
 \frac{|\Gamma_\theta|}{|[3,6)\times B_{\sqrt3}|}
 \geq\frac{1}{4\cdot 3^{(d+2)/2}}.
\]
For \(d\geq2\), apply the cited corollary with radius \(\sqrt3\)
and interior parameter \(3/4\).  The initial state lies in \(B_1\),
hence in \(B_{(3/4)\sqrt3}\).  Since
\(\abs{\beta_u^{\alpha,P}}\leq B\) before rescaling, the rescaled
drift is bounded by \(\rho B\).

We verify the drift hypotheses in that reference by localisation.
After exit from \([0,6)\times B_{\sqrt3}\), set the drift to zero
and continue the martingale part with a fixed nondegenerate covariance.
The stopped drift, expressed along the increment process
\(\widehat Y_u-\widehat Y_0\), is dominated by
\[
 g(u,z):=\rho B\,\mathbf 1_{[0,6]\times B_4}(u,z).
\]
Taking both mixed-norm exponents in the reference to be \(d+1\),
one may choose its local drift-norm bound as
\[
 \overline b_r=C_d(\rho B)^{d+1}\min\{r^{d+1},1\},
 \qquad r>0.
\]
Indeed, for \(r\leq1\) the integral of \(|g|^{d+1}\) over a
parabolic cylinder of radius \(r\) is at most
\(C_d(\rho B)^{d+1}r^{d+2}\), and for larger \(r\) the total
integral is bounded by compact support.
Thus Assumptions~2.1--2.2 of the reference hold with a sufficiently
small radius depending only on \(d,\lambda,\Lambda,\rho B\).
Corollary~3.6 therefore gives a positive lower bound, depending only
on these parameters, for hitting \(\Gamma_\theta\) before exit from
\([0,6)\times B_{\sqrt3}\).

This hitting event occurs before exit from our larger cylinder:
\(B_{\sqrt3}\subset B_3\), and its remaining-time bottom is reached
at elapsed time \(\theta+9>8\).  It therefore implies the event in
\eqref{eq:controlled-uniform-hitting}.
For \(d=1\), apply the same two-dimensional estimate after adjoining
an independent Brownian coordinate and replacing the target by
\[
 \{(u,y,z):(u,y)\in\Gamma_\theta,\ |z|\leq1\}.
\]
This target is contained in \([3,6)\times B_{\sqrt3}\) in two spatial
dimensions and has a uniformly positive volume fraction.  Its hitting
implies hitting the original target.  This auxiliary enlargement is
used only for the probability estimate and does not change the original
admissible-law family.
The resulting bounds are uniform in the initial point, the feedback
rule and the admissible law.  Taking their infimum proves
\eqref{eq:controlled-uniform-hitting}.  When \(\rho\leq r\), use the
common drift bound \(rB\) to obtain the final assertion of the lemma.
\end{proof}

\begin{lemma}
\label{lem:controlled-KS-smoothing}
Let \(D\subset(0,\infty)\times\R^d\) be open and let
\(v:D\to\R\) be bounded and measurable.
Fix \(z_0=(\ell_0,x_0)\) and \(r>0\) such that
\[
 \overline{\mathcal C_{4r}(z_0)}\subset D.
\]
For \(z=(\ell,x)\in D\), write
\[
 Z_q^{z,\alpha}
 :=(\ell-q,\mathsf X_q^{\ell,x,\alpha}).
\]
Suppose that, on every parabolic cylinder
\(\mathcal C_R(\zeta)\) with
\(\overline{\mathcal C_R(\zeta)}\subset D\),
the local characteristic assumptions stated before
Lemma~\ref{lem:controlled-uniform-hitting} hold uniformly in
\(\alpha\) and \(P\), with constants \(B,\lambda,\Lambda\).
Assume also that, for every such cylinder, every starting point
\(z=(\ell,x)\in\mathcal C_R(\zeta)\), and every
feedback-compatible stopping rule \(\theta^\alpha\) bounded by
the exit time of \(Z^{z,\alpha}\) from \(\mathcal C_R(\zeta)\),
\begin{equation}\label{eq:abstract-controlled-stopped-DPP}
 v(\ell,x)
=\inf_{\alpha\in\cA}\sup_{P\in\cP_{\cQ}}
 E^P\!\left[
  v(\ell-\theta^\alpha,
    \mathsf X_{\theta^\alpha}^{\ell,x,\alpha})
 \right].
\end{equation}
Then there are \(\eta\in(0,1)\) and \(C<\infty\), depending only
on \(d,\lambda,\Lambda,rB\) and the fixed cylinder geometry,
such that
\begin{equation}\label{eq:controlled-parabolic-Holder}
 \abs{v(z)-v(z')}
 \leq C\,\osc_{\mathcal C_{4r}(z_0)}v
 \left(\frac{d_{\rm p}(z,z')}{r}\right)^\eta,
 \qquad z,z'\in\mathcal C_r(z_0).
\end{equation}
\end{lemma}

\begin{proof}
We divide the proof into a one-step oscillation reduction and its
iteration.  Fix a cylinder
\(\overline{\mathcal C_{4\rho}(\zeta)}
 \subset \mathcal C_{4r}(z_0)\), where
\(\zeta=(k_0,y_0)\), and use the slab \(A_\rho(\zeta)\) introduced before
Lemma~\ref{lem:controlled-uniform-hitting}.  Set
\[
 M:=\sup_{\mathcal C_{4\rho}(\zeta)}v,
 \qquad
 m:=\inf_{\mathcal C_{4\rho}(\zeta)}v.
\]
If \(M=m\), there is nothing to prove, so assume \(M>m\).

Choose a median \(\mu\in[m,M]\) of the measurable function \(v\) on
\(A_\rho(\zeta)\); thus
\[
 \abs{A_\rho(\zeta)\cap\{v\leq\mu\}}
 \geq \frac12\abs{A_\rho(\zeta)},
 \qquad
 \abs{A_\rho(\zeta)\cap\{v\geq\mu\}}
 \geq \frac12\abs{A_\rho(\zeta)}.
\]
By inner regularity of Lebesgue measure, there are compact sets
\[
 \Gamma_-\subset A_\rho(\zeta)\cap\{v\leq\mu\},
 \qquad
 \Gamma_+\subset A_\rho(\zeta)\cap\{v\geq\mu\},
\]
such that
\[
 \abs{\Gamma_-},\abs{\Gamma_+}
 \geq \frac14\abs{A_\rho(\zeta)}.
\]
Since \(\rho\leq r\), Lemma~\ref{lem:controlled-uniform-hitting}
provides a constant
\[
 p_0=p_0(d,\lambda,\Lambda,rB)\in(0,1)
\]
which works simultaneously for \(\Gamma_-\) and \(\Gamma_+\), uniformly
in the starting point, the feedback rule and the admissible law.

Fix \(z=(\ell,x)\in\mathcal C_\rho(\zeta)\).  For
\(\varepsilon\in\{-,+\}\), let
\[
 \theta_\varepsilon^{z,\alpha}
 :=\sigma_{\Gamma_\varepsilon}^{z,\alpha}
   \wedge\tau_{3\rho}^{z,\alpha},
 \qquad
 H_\varepsilon^{z,\alpha}
 :=\{\sigma_{\Gamma_\varepsilon}^{z,\alpha}
       <\tau_{3\rho}^{z,\alpha}\}.
\]
By Lemma~\ref{lem:feedback-compatible-exit-rules}, each family
\((\theta_\varepsilon^{z,\alpha})_{\alpha\in\mathcal A}\) is feedback
compatible, so the stopped DPP
\eqref{eq:abstract-controlled-stopped-DPP} applies.  Moreover,
\[
 P(H_-^{z,\alpha})\geq p_0,
 \qquad
 P(H_+^{z,\alpha})\geq p_0
\]
for every \(\alpha\in\mathcal A\) and
\(P\in\mathcal P_{\mathcal Q}\).

For the stopping rule associated with \(\Gamma_-\), continuity of the
state paths and closedness of \(\Gamma_-\) imply that on
\(H_-^{z,\alpha}\) the stopped point belongs to \(\Gamma_-\), whereas on
the complementary event it lies in
\(\overline{\mathcal C_{3\rho}(\zeta)}\).  Hence, for every \(\alpha\)
and \(P\),
\begin{align*}
 E^P\!\left[
  v\bigl(Z_{\theta_-^{z,\alpha}}^{z,\alpha}\bigr)
 \right]
 &\leq
 \mu P(H_-^{z,\alpha})
 +M\bigl(1-P(H_-^{z,\alpha})\bigr)\\
 &\leq M-p_0(M-\mu).
\end{align*}
The bound is uniform in \(P\) and \(\alpha\); therefore the stopped DPP
gives
\begin{equation}\label{eq:KS-upper-reduction}
 v(z)\leq M-p_0(M-\mu),
 \qquad z\in\mathcal C_\rho(\zeta).
\end{equation}

Applying instead the stopping rule associated with \(\Gamma_+\), we
obtain, again for every \(\alpha\) and \(P\),
\begin{align*}
 E^P\!\left[
  v\bigl(Z_{\theta_+^{z,\alpha}}^{z,\alpha}\bigr)
 \right]
 &\geq
 \mu P(H_+^{z,\alpha})
 +m\bigl(1-P(H_+^{z,\alpha})\bigr)\\
 &\geq m+p_0(\mu-m).
\end{align*}
This lower bound is also uniform in \(P\) and \(\alpha\), so it survives
first the supremum over \(P\) and then the infimum over \(\alpha\) in the
stopped DPP.  Consequently,
\begin{equation}\label{eq:KS-lower-reduction}
 v(z)\geq m+p_0(\mu-m),
 \qquad z\in\mathcal C_\rho(\zeta).
\end{equation}

Combining \eqref{eq:KS-upper-reduction} and
\eqref{eq:KS-lower-reduction}, for arbitrary
\(z,z'\in\mathcal C_\rho(\zeta)\),
\begin{align*}
 v(z)-v(z')
 &\leq
 M-p_0(M-\mu)-m-p_0(\mu-m)\\
 &=(1-p_0)(M-m).
\end{align*}
Interchanging \(z\) and \(z'\) gives
\begin{equation}\label{eq:controlled-oscillation-decay}
 \osc_{\mathcal C_\rho(\zeta)}v
 \leq q\,\osc_{\mathcal C_{4\rho}(\zeta)}v,
 \qquad
 q:=1-p_0\in(0,1).
\end{equation}
The use of a median is what makes both one-sided estimates available
simultaneously; with the midpoint \((M+m)/2\), only one of the two level
sets is guaranteed to have the required measure.

We now iterate \eqref{eq:controlled-oscillation-decay}.  Choose a fixed
geometric constant \(c_0\in(0,1)\), for instance \(c_0=1/4\), such that
\[
 \overline{\mathcal C_{4c_0r}(\zeta)}
 \subset\mathcal C_{4r}(z_0),
 \qquad \zeta\in\mathcal C_r(z_0).
\]
Set \(R_0=c_0r\) and \(R_n=4^{-n}R_0\).  Since
\(\mathcal C_{4R_{n+1}}(\zeta)=\mathcal C_{R_n}(\zeta)\), repeated
application of \eqref{eq:controlled-oscillation-decay} yields
\begin{equation}\label{eq:controlled-oscillation-iteration}
 \osc_{\mathcal C_{R_n}(\zeta)}v
 \leq q^{n+1}
 \osc_{\mathcal C_{4R_0}(\zeta)}v
 \leq q^{n+1}
 \osc_{\mathcal C_{4r}(z_0)}v.
\end{equation}
Choose
\[
 \eta:=\min\left\{\frac12,
       \frac{-\log q}{\log4}\right\}\in(0,1).
\]
Then \(q\leq4^{-\eta}\).

Let \(z=(\ell,x),z'=(k,y)\in\mathcal C_r(z_0)\), and put
\(\delta=d_{\rm p}(z,z')\).  If \(\delta<R_0/2\), choose \(n\geq0\)
so that
\[
 R_{n+1}<2\delta\leq R_n
\]
and set \(\zeta_*=(\max\{\ell,k\},x)\).  Then
\(z,z'\in\mathcal C_{R_n}(\zeta_*)\), while
\(\overline{\mathcal C_{4R_0}(\zeta_*)}\subset
  \mathcal C_{4r}(z_0)\).  Hence
\eqref{eq:controlled-oscillation-iteration} gives
\begin{align*}
 \abs{v(z)-v(z')}
 &\leq q^{n+1}
 \osc_{\mathcal C_{4r}(z_0)}v\\
 &\leq 4^{-(n+1)\eta}
 \osc_{\mathcal C_{4r}(z_0)}v\\
 &< C\,
 \osc_{\mathcal C_{4r}(z_0)}v
 \left(\frac{\delta}{r}\right)^\eta,
\end{align*}
where \(C\) depends only on the fixed cylinder geometry.  If
\(\delta\geq R_0/2\), the same inequality follows from the trivial bound
by \(\osc_{\mathcal C_{4r}(z_0)}v\), after increasing \(C\).  This proves
\eqref{eq:controlled-parabolic-Holder}.
\end{proof}

\begin{theorem}
\label{thm:controlled-positive-time}
Suppose Assumption~\ref{ass:uniform-parabolic} holds.
Fix \(t_*\in\R\) and \(\varphi\in\BUC_b(\R^d)\), and define
\begin{equation}\label{eq:controlled-remaining-value}
 v^{t_*}(\ell,x)
 :=\bS_{t_*-\ell,t_*}\varphi(x),
 \qquad \ell>0.
\end{equation}
Then \(v^{t_*}\) is locally parabolically H\"older continuous on
\((0,\infty)\times\R^d\).  More precisely, for every \(\delta>0\)
and \(R<\infty\), there are \(\eta\in(0,1)\) and
\(C_{\delta,R}<\infty\), independent of \(t_*\) and \(\varphi\),
such that
\begin{align}
 &\abs{v^{t_*}(\ell,x)-v^{t_*}(\ell',y)}
 \notag\\
 &\qquad\leq
 C_{\delta,R}\norm{\varphi}_\infty
 \left[
  \bigl(\abs{x-y}+\abs{\ell-\ell'}^{1/2}\bigr)\wedge1
 \right]^\eta
 \label{eq:controlled-uniform-Holder}
\end{align}
for \(\ell,\ell'\geq\delta\) and \(x,y\in B_R\).
\end{theorem}

\begin{proof}
Write \(v=v^{t_*}\).  Proposition~\ref{prop:controlled-continuity} implies
that \(v\) is Borel measurable on \((0,\infty)\times\R^d\), and the
infimum--supremum representation gives
\[
 \abs{v(\ell,x)}\leq\norm{\varphi}_\infty,
 \qquad \ell>0,\ x\in\R^d.
\]
Only this measurability and boundedness, rather than the continuity modulus
of \(\varphi\), will be used below.

Fix \((\ell,x)\in(0,\infty)\times\R^d\), and put
\[
 s_\ell:=t_*-\ell,
 \qquad
 \widetilde X_q^\alpha
 :=X_{s_\ell+q}^{s_\ell,x,\alpha},
 \qquad 0\leq q\leq\ell.
\]
Let \((\theta^\alpha)_{\alpha\in\cA}\) be any feedback-compatible
stopping rule bounded by the first exit of the backward space--time process
\[
 Z_q^\alpha:=(\ell-q,\widetilde X_q^\alpha)
\]
from a cylinder whose closure is contained in
\((0,\infty)\times\R^d\).  In particular,
\(0\leq\theta^\alpha\leq\ell\).  Applying
Proposition~\ref{prop:compatible-stopped-DPP} with initial time \(s_\ell\)
and terminal time \(t_*\) gives
\begin{align*}
 v(\ell,x)
 &=\inf_{\alpha\in\cA}\sup_{P\in\cP_{\cQ}}
 E^P\!\left[
  \bS_{s_\ell+\theta^\alpha,t_*}\varphi
  \bigl(\widetilde X_{\theta^\alpha}^\alpha\bigr)
 \right]\\
 &=\inf_{\alpha\in\cA}\sup_{P\in\cP_{\cQ}}
 E^P\!\left[
  v\bigl(\ell-\theta^\alpha,
    \widetilde X_{\theta^\alpha}^\alpha\bigr)
 \right],
\end{align*}
because \(s_\ell+\theta^\alpha=t_*-(\ell-\theta^\alpha)\).  This is
exactly \eqref{eq:abstract-controlled-stopped-DPP}.  In particular, it
holds for the compact-set hitting-before-exit rules used in
Lemma~\ref{lem:controlled-KS-smoothing}.

Let \(O\Subset\R^d\) be bounded and open.  Under
\(P\in\cP_{\cQ}\), use the elapsed-time density
\[
 Q_u^P:=\widehat Q_u,\qquad
 d\langle B\rangle_u=Q_u^P\,du,\qquad
 Q_u^P\in\cQ\quad du\times P\text{-a.e.},
\]
as in \eqref{eq:pathwise-covariance-identity} and \eqref{eq:qv}.
Equation \eqref{eq:controlled-SDE}, with initial calendar time
\(s_\ell\), directly gives, up to the exit of
\(\widetilde X^\alpha\) from \(O\),
\[
 \widetilde X_q^\alpha
 =x+\int_0^q
 b_{\alpha_u,Q_u^P}(s_\ell+u,\widetilde X_u^\alpha)\,du
 +M_q^{\alpha,P},
\]
where the continuous local \(P\)-martingale is
\[
 M_q^{\alpha,P}
 :=\sum_{j=1}^m\int_0^q
 \sigma_j(s_\ell+u,\widetilde X_u^\alpha,\alpha_u)\,dB_u^j.
\]
Its quadratic-variation density satisfies
\[
 \frac{d\langle M^{\alpha,P}\rangle_u}{du}
 ={c_{\alpha_u,Q_u^P}}
   (s_\ell+u,\widetilde X_u^\alpha),
 \qquad du\times P\text{-a.e.}
\]
Assumption~\ref{ass:uniform-parabolic} therefore yields
\[
 \abs{b_{\alpha_u,Q_u^P}
   (s_\ell+u,\widetilde X_u^\alpha)}\leq B_O,
 \qquad
 \underline\lambda I_d
 \leq\frac{d\langle M^{\alpha,P}\rangle_u}{du}
 \leq\overline\lambda I_d,
\]
uniformly in \(\alpha\), \(P\), \(t_*\) and \(\ell\).  Hence all the
hypotheses of Lemma~\ref{lem:controlled-KS-smoothing} hold locally, with
\(B=B_O\), \(\lambda=\underline\lambda\) and
\(\Lambda=\overline\lambda\).

Fix now \(\delta>0\) and \(R<\infty\).  Take
\[
 O:=B_{R+2},
 \qquad
 r_0:=\min\left\{\frac14,\frac{\sqrt\delta}{8}\right\}.
\]
Then, for every \(k\geq\delta\) and every \(x_0\in B_R\),
\[
 \overline{\mathcal C_{4r_0}(k,x_0)}
 \subset(0,\infty)\times O.
\]
Indeed, \(16r_0^2\leq\delta/4\) and \(4r_0\leq1\).  Therefore
Lemma~\ref{lem:controlled-KS-smoothing}, together with
\(\osc v\leq2\norm{\varphi}_\infty\), gives constants
\(\eta\in(0,1)\) and \(C_0<\infty\), depending only on \(\delta\),
\(R\) and the structural coefficient bounds, such that
\begin{equation}\label{eq:positive-time-local-holder-step}
 \abs{v(z)-v(z')}
 \leq
 2C_0\norm{\varphi}_\infty
 \left(\frac{d_{\rm p}(z,z')}{r_0}\right)^\eta
\end{equation}
whenever \(z,z'\in\mathcal C_{r_0}(k,x_0)\), with \(k\geq\delta\) and
\(x_0\in B_R\).  The same constants work for every \(t_*\).

Finally let \(\ell,\ell'\geq\delta\), \(x,y\in B_R\), and set
\[
 \Delta:=\abs{x-y}+\abs{\ell-\ell'}^{1/2}.
\]
If \(\Delta<r_0\), set \(k=\max\{\ell,\ell'\}\) and
\(\zeta=(k,x)\).  Then both \((\ell,x)\) and \((\ell',y)\) belong to
\(\mathcal C_{r_0}(\zeta)\), so
\eqref{eq:positive-time-local-holder-step} gives
\[
 \abs{v(\ell,x)-v(\ell',y)}
 \leq
 2C_0r_0^{-\eta}\norm{\varphi}_\infty\Delta^\eta.
\]
If \(\Delta\geq r_0\), then \(r_0\leq1\) and boundedness give
\[
 \abs{v(\ell,x)-v(\ell',y)}
 \leq2\norm{\varphi}_\infty
 \leq
 2r_0^{-\eta}\norm{\varphi}_\infty
 (\Delta\wedge1)^\eta.
\]
Combining the two cases proves
\eqref{eq:controlled-uniform-Holder}, for example with
\[
 C_{\delta,R}
 :=2\max\{C_0,1\}r_0^{-\eta}.
\]
This constant is independent of \(t_*\) and \(\varphi\).  Since every
compact subset of \((0,\infty)\times\R^d\) is contained in
\([\delta,\infty)\times B_R\) for suitable \(\delta>0\) and \(R<\infty\),
the asserted local parabolic H\"older continuity follows.
\end{proof}

\begin{remark}
(i) The estimate \eqref{eq:controlled-uniform-Holder} is
uniform only when the remaining times are bounded away from zero.  No
H\"older regularity at \(\ell=0\), and no uniform bound on the H\"older
constants as \(\delta\downarrow0\), is asserted.

(ii)
The estimate in Theorem~\ref{thm:controlled-positive-time} is independent
of the modulus of continuity of the terminal payoff.  The theorem is
nevertheless stated for \(\varphi\in\BUC_b(\R^d)\), because this is the
function domain on which Proposition~\ref{prop:compatible-stopped-DPP}
has been established in the present canonical control construction.

{
A bounded-Borel extension is not needed for the results of this paper.
If, for every \(\varphi\in B_b(\R^d)\) and \(t_*\in\R\), the value
\(v^{t_*}\) defined by the same infimum--supremum formula is jointly
universally measurable and satisfies
\eqref{eq:abstract-controlled-stopped-DPP} for every feedback-compatible
stopping rule bounded by a local cylinder exit, including the compact-set
hitting-before-exit rules, then the proof of
Theorem~\ref{thm:controlled-positive-time} gives
\[
 \bS_{s,t}\bigl(B_b(\R^d)\bigr)\subset C_b(\R^d),\qquad s<t.
\]
Expectations of universally measurable stopped payoffs are interpreted
under the corresponding completed laws.  These measurability and
stopped-DPP prerequisites must be established separately; a standard-Borel
encoding of controls alone does not supply them.  Once they are available,
no additional Krylov--Safonov estimate is required.
}
\end{remark}

\section[Local Harnack inequality and Harris oscillation contraction]{Local Harnack inequality and Harris oscillation contraction} 
\subsection{Bellman--Isaacs equation and local Harnack inequality}

We use the Frobenius norm
\[
 \|Q\|_{\mathrm F}
 :=\bigl(\tr(Q^{\mathsf T}Q)\bigr)^{1/2}.
\]
For \((q,x,p,X)\in\R\times\R^d\times\R^d\times\mathbb S^d\),
define
\begin{align}
 H_0(q,x,p,X)
 &:=\inf_{a\in A}\sup_{Q\in\cQ}
 \left\{b_{a,Q}(q,x)\cdot p
       +\frac12\tr\bigl(c_{a,Q}(q,x)X\bigr)\right\},
 \label{eq:lower-Isaacs-operator}\\
 H_0^\dagger(q,x,p,X)
 &:=\sup_{a\in A}\inf_{Q\in\cQ}
 \left\{b_{a,Q}(q,x)\cdot p
       +\frac12\tr\bigl(c_{a,Q}(q,x)X\bigr)\right\}.
 \label{eq:conjugate-Isaacs-operator}
\end{align}
The subscript \(0\) indicates zero running cost, and
\(X\in\mathbb S^d\) denotes the matrix argument corresponding
to the spatial Hessian.
Here \(q\) is calendar time, \(b_{a,Q}\) and \(c_{a,Q}\) are defined by \eqref{eq:controlled-effective-drift} and \eqref{eq:controlled-effective-covariance} respectively.  Compactness of \(A\) and \(\cQ\),
together with continuity of the coefficients, implies continuity of both
operators.  No equality between \(H_0\) and \(H_0^\dagger\) is assumed.

Put
\[
 \mathsf H_a(q,x,p):=\bigl(h_{ij}(q,x,a)\cdot p\bigr)_{i,j=1}^m.
\]
Using the \(G\)-operator introduced in
Remark~\ref{rem:Where-the-G-Brownian-motion-enters}, 
\begin{equation}\label{eq:controlled-G-form}
 H_0(q,x,p,X)
 =\inf_{a\in A}
 \left\{b(q,x,a)\cdot p
 +G\!\left(\sigma(q,x,a)^{\mathsf T}X\sigma(q,x,a)
             +2\mathsf H_a(q,x,p)\right)\right\}.
\end{equation}
Thus, for each fixed control action \(a\), the worst-case covariance matrix \(Q\in\cQ\) is selected in the inner supremum.
The following two lemmas provide a continuous approximate covariance
response and its realisation by an admissible canonical law.

\begin{lemma}
\label{lem:continuous-approximate-covariance-selector}
Let \(Z\) be a metric space and let
\(\mathcal K_0\subset\mathbb S_+^m\) be nonempty and compact.
For each \(z\in Z\), let \(\Gamma(z)\) be a nonempty compact convex
subset of \(\mathcal K_0\).  Assume that the set-valued map
\(z\mapsto\Gamma(z)\) is continuous in the Hausdorff distance, where
\[
 d_{\mathrm H}(K,L)
 :=\max\left\{
    \sup_{Q\in K}\inf_{Q'\in L}\|Q-Q'\|_{\mathrm F},
    \sup_{Q'\in L}\inf_{Q\in K}\|Q-Q'\|_{\mathrm F}
 \right\}.
\]
If \(C:Z\to\mathbb S^m\) is continuous, then for every
\(\varepsilon>0\) there is a continuous ordinary map
\(q^\varepsilon:Z\to\mathcal K_0\) such that
\begin{align}
 q^\varepsilon(z)&\in\Gamma(z),
 \label{eq:continuous-selector-feasible}\\
 \tr\!\bigl(C(z)q^\varepsilon(z)\bigr)
 &\geq\max_{Q\in\Gamma(z)}\tr\!\bigl(C(z)Q\bigr)-\varepsilon,
 \qquad z\in Z.
 \label{eq:continuous-selector-near-max}
\end{align}
In the present fixed-ambiguity model, \(\Gamma(z)\equiv\cQ\).
\end{lemma}

\begin{proof}
Let
\[
 R_0^2:=\max_{Q\in\mathcal K_0}\|Q\|_{\mathrm F}^2,
 \qquad \gamma:=\frac{\varepsilon}{1+R_0^2},
\]
and define
\[
 q^\varepsilon(z)
 :=\operatorname*{arg\,max}_{Q\in\Gamma(z)}
       \bigl\{\tr\!\bigl(C(z)Q\bigr)-\gamma\|Q\|_{\mathrm F}^2\bigr\}.
\]
A maximiser exists by compactness.  It is unique because the objective
is strictly concave and \(\Gamma(z)\) is convex.

We verify continuity directly.  Let \(z_n\to z\).  By compactness of
\(\mathcal K_0\), any subsequence of \(q^\varepsilon(z_n)\) has a
further subsequence converging to some \(\overline Q\).
Hausdorff continuity implies \(\overline Q\in\Gamma(z)\).
For each \(Q\in\Gamma(z)\), choose \(Q_n\in\Gamma(z_n)\) with
\(Q_n\to Q\).  By the optimality of \(q^\varepsilon(z_n)\) for the maximization
problem at \(z_n\),
\[
 \tr\!\bigl(C(z_n)q^\varepsilon(z_n)\bigr)
       -\gamma\|q^\varepsilon(z_n)\|_{\mathrm F}^2
 \geq \tr\!\bigl(C(z_n)Q_n\bigr)-\gamma\|Q_n\|_{\mathrm F}^2.
\]
Passing to the selected subsequence shows that \(\overline Q\)
maximises the regularised objective at \(z\).
Uniqueness implies \(\overline Q=q^\varepsilon(z)\).
Every convergent subsequence has this limit, so
\(q^\varepsilon(z_n)\to q^\varepsilon(z)\).

Finally, let \(Q^*\in\Gamma(z)\) maximise \(\tr\!\bigl(C(z)Q\bigr)\).
The defining inequality for \(q^\varepsilon(z)\) yields
\begin{align*}
 \tr\!\bigl(C(z)q^\varepsilon(z)\bigr)
 &\geq \tr\!\bigl(C(z)Q^*\bigr)
       +\gamma\bigl(\|q^\varepsilon(z)\|_{\mathrm F}^2
                    -\|Q^*\|_{\mathrm F}^2\bigr)\\
 &\geq\max_{Q\in\Gamma(z)}\tr\!\bigl(C(z)Q\bigr)-\gamma R_0^2\\
 &\geq\max_{Q\in\Gamma(z)}\tr\!\bigl(C(z)Q\bigr)-\varepsilon.
\end{align*}
This proves both assertions.
\end{proof}

For later use, fix a terminal calendar time \(t_*\).
If \(\psi=\psi(\ell,x)\) is a smooth test function written in
remaining-time coordinates, define
\begin{equation}\label{eq:test-covariance-matrix}
 \begin{split}
 C_\psi(q,x,a)
 :=\Bigl(&h_{ij}(q,x,a)\cdot D_x\psi(t_*-q,x)\\
 &+\frac12\bigl[
   \sigma(q,x,a)^{\mathsf T}D_x^2\psi(t_*-q,x)
   \sigma(q,x,a)\bigr]_{ij}\Bigr)_{i,j=1}^m.
 \end{split}
\end{equation}
Then
\[
 \mathcal L_{a,Q}^{q}[\psi(t_*-q,\cdot)](x)
 =b(q,x,a)\cdot D_x\psi(t_*-q,x)+{\tr\!\bigl(C_\psi(q,x,a)Q\bigr)}.
\]
Consequently, applying the preceding lemma to \(C_\psi\) gives a
continuous approximate maximiser of the generator over \(Q\in\cQ\).
A locally defined test function will first be extended, using a smooth
cutoff, without changing it on the compact region where it is used.

\begin{lemma}
\label{lem:elementary-covariance-realisation}
Suppose Assumption~\ref{ass:uniform-parabolic} holds, and let
\(\widehat q:\R\times\R^d\times A\to\cQ\) be continuous.
For every \(s\in\R\), \(x\in\R^d\), \(0<T<\infty\), and
\(\alpha\in\cA\), there exists
\(P^{s,x,\alpha,\widehat q,T}\in\cP_{\cQ}\) such that
\begin{equation}\label{eq:elementary-covariance-realisation}
 \frac{d\langle B\rangle_u}{du}
 =\widehat q\bigl(s+u,X_{s+u}^{s,x,\alpha},\alpha_u\bigr),
 \quad du\times P^{s,x,\alpha,\widehat q,T}\text{-a.e. on }[0,T].
\end{equation}
Here \(\alpha_u\) is the elapsed-time feedback rule, evaluated on the
canonical noise path.  
\end{lemma}

\begin{proof}
We construct a joint law of the canonical noise and the state, and then
identify its state coordinate with the canonical state map.

\smallskip
\noindent\emph{Step 1: a frozen action.}
Fix a starting calendar time \(r\), a state \(y\), and an action
\(a\in A\).  On an auxiliary filtered probability space, consider the
joint weak equation
\begin{equation}\label{eq:joint-response-SDE}
 \begin{aligned}
 d\mathsf B_u
   &=\widehat q(r+u,Y_u,a)^{1/2}\,dW_u,\qquad \mathsf B_0=0,\\
 dY_u
   &=b_{a,\widehat q(r+u,Y_u,a)}(r+u,Y_u)\,du
     +\sigma(r+u,Y_u,a)\,d\mathsf B_u,
     \qquad Y_0=y.
 \end{aligned}
\end{equation}
The drift and diffusion matrix of \(U=(\mathsf B,Y)\) are
\[
 \mu_a(q,y)
 =\begin{pmatrix}0\\ b_{a,\widehat q(q,y,a)}(q,y)\end{pmatrix},
 \qquad
 \Sigma_a(q,y)
 =\begin{pmatrix}I_m\\ \sigma(q,y,a)\end{pmatrix}
       \widehat q(q,y,a)^{1/2}.
\]
Both are continuous.  Compactness of \(\cQ\) and the coefficient
bounds imply a common linear-growth bound.  The usual weak-existence
theorem for continuous coefficients of linear growth therefore applies;
see \cite{StroockVaradhan1979}.  This theorem does not require the joint
diffusion matrix to be nondegenerate.

For completeness, the compactness argument uses Euler approximations
and, on every finite horizon \(H\), the uniform bounds
\[
 E\!\left[\sup_{u\leq H}|U_u|^4\right]
 \leq C_H(1+|y|^4),
 \qquad
 E[|U_t-U_u|^4]
 \leq C_H(1+|y|^4)|t-u|^2.
\]
These follow from the Burkholder--Davis--Gundy inequality and
Gronwall's lemma, and also hold for the approximations.
They give tightness on continuous-path space.
After localisation to a compact set, continuity of \(\mu_a\) and
\(\Sigma_a\) identifies every limit as a solution of the associated
martingale problem.  The moment bound removes the localisation and
precludes explosion.

\smallskip
\noindent\emph{Step 2: a measurable choice of the frozen-action law.}
Let
\[
 \mathcal E:=C([0,\infty);\R^{m+d}),
\]
and let \(\mathfrak S(r,y,a)\) be the set of laws on \(\mathcal E\)
solving the martingale problem of
\eqref{eq:joint-response-SDE}, with initial point \((0,y)\).
Its sections are nonempty by Step~1 and compact in the weak topology:
the same moment estimates give tightness, and the martingale problem is
closed under weak convergence because its coefficients are continuous.

We spell out its measurability.  For \(f\in C_c^\infty(\R^{m+d})\),
write
\[
 \mathscr L_a^q f(b,y)
 :=\mu_a(q,y)\cdot Df(b,y)
   +\frac12\tr\!\left(
     \Sigma_a(q,y)\Sigma_a(q,y)^{\mathsf T}D^2f(b,y)\right).
\]
The martingale problem is characterised by the initial condition and
identities
\[
 \int_{\mathcal E}
 \left[f(U_t)-f(U_u)
       -\int_u^t\mathscr L_a^{r+v}f(U_v)\,dv\right]
 \Phi(U)\,\mathsf P(dU)=0,
\]
where one may restrict to rational \(u<t\), a countable family of
smooth compactly supported test functions dense together with their
first two derivatives on compact sets, and a countable determining
algebra of bounded continuous cylinder functions \(\Phi\) depending
only on the path up to \(u\).
Continuity of paths and a monotone-class argument recover all the
martingale identities from this countable family.
Each displayed integral is a Borel function of \((r,y,a,\mathsf P)\).
Thus
\[
 \{(r,y,a,\mathsf P):\mathsf P\in\mathfrak S(r,y,a)\}
\]
is Borel.  The Borel uniformisation theorem for nonempty compact
sections gives a Borel choice
\[
 (r,y,a)\longmapsto\mathsf K(r,y,a)\in\mathfrak S(r,y,a);
\]
see the Arsenin--Kunugui selection theorem in
\cite[Section~18]{Kechris1995}.
In particular, evaluating \(\mathsf K\) at a random starting time,
state, and previously selected action gives a measurable continuation
kernel.  No uniqueness in law of the joint equation is needed.

\smallskip
\noindent\emph{Step 3: successive holding intervals.}
Let \((\tau_k,\xi_k)\) describe an elementary rule \(\alpha\).
At \(\tau_k\), the noise history and the state \(Y_{\tau_k}\)
have already been constructed, and the action \(\xi_k\) is measurable
with respect to the available noise history.  Use the kernel
\[
 \mathsf K(s+\tau_k,Y_{\tau_k},\xi_k)
\]
to generate a continuation with this action frozen.
Add its noise increment to the noise history already constructed, and
stop the continuation at the next switching time \(\tau_{k+1}\),
evaluated on that concatenated noise path, or at \(T\), whichever
comes first.  The next switching time is an \(\mathbb F\)-stopping
time on the canonical noise space; the
stopping and concatenation maps on continuous-path space are Borel.
Hence this procedure defines a measurable kernel at every step.

Iterated integration of these kernels constructs the joint law up to
each successive switching time.  Conditioning on the past shows that
the martingale identities on adjacent holding intervals concatenate;
this is the usual conditioning-and-pasting argument for martingale
problems, as in \cite{ElKarouiTan2015}.
The coefficient and moment bounds of Step~1 are uniform in the chosen
action, so they also hold after any finite number of switches.
They give continuous nonexplosive limits at any putative finite
accumulation time.  Such an accumulation would produce a continuous
noise path on which the prescribed decision grid has infinitely many
switches before a finite time, contrary to its pathwise local
finiteness.  Thus the construction reaches \(T\).
The same procedure applies on each event of a countable
\(\F_\tau\)-measurable pasting partition.  Only the selected branch
is followed, and each branch has a locally finite grid.  This covers
the countable eventwise pastings used in the definition of \(\cA\).

Let \(\mathbf P\) denote the resulting joint law up to \(T\).
Under this law, \(\mathsf B\) is a continuous martingale and
\[
 d\langle\mathsf B\rangle_u
 =\widehat q(s+u,Y_u,\alpha_u(\mathsf B))\,du.
\]
Moreover, with stochastic integrals taken in the joint filtration,
\begin{align*}
 Y_u
 ={}&x+\int_0^u b(s+v,Y_v,\alpha_v(\mathsf B))\,dv\\
 &+\sum_{i,j=1}^m\int_0^u
   h_{ij}(s+v,Y_v,\alpha_v(\mathsf B))
   \,d\langle\mathsf B^i,\mathsf B^j\rangle_v\\
 &+\int_0^u\sigma(s+v,Y_v,\alpha_v(\mathsf B))\,d\mathsf B_v,
 \qquad 0\leq u\leq T.
\end{align*}
These identities also follow directly from the joint martingale
problem: the cross-variation of the noise and state martingales is
\(\widehat q\sigma^{\mathsf T}\,du\), so the difference between the
state martingale and \(\int\sigma\,d\mathsf B\) has zero quadratic
variation and starts at zero.

\smallskip
\noindent\emph{Step 4: projection to the canonical noise space.}
Project \(\mathbf P\) onto its noise coordinate.  Beyond \(T\),
append independent Brownian increments with any fixed covariance
\(Q_0\in\cQ\), and denote the resulting full-path law by \(P\).
The noise is a continuous martingale in its natural filtration, and
its absolutely continuous quadratic variation has density in \(\cQ\)
both before and after \(T\).  Therefore \(P\in\cP_{\cQ}\).

Under \(P\), the canonical state map of
Lemma~\ref{lem:controlled-canonical-state-map} solves the original
controlled equation.  Pull this solution back to the joint space of
\(\mathbf P\).  It and \(Y\) have the same initial state, the same
noise, and the same prescribed process \(\alpha(\mathsf B)\).
The state-Lipschitz bounds, the bounded quadratic-variation density,
and Gronwall's lemma give pathwise uniqueness for that equation.
Consequently,
\[
 Y_u=X_{s+u}^{s,x,\alpha}(\mathsf B),
 \qquad 0\leq u\leq T,\quad \mathbf P\text{-a.s.}
\]
Here uniqueness is used for the original equation with its fixed noise
and control, not for the possibly non-Lipschitz joint response equation.
Substituting this identity into the quadratic-variation formula proves
\eqref{eq:elementary-covariance-realisation}.
\end{proof}

\begin{proposition}
\label{prop:controlled-value-PDE}
Suppose Assumption~\ref{ass:uniform-parabolic} holds.
Fix \(t_*\in\R\) and \(\varphi\in\BUC_b(\R^d)\).
Let \(v=v^{t_*}\) be defined by
\eqref{eq:controlled-remaining-value} for \(\ell>0\), and set
\begin{equation}\label{eq:controlled-remaining-initial-condition}
 v(0,x):=\varphi(x),\qquad x\in\R^d.
\end{equation}
Then, for every \(L>0\), \(v\) is the unique bounded continuous
viscosity solution on \([0,L]\times\R^d\) of
\begin{equation}\label{eq:controlled-remaining-HJB}
 \partial_\ell v-H_0(t_*-\ell,x,Dv,D^2v)=0
 \quad\text{on }(0,L)\times\R^d,
\end{equation}
with initial condition
\eqref{eq:controlled-remaining-initial-condition}.
Continuity at \(\ell=0\) does not assert H\"older regularity of the
initial datum.

If \(M\in\R\) is constant and \(w:=M-v\), then \(w\) satisfies
\begin{equation}\label{eq:controlled-conjugate-HJB}
 \partial_\ell w-H_0^\dagger(t_*-\ell,x,Dw,D^2w)=0
 \quad\text{on }(0,L)\times\R^d,
\end{equation}
in the viscosity sense, with \(w(0,x)=M-\varphi(x)\).
\end{proposition}

\begin{proof}
\emph{Step 1: boundedness, continuity, and the initial value.}
The value representation gives \(|v|\leq\|\varphi\|_\infty\),
and Theorem~\ref{thm:controlled-positive-time} gives continuity for
\(\ell>0\).  Let
\[
 \omega_\varphi(r)
 :=\sup_{|x-y|\leq r}|\varphi(x)-\varphi(y)|.
\]
For every compact \(K\subset\R^d\), the short-time increment estimate
obtained from the controlled equation gives
\[
 \sup_{\substack{x\in K,\ \alpha\in\cA\\P\in\cP_{\cQ}}}
 E^P\!\left[
  |X_{t_*}^{t_*-\ell,x,\alpha}-x|^2\right]
 \leq C_K\ell,\qquad 0<\ell\leq1.
\]
Indeed, the drift contribution is bounded by
\(C_K\ell^2\) using the moment estimate, and the martingale
contribution by \(C_K\ell\).
For every \(r>0\), splitting according to whether the increment is
at most \(r\) yields
\[
 \sup_{x\in K}|v(\ell,x)-\varphi(x)|
 \leq\omega_\varphi(r)
       +2\|\varphi\|_\infty\frac{C_K\ell}{r^2}.
\]
Letting first \(\ell\downarrow0\) and then \(r\downarrow0\) proves
local uniform convergence to \(\varphi\), and hence joint continuity
at every point of \(\{0\}\times\R^d\).

\smallskip
\noindent\emph{Step 2: the local stopped DPP and its time scale.}
Fix \(z_0=(\ell_0,x_0)\) with \(\ell_0>0\), and put
\[
 s_0:=t_*-\ell_0,\qquad
 Y_u^\alpha:=X_{s_0+u}^{s_0,x_0,\alpha},\qquad
 Z_u^\alpha:=(\ell_0-u,Y_u^\alpha).
\]
Choose \(r>0\) with \(r^2<\ell_0\), and set
\[
 K_r:=[\ell_0-r^2,\ell_0]\times\overline{B_r(x_0)}.
\]
Let \(\tau_r^\alpha\) be the law-independent
\(\mathbb F\)-stopping-time representative of the capped spatial exit rule,
supplied by Lemma~\ref{lem:feedback-compatible-exit-rules}; thus
\[
 \tau_r^\alpha
 =r^2\wedge\inf\{u\geq0:|Y_u^\alpha-x_0|\geq r\}
 \quad P\text{-a.s. for every }P\in\cP_{\cQ}.
\]
For \(0<h<r^2\), let \(\theta_h^\alpha:=h\wedge\tau_r^\alpha\).
These are feedback-compatible \(\mathbb F\)-stopping times by the same lemma.
Proposition~\ref{prop:compatible-stopped-DPP} gives
\begin{equation}\label{eq:section8-local-stopped-DPP}
 v(z_0)
 =\inf_{\alpha\in\cA}\sup_{P\in\cP_{\cQ}}
 E^P[v(Z_{\theta_h^\alpha}^\alpha)].
\end{equation}
The remaining-time argument follows from
\(s_0+\theta_h^\alpha=t_*-(\ell_0-\theta_h^\alpha)\).

The local drift and covariance bounds imply
\[
 \sup_{\alpha,P}E^P\!\left[
   \sup_{u\leq h}|Y_{u\wedge\tau_r^\alpha}^\alpha-x_0|^2
 \right]\leq C_rh,
 \qquad 0<h<\min\{r^2,1\}.
\]
By continuity of the paths and Markov's inequality, after changing
\(C_r\),
\[
 \sup_{\alpha,P}P(\tau_r^\alpha\leq h)\leq C_rh.
\]
Consequently,
\begin{equation}\label{eq:section8-stopped-duration}
 \inf_{\alpha,P}E^P[\theta_h^\alpha]
 \geq h(1-C_rh)\geq\frac h2
\end{equation}
for all sufficiently small \(h\).

Let \(\psi\in C^{1,2}\) be defined on a neighbourhood of \(K_r\).
Write
\begin{align*}
 g_a(\ell,y)
 &:=-\partial_\ell\psi(\ell,y)
   +\max_{Q\in\cQ}
       \mathcal L_{a,Q}^{t_*-\ell}[\psi(\ell,\cdot)](y),\\
 g(\ell,y)
 &:=\min_{a\in A}g_a(\ell,y)\\
 &=-\partial_\ell\psi(\ell,y)
   +H_0(t_*-\ell,y,D\psi(\ell,y),D^2\psi(\ell,y)).
\end{align*}
These functions are continuous, jointly in the displayed variables.
For \(P\in\cP_{\cQ}\), write \(d\langle B\rangle_u=Q_u^P\,du\).
It\^o's formula, stopped at \(\theta_h^\alpha\), gives
\begin{equation}\label{eq:section8-stopped-Ito}
 \begin{split}
 &E^P[\psi(Z_{\theta_h^\alpha}^\alpha)]-\psi(z_0)\\
 &\quad=E^P\int_0^{\theta_h^\alpha}
 \Bigl[-\partial_\ell\psi(\ell_0-u,Y_u^\alpha)
  +\mathcal L_{\alpha_u,Q_u^P}^{s_0+u}
       [\psi(\ell_0-u,\cdot)](Y_u^\alpha)\Bigr]du.
 \end{split}
\end{equation}
The stopped stochastic integral has zero expectation: its quadratic
variation is integrable because the spatial derivatives of \(\psi\)
and the effective covariance are bounded on the stopping region.

\smallskip
\noindent\emph{Step 3: the viscosity subsolution inequality.}
Suppose \(\psi\) touches \(v\) from above at \(z_0\), with
\(\psi(z_0)=v(z_0)\).  Reduce \(r\) so that \(v\leq\psi\) on
\(K_r\).  We must prove \(g(z_0)\geq0\).
If instead \(g(z_0)<0\), choose \(a_0\in A\) and \(\kappa>0\)
with \(g_{a_0}(z_0)<-2\kappa\).
Continuity permits a further reduction of \(r\) such that
\[
 g_{a_0}(\ell,y)\leq-\kappa,
 \qquad (\ell,y)\in K_r.
\]
Use the constant action \(a_0\) until \(\theta_h^{a_0}\).
For every admissible law, the integrand in
\eqref{eq:section8-stopped-Ito} is at most \(-\kappa\), since its
covariance lies in \(\cQ\).
Therefore \eqref{eq:section8-local-stopped-DPP}, 
\eqref{eq:section8-stopped-duration} and \eqref{eq:section8-stopped-Ito} imply
\begin{align*}
 v(z_0)
 &\leq\sup_P E^P[v(Z_{\theta_h^{a_0}}^{a_0})]\\
 &\leq\sup_P E^P[\psi(Z_{\theta_h^{a_0}}^{a_0})]\\
 &\leq v(z_0)-\kappa\inf_P E^P[\theta_h^{a_0}]\\
 &\leq v(z_0)-\frac{\kappa h}{2},
\end{align*}
a contradiction.  Thus
\[
 \partial_\ell\psi(z_0)
 -H_0(t_*-\ell_0,x_0,D\psi(z_0),D^2\psi(z_0))\leq0.
\]
Notice that the integrand was bounded under every admissible law;
no restriction of the covariance supremum to constant-covariance laws
was made.

\smallskip
\noindent\emph{Step 4: the viscosity supersolution inequality.}
Suppose \(\psi\) touches \(v\) from below at \(z_0\), with equal
values there, and choose \(r\) so that \(v\geq\psi\) on \(K_r\).
We must prove \(g(z_0)\leq0\).
If \(g(z_0)>0\), choose \(\kappa>0\) and reduce \(r\) so that
\[
 g(\ell,y)\geq2\kappa,\qquad (\ell,y)\in K_r.
\]
Multiply the local test function by a smooth cutoff equal to one on
a neighbourhood of \(K_r\), with compact support in its domain of
definition, and extend the product by zero.
This gives a global \(C^{1,2}\) function agreeing with \(\psi\),
including its relevant derivatives, on \(K_r\).

Apply Lemma~\ref{lem:continuous-approximate-covariance-selector}
to \eqref{eq:test-covariance-matrix}, with error \(\kappa\) and
\(\Gamma\equiv\cQ\).  It gives a continuous map
\(\widehat q:\R\times\R^d\times A\to\cQ\) such that
\begin{equation}\label{eq:section8-positive-response}
 \begin{split}
 &-\partial_\ell\psi(\ell,y)
  +\mathcal L_{a,\widehat q(t_*-\ell,y,a)}^{t_*-\ell}
       [\psi(\ell,\cdot)](y)\\
 &\qquad\geq g_a(\ell,y)-\kappa
 \geq g(\ell,y)-\kappa\geq\kappa
 \end{split}
\end{equation}
for every \((\ell,y)\in K_r\) and \(a\in A\).

For each small \(h>0\), choose an \(h^2\)-optimal rule
\(\alpha^h\) in \eqref{eq:section8-local-stopped-DPP}, so that
\begin{equation}\label{eq:section8-near-optimal-control}
 \sup_{P\in\cP_{\cQ}}
 E^P[v(Z_{\theta_h^{\alpha^h}}^{\alpha^h})]
 \leq v(z_0)+h^2.
\end{equation}
Lemma~\ref{lem:elementary-covariance-realisation}, with horizon \(h\),
provides an admissible law \(P^h\) satisfying
\[
 Q_u^{P^h}
 =\widehat q(s_0+u,Y_u^{\alpha^h},\alpha_u^h)
 \quad du\times P^h\text{-a.e. on }[0,h].
\]
Use this law in \eqref{eq:section8-near-optimal-control}.
The lower contact inequality, \eqref{eq:section8-stopped-Ito}, and
\eqref{eq:section8-positive-response} give
\begin{align*}
 h^2
 &\geq E^{P^h}[v(Z_{\theta_h^{\alpha^h}}^{\alpha^h})]-v(z_0)\\
 &\geq E^{P^h}[\psi(Z_{\theta_h^{\alpha^h}}^{\alpha^h})]-\psi(z_0)\\
 &\geq\kappa E^{P^h}[\theta_h^{\alpha^h}]
 \geq\frac{\kappa h}{2},
\end{align*}
which is impossible for sufficiently small \(h\).
Hence
\[
 \partial_\ell\psi(z_0)
 -H_0(t_*-\ell_0,x_0,D\psi(z_0),D^2\psi(z_0))\geq0.
\]
This argument first fixes the approximately minimising feedback rule
and then realises a covariance response to that rule.  It therefore
retains the order \(\inf_a\sup_Q\), without an Isaacs interchange.

\smallskip
\noindent\emph{Step 5: comparison and uniqueness on finite strips.}
Fix \(L>0\).  We first check the structural condition for bounded-domain
parabolic comparison.  Put
\[
 \Sigma_{a,Q}(q,x):=\sigma(q,x,a)Q^{1/2},
 \qquad c_{a,Q}=\Sigma_{a,Q}\Sigma_{a,Q}^{\mathsf T}.
\]
Both \(b_{a,Q}\) and \(\Sigma_{a,Q}\) are Lipschitz in \(x\), with
constants uniform in \(q,a,Q\).
If \(\nu>0\), \(p=\nu(x-y)\), and \(X_1,X_2\in\mathbb S^d\)
satisfy the block inequality
\[
 -3\nu I_{2d}
 \leq\begin{pmatrix}X_1&0\\0&-X_2\end{pmatrix}
 \leq3\nu\begin{pmatrix}I_d&-I_d\\-I_d&I_d\end{pmatrix},
\]
then
\begin{align*}
 &(b_{a,Q}(q,x)-b_{a,Q}(q,y))\cdot p
       \leq C\nu|x-y|^2,\\
 &\frac12\tr\!\left(
       c_{a,Q}(q,x)X_1-c_{a,Q}(q,y)X_2\right)
       \leq\frac{3\nu}{2}
          \|\Sigma_{a,Q}(q,x)-\Sigma_{a,Q}(q,y)\|_{\mathrm F}^2
       \leq C\nu|x-y|^2.
\end{align*}
Taking the two optimisations preserves this common estimate, giving
\[
 H_0(q,x,p,X_1)-H_0(q,y,p,X_2)\leq C\nu|x-y|^2.
\]
The operator \(-H_0(t_*-\ell,x,p,X)\) is independent of the
solution value and uniformly elliptic in \(X\), and hence is
proper in the sense of \cite{CIL1992}.  Together with continuity
and the preceding structural estimate, these properties verify
the hypotheses of the bounded-domain parabolic comparison
theorem \cite[Theorem~8.2]{CIL1992}.

To handle spatial infinity, the linear-growth drift and bounded
effective covariance give a constant \(C_*\geq0\) such that
\[
 2b_{a,Q}(q,x)\cdot x+\tr(c_{a,Q}(q,x))
 \leq C_*(1+|x|^2)
\]
for every \(q,x,a,Q\).  Define
\[
 \chi(\ell,x):=e^{(C_*+1)\ell}(1+|x|^2).
\]
Then
\begin{equation}\label{eq:section8-comparison-barrier}
 \partial_\ell\chi
 -\sup_{a,Q}\mathcal L_{a,Q}^{t_*-\ell}\chi\geq\chi.
\end{equation}
Let \(U,V\) be bounded continuous sub- and supersolutions on
\([0,L]\times\R^d\), with \(U(0,\cdot)\leq V(0,\cdot)\).
For \(\varepsilon>0\), the function
\(V_\varepsilon:=V+\varepsilon\chi\) is a supersolution.
Indeed, the elementary increment bound
\[
 \begin{split}
 &H_0(q,x,p+\varepsilon D\chi,X+\varepsilon D^2\chi)
       -H_0(q,x,p,X)\\
 &\qquad\leq\varepsilon
       \sup_{a,Q}\mathcal L_{a,Q}^{q}\chi
 \end{split}
\]
and a lower touching test for \(V_\varepsilon\), with
\(\varepsilon\chi\) subtracted, reduce this assertion to
\eqref{eq:section8-comparison-barrier} and the supersolution property
of \(V\).
For every sufficiently large radius \(R\), boundedness gives
\(U\leq V_\varepsilon\) on \([0,L]\times\partial B_R\);
the same inequality holds at \(\ell=0\).
Bounded-domain comparison implies \(U\leq V_\varepsilon\) in
\([0,L)\times B_R\).  Such radii can be chosen arbitrarily large.
Letting \(\varepsilon\downarrow0\), and using continuity at \(L\),
proves comparison on \([0,L]\times\R^d\).
This proves uniqueness in the stated class.

\smallskip
\noindent\emph{Step 6: the conjugate equation.}
The definitions imply
\[
 -H_0(q,x,-p,-X)=H_0^\dagger(q,x,p,X).
\]
If \(\eta\) touches \(w=M-v\) from above, then \(M-\eta\) touches
\(v\) from below.  The supersolution inequality for \(v\) therefore
implies
\[
 \partial_\ell\eta
 -H_0^\dagger(t_*-\ell,x,D\eta,D^2\eta)\leq0.
\]
A lower touching test for \(w\) gives the reverse inequality by using
the subsolution property of \(v\).  Thus \(w\) solves
\eqref{eq:controlled-conjugate-HJB} in the viscosity sense, and its
initial value follows directly from \(v(0,x)=\varphi(x)\).
\end{proof}

\begin{proposition}
\label{prop:controlled-local-Harnack}
Suppose Assumption~\ref{ass:uniform-parabolic} holds, and let
\(O\subset\R^d\) be bounded and open.
Consider cylinders
\begin{align*}
 \mathcal C_{2\rho}
 &:=(\ell_0-4\rho^2,\ell_0]\times B_{2\rho}(x_0),\\
 \mathcal C_\rho^-
 &:=(\ell_0-3\rho^2,\ell_0-2\rho^2)\times B_\rho(x_0),\\
 \mathcal C_\rho^+
 &:=(\ell_0-\rho^2,\ell_0]\times B_\rho(x_0),
\end{align*}
with \(\overline{\mathcal C_{2\rho}}\subset(0,\infty)\times O\).
There is a constant \(C_H<\infty\), depending only on
\[
 d,\quad\underline\lambda,\quad\overline\lambda,
 \quad\rho B_O,
\]
and the fixed relative cylinder geometry, such that every continuous
viscosity solution \(z\) of either
\eqref{eq:controlled-remaining-HJB} or
\eqref{eq:controlled-conjugate-HJB}, defined on a neighbourhood of
\(\overline{\mathcal C_{2\rho}}\) and nonnegative on
\(\mathcal C_{2\rho}\), satisfies
\begin{equation}\label{eq:controlled-local-Harnack}
 \sup_{\mathcal C_\rho^-}z
 \leq C_H\inf_{\mathcal C_\rho^+}z.
\end{equation}
The constant is independent of \(t_*\) and of the particular solution;
in particular, it is independent of the terminal payoff for the value
solutions above.
\end{proposition}

\begin{proof}
Put \(\lambda_0=\underline\lambda/2\) and
\(\Lambda_0=\overline\lambda/2\), and define the Pucci operators by
\[
 \mathcal M^-_{\lambda_0,\Lambda_0}(X)
 :=\inf_{\lambda_0 I_d\leq A_0\leq\Lambda_0 I_d}\tr(A_0X),
 \qquad
 \mathcal M^+_{\lambda_0,\Lambda_0}(X)
 :=\sup_{\lambda_0 I_d\leq A_0\leq\Lambda_0 I_d}\tr(A_0X).
\]
For every \(x\in O\), \(q\in\R\), \(a\in A\), and \(Q\in\cQ\),
\[
 \mathcal M^-_{\lambda_0,\Lambda_0}(X)-B_O|p|
 \leq b_{a,Q}(q,x)\cdot p+\tfrac12\tr(c_{a,Q}(q,x)X)
 \leq\mathcal M^+_{\lambda_0,\Lambda_0}(X)+B_O|p|.
\]
Both \(\inf_a\sup_Q\) and \(\sup_a\inf_Q\) preserve these bounds.
Thus, writing \(\mathcal H\) for either \(H_0\) or \(H_0^\dagger\),
\[
 \mathcal M^-_{\lambda_0,\Lambda_0}(X)-B_O|p|
 \leq\mathcal H(q,x,p,X)
 \leq\mathcal M^+_{\lambda_0,\Lambda_0}(X)+B_O|p|.
\]
If \(\psi\) touches \(z\) from above, its viscosity subsolution
inequality gives
\[
 \partial_\ell\psi
 \leq\mathcal M^+_{\lambda_0,\Lambda_0}(D^2\psi)+B_O|D\psi|.
\]
If it touches \(z\) from below, the supersolution inequality gives
\[
 \partial_\ell\psi
 \geq\mathcal M^-_{\lambda_0,\Lambda_0}(D^2\psi)-B_O|D\psi|.
\]
Hence \(z\) belongs to the same two-sided parabolic Pucci class for
both equations.

Make the change of variables
\[
 \widehat z(\tau,y):=z(\ell_0+\rho^2\tau,x_0+\rho y).
\]
The outer cylinder becomes \((-4,0]\times B_2\), and the two
comparison cylinders become \((-3,-2)\times B_1\) and
\((-1,0]\times B_1\).  In the viscosity sense,
\begin{align*}
 \partial_\tau\widehat z
   -\mathcal M^+_{\lambda_0,\Lambda_0}(D_y^2\widehat z)
   -\rho B_O|D_y\widehat z|&\leq0,\\
 \partial_\tau\widehat z
   -\mathcal M^-_{\lambda_0,\Lambda_0}(D_y^2\widehat z)
   +\rho B_O|D_y\widehat z|&\geq0.
\end{align*}
The scaled function belongs to
\(\mathcal S(\lambda_0,\Lambda_0,\rho B_O)\), the continuous viscosity
Pucci class in
\cite[Definition~2.2 and (2.5), p.~4, of the arXiv version]{BanerjeeGarofalo2015}.
The interior Harnack inequality
\cite[Theorem~2.10, p.~7, of the arXiv version]{BanerjeeGarofalo2015}
applies to nonnegative functions in this class, including the bounded
first-order term.  See also
\cite{Wang1992,Krylov1987,Lieberman1996} for the classical background.

For completeness, its two-time-slice formulation gives our slab
estimate by a uniformly finite interior chain. 
Fix an integer \(N>\max\{16,12d\}\), independently of the
points to be compared.  For any
\((u,x)\in(-3,-2)\times B_1\) and
\((v,y)\in(-1,0)\times B_1\), put
\[
 \begin{gathered}
 \Delta:=v-u\in(1,3),\qquad r:=\sqrt{\Delta/N},\\
 t_j:=u+j\Delta/N,\qquad x_j:=x+j(y-x)/N,
 \quad j=0,\ldots,N.
 \end{gathered}
\]
For \(1\leq j\leq N\), apply the cited theorem with time parameter
\(2\) and interior parameter \(1/2\), in the spatial cube of
half-side \(r\) centred at \(x_j\), with bottom time
\(t_{j-1}-r^2\).  The two comparison times are \(t_{j-1}\) and
\(t_j\), and \(|x_j-x_{j-1}|<2/N<r/2\).
Every cube lies in \(B_2\), since its centre lies in \(B_1\) and
\(\sqrt d\,r<1/2\); all bottom times exceed \(-4\).
Thus all these cylinders lie in \((-4,0)\times B_2\).
The local Harnack inequality, with a common constant
\(C_{\mathrm{loc}}\geq1\) independent of the chosen points,
gives by iteration
\[
 \widehat z(u,x)
 \leq C_{\mathrm{loc}}\widehat z(t_1,x_1)
 \leq\cdots\leq
 C_{\mathrm{loc}}^N\widehat z(v,y).
\]
Setting \(C_H:=C_{\mathrm{loc}}^N\), taking the supremum over
\((u,x)\in(-3,-2)\times B_1\), and then the infimum over
\((v,y)\in(-1,0)\times B_1\), we obtain
\[
 \sup_{(-3,-2)\times B_1}\widehat z
 \leq C_H\inf_{(-1,0)\times B_1}\widehat z.
\]
The constant depends only on
\(d,\lambda_0,\Lambda_0,\rho B_O\) and the fixed geometry.
Continuity at the top time slice permits its inclusion in the infimum.
Returning to the original variables proves
\eqref{eq:controlled-local-Harnack}.
All bounds used above are uniform in calendar time and in the action
and covariance indices, which proves the asserted uniformity.
\end{proof}


\subsection{Consequence: controlled local oscillation contraction}

We now pass from the local Harnack inequality to a strict oscillation
contraction on an arbitrary compact set.  The construction below fixes
both the comparison region and the number of Harnack steps independently
of the calendar-time location and the terminal payoff.

\begin{proposition}
\label{prop:controlled-local-oscillation}
Suppose Assumption~\ref{ass:uniform-parabolic} holds.  Let \(a>0\),
let \(O\subset\R^d\) be bounded, open, and connected, and let
\(K\subset O\) be nonempty and compact.  There exist a compact set
\[
 \mathfrak C_{a,K,O}\subset(0,a]\times O,
 \qquad
 \{a\}\times K\subset\mathfrak C_{a,K,O},
\]
and a constant \(\vartheta_{a,K,O}\in(0,1)\) such that, for every
\(s\in\R\) and \(\varphi\in\BUC_b(\R^d)\),
\begin{equation}\label{eq:controlled-local-oscillation}
 \osc_K(\bS_{s,s+a}\varphi)
 \leq\vartheta_{a,K,O}
       \osc_{\mathfrak C_{a,K,O}}v_s,
 \qquad
 v_s(\ell,x):=\bS_{s+a-\ell,s+a}\varphi(x).
\end{equation}
The comparison set depends only on \(a,K,O\), and the contraction
constant depends only on these data and the structural coefficient
bounds.  Both are independent of \(s\) and \(\varphi\).
\end{proposition}

\begin{proof}
\emph{Step 1: a common compact comparison region.}
Fix \(x_*\in O\).  By compactness, choose finitely many balls
\(B_{r_i}(x_i)\), \(1\leq i\leq m_0\), covering \(K\), with
\[
 \overline{B_{r_i}(x_i)}\subset O.
\]
Since an open connected subset of \(\R^d\) is polygonally connected,
for each \(i\) there is a polygonal path \(\gamma_i\) in \(O\)
joining \(x_*\) to \(x_i\).  Define
\[
 E:=\bigcup_{i=1}^{m_0}\gamma_i([0,1])
       \;\cup\;
       \bigcup_{i=1}^{m_0}\overline{B_{r_i}(x_i)}.
\]
Then \(E\) is compact, \(\{x_*\}\cup K\subset E\subset O\), and
\[
 d_0:=\operatorname{dist}(E,O^c)>0.
\]
If \(x\in K\cap B_{r_i}(x_i)\), concatenate \(\gamma_i\) with
the line segment from \(x_i\) to \(x\).  This path lies in \(E\)
and has length at most
\[
 L:=\max_{1\leq i\leq m_0}
       \bigl\{\operatorname{length}(\gamma_i)+r_i\bigr\}<\infty.
\]
Parametrising each such path at constant speed gives, for every
\(x\in K\), a path \(\gamma^x:[0,1]\to E\) satisfying
\[
 \gamma^x(0)=x_*,\qquad \gamma^x(1)=x,
 \qquad
 |\gamma^x(t)-\gamma^x(u)|\leq L|t-u|.
\]
A path of zero length is taken to be constant.

Let
\[
 E_{d_0/2}:=
 \{y\in\R^d:\operatorname{dist}(y,E)\leq d_0/2\},
\]
and set
\begin{equation*}
 \mathfrak C_{a,K,O}:=[a/4,a]\times E_{d_0/2}.
\end{equation*}
The set \(E_{d_0/2}\) is compact and contained in \(O\).
Consequently, \(\mathfrak C_{a,K,O}\) has the compactness and
containment properties asserted in the proposition.  In particular,
it stays a positive distance from remaining time zero.  All choices
in this step are made independently of \(s\) and \(\varphi\).

\smallskip
\noindent\emph{Step 2: uniformly finite forward Harnack chains.}
Choose an integer \(N\geq2\) sufficiently large that the radius
\(\rho:=\sqrt{a/(5N)}\) satisfies
\[
 2\rho<\frac{d_0}{2},\qquad \frac{L}{N}<\rho.
\]
Define the remaining-time levels
\[
 \ell_j:=\frac a2+\frac{ja}{2N},\qquad j=0,\ldots,N.
\]
Thus
\[
 \ell_0=\frac a2,\qquad \ell_N=a,\qquad
 \ell_j-\ell_{j-1}=\frac52\rho^2.
\]
Fix \(x\in K\), and put
\[
 x_j:=\gamma^x(j/N),\qquad z_j:=(\ell_j,x_j),
 \qquad j=0,\ldots,N.
\]
Then
\[
 x_0=x_*,\qquad x_N=x,\qquad
 |x_j-x_{j-1}|\leq L/N<\rho.
\]
For \(1\leq j\leq N\), introduce the cylinders
\begin{align*}
 Q_j&:=(\ell_j-4\rho^2,\ell_j]\times B_{2\rho}(x_j),\\
 Q_j^-&:=(\ell_j-3\rho^2,\ell_j-2\rho^2)
                      \times B_\rho(x_j),\\
 Q_j^+&:=(\ell_j-\rho^2,\ell_j]\times B_\rho(x_j).
\end{align*}
Their closures satisfy
\begin{equation}\label{eq:section9-cylinder-containment}
 \overline{Q_j}\subset\mathfrak C_{a,K,O}.
\end{equation}
Indeed, \(x_j\in E\) and \(2\rho<d_0/2\) give the spatial
containment.  For the time coordinate, using \(N\geq2\),
\[
 \ell_j-4\rho^2
 \geq\ell_1-4\rho^2
 =\frac a2-\frac{3a}{10N}
 \geq\frac{7a}{20}>\frac a4,
 \qquad \ell_j\leq a.
\]
The time and space increments also imply
\[
 z_{j-1}\in Q_j^-,\qquad z_j\in Q_j^+.
\]
In particular, each step respects the earlier-to-later time order
in Proposition~\ref{prop:controlled-local-Harnack}.

Let \(u\) be a continuous viscosity solution of either
\eqref{eq:controlled-remaining-HJB} or
\eqref{eq:controlled-conjugate-HJB}, with any fixed terminal calendar
time, on a neighbourhood of \(\mathfrak C_{a,K,O}\).
Suppose that \(u\geq0\) on \(\mathfrak C_{a,K,O}\).
By \eqref{eq:section9-cylinder-containment}, the local Harnack
inequality applies in every \(Q_j\), and gives
\[
 u(z_{j-1})
 \leq\sup_{Q_j^-}u
 \leq C_H\inf_{Q_j^+}u
 \leq C_Hu(z_j).
\]
All these cylinders have the same radius and relative geometry.
Thus one may use the same constant \(C_H\geq1\), depending only on
\(d,\underline\lambda,\overline\lambda,\rho B_O\), for every
\(j\), every target \(x\in K\), and either equation.
The inclusion of the top point \(z_j\) is justified by the
continuity stipulated in the local Harnack proposition.

Iterating the preceding inequalities yields
\[
 u(a/2,x_*)\leq C_H^N u(a,x),\qquad x\in K.
\]
The paths may depend on \(x\), but \(N\), \(\rho\), and \(C_H\)
do not.  Define
\[
 z_*:=(a/2,x_*),\qquad
 H_{a,K,O}:=\max\{2,C_H^N\}.
\]
Since \(u\) is nonnegative on the target slice, taking the infimum
in \(x\) gives the common chain estimate
\begin{equation}\label{eq:section9-common-Harnack-chain}
 u(z_*)\leq H_{a,K,O}\inf_{x\in K}u(a,x).
\end{equation}
No nonnegativity assumption outside \(\mathfrak C_{a,K,O}\) has
been used. 
Note the number \(2\) in the definition of \(H_{a,K,O}\) is included
only to ensure that \(H_{a,K,O}>1\).

\smallskip
\noindent\emph{Step 3: application to the value and its conjugate.}
Fix \(s\in\R\) and \(\varphi\in\BUC_b(\R^d)\).
Theorem~\ref{thm:controlled-positive-time} and
Proposition~\ref{prop:controlled-value-PDE}, with terminal calendar
time \(t_*=s+a\), show that \(v_s\) is continuous and solves
\eqref{eq:controlled-remaining-HJB} for \(\ell>0\).
Set
\[
 m_s:=\min_{\mathfrak C_{a,K,O}}v_s,
 \qquad
 M_s:=\max_{\mathfrak C_{a,K,O}}v_s.
\]
The extrema are finite and attained because the comparison region
is compact and contained in the positive-time domain.

The functions
\[
 u_s:=v_s-m_s,\qquad w_s:=M_s-v_s
\]
are nonnegative on \(\mathfrak C_{a,K,O}\).
Since the value equation has no dependence on the undifferentiated
solution, subtracting a constant shows that \(u_s\) solves
\eqref{eq:controlled-remaining-HJB}.
By the conjugate-equation assertion of
Proposition~\ref{prop:controlled-value-PDE}, \(w_s\) solves
\eqref{eq:controlled-conjugate-HJB}.
Both functions are defined and solve their respective equations
on a neighbourhood of every cylinder closure in Step~2.
It is not necessary that they solve the same nonlinear equation:
Proposition~\ref{prop:controlled-local-Harnack} supplies the same
Harnack bound for the two conjugate equations.

Applying \eqref{eq:section9-common-Harnack-chain} to each function
gives
\begin{align}
 v_s(z_*)-m_s
 &\leq H_{a,K,O}
       \left(\inf_{x\in K}v_s(a,x)-m_s\right),
 \label{eq:section9-lower-gap}\\
 M_s-v_s(z_*)
 &\leq H_{a,K,O}
       \left(M_s-\sup_{x\in K}v_s(a,x)\right).
 \label{eq:section9-upper-gap}
\end{align}

\smallskip
\noindent\emph{Step 4: oscillation contraction.}
Write
\[
 I_s:=\inf_{x\in K}v_s(a,x),\qquad
 J_s:=\sup_{x\in K}v_s(a,x).
\]
Adding \eqref{eq:section9-lower-gap} and
\eqref{eq:section9-upper-gap} yields
\begin{align*}
 M_s-m_s
 &\leq H_{a,K,O}\bigl[(I_s-m_s)+(M_s-J_s)\bigr]\\
 &=H_{a,K,O}\bigl[(M_s-m_s)-(J_s-I_s)\bigr].
\end{align*}
Rearranging, we obtain
\[
 \osc_Kv_s(a,\cdot)=J_s-I_s
 \leq\left(1-H_{a,K,O}^{-1}\right)(M_s-m_s).
\]
Since
\(
v_s(a,x)=\bS_{s,s+a}\varphi(x),\qquad
 M_s-m_s=\osc_{\mathfrak C_{a,K,O}}v_s,
\)
this is \eqref{eq:controlled-local-oscillation}, with
\(
 \vartheta_{a,K,O}:=1-H_{a,K,O}^{-1}\in(0,1).
\)
The calculation also covers the case \(M_s=m_s\), since no division
by the oscillation is made.

Finally, the comparison region, the radius \(\rho\), and the number
of steps \(N\) were fixed using only \(a,K,O\).
The local Harnack constant is uniform in calendar time by
Assumption~\ref{ass:uniform-parabolic} and
Proposition~\ref{prop:controlled-local-Harnack}.
Consequently, \(H_{a,K,O}\) and \(\vartheta_{a,K,O}\) have the
asserted dependence and are independent of \(s\) and \(\varphi\).
\end{proof}

\subsection{Verification for the controlled nonstationary diffusion}

{
For a continuous weight \(V:\R^d\to[0,\infty)\), define
\[
 \mathscr C_V(\R^d)
 :=\left\{\varphi\in C(\R^d):
       \sup_{x\in\R^d}\frac{\abs{\varphi(x)}}{1+V(x)}<\infty\right\}.
\]
Thus \(\mathscr C_V(\R^d)\) is the continuous subspace of
\(\mathscr B_V(\R^d)\).
}

The abstract theorem is now applied with
\(E=\R^d\), \(\mathscr X_V=\mathscr C_V(\R^d)\), and with the
all-policy moment bound
\begin{equation}\label{eq:all-policy-envelope}
 \mathsf M_{s,t}V(x)
 :=\sup_{\alpha\in\cA}\sup_{P\in\cP_{\cQ}}
 E^P[V(X_t^{s,x,\alpha})].
\end{equation}
This quantity dominates the lower value \(\bS_{s,t}V\) and controls terminal
perturbations uniformly over all admissible controls and covariance laws.

\begin{lemma}
\label{lem:controlled-symmetric-bound}
For terminal functions for which the expectations are finite,
\begin{equation}\label{eq:controlled-terminal-perturbation}
 \abs{\bS_{s,t}\varphi(x)-\bS_{s,t}\psi(x)}
 \leq\sup_{\alpha\in\cA}\sup_{P\in\cP_{\cQ}}
 E^P\!\left[\abs{\varphi-\psi}(X_t^{s,x,\alpha})\right].
\end{equation}
Consequently, if \(\abs{\varphi-\psi}\leq c+dV\), then
\begin{equation}\label{eq:controlled-model-domination}
 \abs{\bS_{s,t}\varphi(x)-\bS_{s,t}\psi(x)}
 \leq c+d\,\mathsf M_{s,t}V(x).
\end{equation}
In particular, since \(\bS_{s,t}0=0\), if
\(\abs{\varphi}\leq1+\beta V\), then
\begin{equation}\label{eq:controlled-symmetric-bound}
 \abs{\bS_{s,t}\varphi(x)}
 \leq1+\beta\mathsf M_{s,t}V(x).
\end{equation}
Thus Assumption~\ref{ass:abstract-terminal-domination} holds with
\(m_{s,t}=\mathsf M_{s,t}V\).
\end{lemma}

\begin{proof}
For each fixed \(\alpha\), the absolute difference of the two suprema over
\(P\) is bounded by the supremum of the absolute difference of the terminal
payoffs.  The absolute difference of the two infima over \(\alpha\) is then
bounded by the supremum over \(\alpha\) of these quantities.  This proves
\eqref{eq:controlled-terminal-perturbation}; the remaining assertions follow
immediately.
\end{proof}

From now on take
\begin{equation}\label{eq:quadratic-V}
 V(x):=1+\abs{x}^2.
\end{equation}

{\renewcommand{\theassumption}{D}
\begin{assumption}
\label{ass:controlled-dissipativity}
There are $c_0>0$ and $c_1\geq0$ such that, for every
$(q,x,a,Q)\in\R\times\R^d\times A\times\cQ$,
\begin{equation}\label{eq:controlled-dissipativity}
 2x\cdot b_{a,Q}(q,x)
 +\tr({c_{a,Q}}(q,x))
 \leq-c_0\abs{x}^2+c_1.
\end{equation}
\end{assumption}
}
\begin{proposition}
\label{prop:controlled-Lyapunov-verification}
Under Assumptions~\ref{ass:uniform-parabolic} and
\ref{ass:controlled-dissipativity},
\begin{equation}\label{eq:all-policy-Lyapunov}
 \mathsf M_{s,t}V(x)
 \leq \exp(-c_0(t-s))V(x)
 +\kappa\bigl(1-\exp(-c_0(t-s))\bigr),
 \qquad
 \kappa:=\frac{c_0+c_1}{c_0}.
\end{equation}
Hence the abstract Lyapunov
Assumption~\ref{ass:controlled-Lyapunov} holds, with
\(m_{s,t}=\mathsf M_{s,t}V\), for every block length \(a>0\) and with
\begin{equation}\label{eq:explicit-Lyapunov-constants}
 \gamma=\exp(-c_0a),
 \qquad K=(1-\gamma)\kappa,
 \qquad \Gamma_a=1,
 \qquad K_a=\kappa.
\end{equation}
\end{proposition}

\begin{proof}
Fix \(\alpha\in\cA\) and \(P\in\cP_{\cQ}\), and write
\(d\langle B\rangle_u=Q_u^P\,du\), with
\(Q_u^P=\widehat Q_u\in\cQ\), \(du\times P\)-almost everywhere.
Apply It\^o's formula to the elapsed-time process
\(u\mapsto X_{s+u}^{s,x,\alpha}\) in \eqref{eq:controlled-SDE}.
By \eqref{eq:controlled-dissipativity},
\[
 \frac{d}{du}E^P[V(X_{s+u}^{s,x,\alpha})]
 \leq-c_0E^P[V(X_{s+u}^{s,x,\alpha})]+c_0+c_1
 \quad\text{for a.e. }u\geq0.
\]
Gronwall's inequality, evaluated at \(u=t-s\), proves
\eqref{eq:all-policy-Lyapunov} for this \((\alpha,P)\).
The estimate is independent of both, so taking the two suprema proves
it for \(\mathsf M_{s,t}V(x)\).  Formula
\eqref{eq:explicit-Lyapunov-constants} follows immediately.
\end{proof}

\begin{lemma}
\label{lem:controlled-weighted-extension}
Under Assumptions~\ref{ass:uniform-parabolic} and
\ref{ass:controlled-dissipativity}, the expectation formula
\eqref{eq:controlled-evolution} defines $\bS_{s,t}$ on
$\mathscr C_V(\R^d)$ with \(V(x)=1+\abs{x}^2\).  This class is preserved by the evolution, and the
deterministic DPP extends to it.  Moreover, for every compact $K$, finite
horizon $H$ and $p>2$,
\begin{equation}\label{eq:controlled-uniform-p-moment}
 \sup_{\substack{s\in\R,\,x\in K,\,s\leq r\leq s+H\\
                  \alpha\in\cA,\,P\in\cP_{\cQ}}}
 E^P[\abs{X_r^{s,x,\alpha}}^p]<\infty,
\end{equation}
and bounded approximations of $\varphi\in\mathscr C_V$ converge after
positive time locally uniformly, uniformly in the initial calendar time.
\end{lemma}

\begin{proof}
The global Lipschitz and linear-growth bounds, BDG and Gronwall give
\eqref{eq:controlled-uniform-p-moment}.  Since
$\abs{\varphi(x)}\leq C(1+V(x))$, this estimate gives uniform integrability
of the terminal costs on compact sets of initial states.  Clip $\varphi$ at
levels $\pm N$, approximate the clipped continuous function on large balls
by bounded uniformly continuous functions, and extend those approximations
to $\R^d$.  First let the approximation error tend to zero, then the ball
radius and $N$ tend to infinity.  Equation
\eqref{eq:controlled-uniform-p-moment} shows that the corresponding values
converge locally uniformly, uniformly in $s$.  Passing to the limit in the
bounded DPP proves the evolution identity on $\mathscr C_V$.  The same
argument, together with Theorem~\ref{thm:controlled-positive-time}, proves
continuity after every positive elapsed time and the asserted invariance of
the weighted domain.
\end{proof}

\begin{proposition}
\label{prop:controlled-difference-LOC}
Suppose Assumptions~\ref{ass:uniform-parabolic} and
\ref{ass:controlled-dissipativity} hold, and let
\(V(x)=1+\abs{x}^2\).  Fix \(a>0\) and \(R<\infty\) such that
\(C_R:=\{V\leq R\}\) is nonempty and compact.  Choose a bounded connected
open set \(O_R\) with \(C_R\Subset O_R\), and let
\(\mathfrak C_{a,C_R,O_R}\) and
\(\vartheta_{a,C_R,O_R}\in(0,1)\) be obtained from
Proposition~\ref{prop:controlled-local-oscillation}.  For \(s\in\R\), set
\begin{equation}\label{eq:controlled-calendar-Harnack-set}
 \mathfrak Q_{s,R}
 :=\{(r,x):(s+a-r,x)\in\mathfrak C_{a,C_R,O_R}\}.
\end{equation}
Then, for every \(\varphi,\psi\in\mathscr C_V(\R^d)\),
\begin{align}
 \osc_{C_R}\!\left(
   d_s^{\varphi,\psi}(s,\cdot)\right)\leq\vartheta_R
 \osc_{\mathfrak Q_{s,R}}\left(d_s^{\varphi,\psi}(\cdot,\cdot)\right),
 \qquad
 \vartheta_R:=\vartheta_{a,C_R,O_R}<1,
 \label{eq:controlled-difference-LOC}
\end{align}
where
\[
 d_s^{\varphi,\psi}(r,x)
 :=\bS_{r,s+a}\varphi(x)-\bS_{r,s+a}\psi(x).
\]
Moreover,
\begin{equation}\label{eq:controlled-explicit-LR}
 L_R:=\sup_{s\in\R}\sup_{(r,x)\in\mathfrak Q_{s,R}}
       \mathsf M_{r,s+a}V(x)<\infty.
\end{equation}
Thus Assumption~\ref{ass:controlled-LOC} holds for the controlled evolution.
\end{proposition}

\begin{proof}
First take bounded uniformly continuous terminal data and fix
\(t_*=s+a\).  Put
\[
 u(\ell,x):=\bS_{t_*-\ell,t_*}\varphi(x),\qquad
 v(\ell,x):=\bS_{t_*-\ell,t_*}\psi(x),
\]
and \(w:=u-v\).  Proposition~\ref{prop:controlled-value-PDE} shows that
\(u\) and \(v\) solve the same uniformly parabolic Bellman--Isaacs
equation.  
Apply the parabolic Crandall--Ishii theorem of sums
\cite[Theorem~8.3]{CIL1992} to \(u\) and \(-v\), using
independent spatial variables \(x\) and \(y\) and the penalty
\(\abs{x-y}^2/(2\varepsilon)\).
Combining the resulting viscosity inequalities with the structural
increment bounds for the Bellman--Isaacs operator
\(H_0\) defined in \eqref{eq:lower-Isaacs-operator} and the spatial Lipschitz estimates
gives, in the viscosity sense on \((0,\infty)\times O\),
\begin{align*}
 \mathcal M^-_{\underline\lambda/2,\overline\lambda/2}(D^2w)
 -B_O\abs{Dw}
 &\leq\partial_\ell w,\\
 \partial_\ell w
 &\leq
 \mathcal M^+_{\underline\lambda/2,\overline\lambda/2}(D^2w)
 +B_O\abs{Dw},
\end{align*}
for every bounded open spatial set \(O\).
Here, for \(X\in\mathbb S^d\),
\[
 \mathcal M^-_{\lambda,\Lambda}(X)
 :=\inf_{\lambda I_d\leq A\leq\Lambda I_d}\tr(AX),
 \qquad
 \mathcal M^+_{\lambda,\Lambda}(X)
 :=\sup_{\lambda I_d\leq A\leq\Lambda I_d}\tr(AX)
\]
are the Pucci extremal operators.  Thus \(w\) belongs to the same
homogeneous parabolic Pucci class as in
Proposition~\ref{prop:controlled-local-Harnack}, with ellipticity constants
\(\underline\lambda/2\) and \(\overline\lambda/2\) and local drift
bound \(B_O\).
The structural increment estimate behind this argument is, for
\(q\in\R\), \(x\in O\), \(p_1,p_2\in\R^d\), and
\(X_1,X_2\in\mathbb S^d\),
\begin{align*}
 &\mathcal M^-_{\underline\lambda/2,\overline\lambda/2}(X_1-X_2)
       -B_O\abs{p_1-p_2}\\
 &\qquad\leq H_0(q,x,p_1,X_1)-H_0(q,x,p_2,X_2)\\
 &\qquad\leq
   \mathcal M^+_{\underline\lambda/2,\overline\lambda/2}(X_1-X_2)
       +B_O\abs{p_1-p_2}.
\end{align*}
It follows by comparing the two generators at the same action and
covariance, then taking the two optimisations.  The state-Lipschitz
bounds control the spatial doubling errors in the viscosity theorem of
sums.  No convexity or concavity of \(H_0\) is used.

The parabolic Harnack theorem used in
Proposition~\ref{prop:controlled-local-Harnack} applies to nonnegative
functions in this Pucci class.  Hence, on every cylinder on which
\(m\leq w\leq M\), both \(w-m\) and \(M-w\) satisfy the same nonnegative
Harnack estimate.  Repeating the finite forward-chain construction in the
proof of Proposition~\ref{prop:controlled-local-oscillation}, with \(w\) in
place of the single value, gives
\eqref{eq:controlled-difference-LOC} with exactly the same chain constant.
Reflection from remaining time to calendar time is precisely
\eqref{eq:controlled-calendar-Harnack-set}.

For \(\varphi,\psi\in\mathscr C_V(\R^d)\), use the bounded uniformly
continuous approximations from Lemma~\ref{lem:controlled-weighted-extension}.
Their controlled values converge locally uniformly after every positive
time, uniformly in the initial calendar time, so the oscillation inequality
passes to the limit.

Finally, the spatial component of
\(\mathfrak C_{a,C_R,O_R}\) is contained in \(O_R\).  With
\[
 \widehat R_R:=\sup_{x\in O_R}V(x)<\infty,
\]
the short-time Lyapunov estimate gives
\[
 L_R\leq\Gamma_a\widehat R_R+K_a<\infty.
\]
This proves the proposition.
\end{proof}
\begin{theorem}
\label{thm:controlled-diffusion-Harris}
Under Assumptions~\ref{ass:uniform-parabolic} and
\ref{ass:controlled-dissipativity}, there exist a block length \(a>0\),
\(\beta>0\), \(q\in(0,1)\), \(C_0<\infty\) and \(\lambda>0\) such that,
for every \(\varphi,\psi\in\mathscr C_V(\R^d)\),
\begin{align}
 [\bS_{s,s+a}\varphi-\bS_{s,s+a}\psi]_{\beta,V}
 &\leq q[\varphi-\psi]_{\beta,V},
 \label{eq:controlled-diffusion-block-contraction-difference}\\
 [\bS_{s,t}\varphi-\bS_{s,t}\psi]_{\beta,V}
 &\leq C_0e^{-\lambda(t-s)}[\varphi-\psi]_{\beta,V},
 \qquad s\leq t.
 \label{eq:controlled-diffusion-exponential-contraction-difference}
\end{align}
The constants are independent of the initial calendar time \(s\).  In
particular, since \(\bS_{s,t}0=0\),
\begin{align}
 [\bS_{s,s+a}\varphi]_{\beta,V}
 &\leq q[\varphi]_{\beta,V},
 \label{eq:controlled-diffusion-block-contraction}\\
 [\bS_{s,t}\varphi]_{\beta,V}
 &\leq C_0e^{-\lambda(t-s)}[\varphi]_{\beta,V}.
 \label{eq:controlled-diffusion-exponential-contraction}
\end{align}
\end{theorem}

\begin{proof}
Fix any \(a>0\).  Proposition~\ref{prop:controlled-Lyapunov-verification}
verifies Assumption~\ref{ass:controlled-Lyapunov} for
\(m_{s,t}=\mathsf M_{s,t}V\), with
\begin{equation*}
 R>2\kappa=\frac{2(c_0+c_1)}{c_0}
\end{equation*}
ensuring \eqref{eq:controlled-R-choice}, because
\(K=(1-\gamma)\kappa\).  Lemma~\ref{lem:controlled-symmetric-bound}
verifies Assumption~\ref{ass:abstract-terminal-domination}.
Proposition~\ref{prop:controlled-difference-LOC} verifies
Assumption~\ref{ass:controlled-LOC}, including the finiteness of \(L_R\).
The deterministic DPP gives the evolution property on
\(\mathscr C_V(\R^d)\) by Lemma~\ref{lem:controlled-weighted-extension},
and the infimum--supremum representation gives translation covariance.
All hypotheses of Theorem~\ref{thm:controlled-weighted-Harris} are therefore
satisfied.  Applying that theorem proves the two difference estimates with
the explicit constants in \eqref{eq:controlled-contraction-constants}.
Taking \(\psi=0\) proves the last two displays.
\end{proof}

\begin{remark}
For an upper bound on \(\bS_{s,t}V\), it would be enough to exhibit one
stabilising feedback rule.  The abstract theorem requires more:
a bound controlling both signs of every centred terminal perturbation,
uniformly over all admissible controller--nature pairs.  The all-policy
moment bound \(\mathsf M_{s,t}V\) provides this control by
Lemma~\ref{lem:controlled-symmetric-bound}. Uniform dissipativity over every
\(a\in A\) and \(Q\in\cQ\) is therefore a clean sufficient condition for
contraction on the full two-sided quotient space \(\mathscr C_V/\R\).
Weaker policy-dependent Lyapunov hypotheses would require a different
one-sided framework or a restriction of the admissible controls.
\end{remark}

\subsection{Large-time convergence and entrance nonlinear expectations}

For a genuinely nonstationary evolution one should not, in general, expect
convergence to a single time-independent functional.  The correct object is
an entrance family indexed by the terminal calendar time.  Fix a reference
point $x_\star\in\R^d$ and put
\begin{equation}\label{eq:controlled-reference-moment}
 \overline V_\star:=V(x_\star)\vee\kappa,
 \qquad \kappa=\frac{c_0+c_1}{c_0}.
\end{equation}
By \eqref{eq:all-policy-Lyapunov},
\begin{equation}\label{eq:controlled-uniform-reference-moment}
 \sup_{r\leq s}\mathsf M_{r,s}V(x_\star)
 \leq\overline V_\star.
\end{equation}

\begin{theorem}
\label{thm:controlled-large-time}
Suppose Assumptions~\ref{ass:uniform-parabolic} and
\ref{ass:controlled-dissipativity} hold.  For every $t\in\R$ and
$\varphi\in\mathscr C_V(\R^d)$, the pullback limit
\begin{equation}\label{eq:controlled-entrance-expectation}
 \mathcal E_t(\varphi)
 :=\lim_{r\to-\infty}\bS_{r,t}\varphi(x_\star)
\end{equation}
exists and is independent of $x_\star$.  Each $\mathcal E_t$ is monotone,
constant preserving, translation invariant and positively homogeneous.
The family satisfies the entrance relation
\begin{equation}\label{eq:controlled-entrance-relation}
 \mathcal E_s(\bS_{s,t}\varphi)=\mathcal E_t(\varphi),
 \qquad s\leq t.
\end{equation}
Moreover, with the constants in
Theorem~\ref{thm:controlled-diffusion-Harris},
\begin{equation}\label{eq:controlled-large-time-rate}
 \left|\bS_{s,t}\varphi(x)-\mathcal E_t(\varphi)\right|
 \leq C_0e^{-\lambda(t-s)}[\varphi]_{\beta,V}
 \left(2+\beta V(x)+\beta\overline V_\star\right)
\end{equation}
for all $s\leq t$ and $x\in\R^d$.
{In particular, for fixed $s$,
\(\bS_{s,t}\varphi(x)-\mathcal E_t(\varphi)\to0\) as $t\to\infty$,
locally uniformly in $x$ and exponentially fast.}

{Finally, $(\mathcal E_t)_{t\in\R}$ is unique among entrance families
$(\widetilde{\mathcal E}_t)_{t\in\R}$ with these properties satisfying}
\begin{equation}\label{eq:controlled-tempered-entrance-family}
 \sup_{t\in\R}\sup_{\norm{\psi}_{\beta,V}\leq1}
 \left|\widetilde{\mathcal E}_t(\psi)\right|<\infty.
\end{equation}
{
Here and below, ``tempered'' refers to the uniform weighted
bound \eqref{eq:controlled-tempered-entrance-family}.
}
\end{theorem}

\begin{proof}
Fix $t$ and set $g_r:=\bS_{r,t}\varphi$.  If $r'<r<t$, the evolution
property and translation invariance give
\[
 \bS_{r',t}\varphi(x_\star)-\bS_{r,t}\varphi(x_\star)
 =\bS_{r',r}\bigl(g_r-g_r(x_\star)\bigr)(x_\star).
\]
By the definition of the weighted seminorm,
\[
 |g_r(y)-g_r(x_\star)|
 \leq[g_r]_{\beta,V}
 \bigl(2+\beta V(y)+\beta V(x_\star)\bigr).
\]
The argument of Lemma~\ref{lem:controlled-symmetric-bound}, now with
the constant \(2+\beta V(x_\star)\) in place of \(1\), gives the first
inequality below.  The second follows from
Theorem~\ref{thm:controlled-diffusion-Harris}, applied with initial
time \(r\) and terminal time \(t\), together with the uniform moment
bound \eqref{eq:controlled-uniform-reference-moment}:
\begin{align*}
 &\left|\bS_{r',t}\varphi(x_\star)
        -\bS_{r,t}\varphi(x_\star)\right|\\
 &\qquad\leq[g_r]_{\beta,V}
 \left(2+\beta V(x_\star)
       +\beta\mathsf M_{r',r}V(x_\star)\right)\\
 &\qquad\leq C_0e^{-\lambda(t-r)}[\varphi]_{\beta,V}
 \left(2+\beta V(x_\star)+\beta\overline V_\star\right).
\end{align*}
The last expression tends to zero as $r\to-\infty$, uniformly in $r'<r$.
Hence the net in \eqref{eq:controlled-entrance-expectation} is Cauchy.

For another reference point $y_\star$,
\[
 |\bS_{r,t}\varphi(x_\star)
   -\bS_{r,t}\varphi(y_\star)|
 \leq C_0e^{-\lambda(t-r)}[\varphi]_{\beta,V}
 \bigl(2+\beta V(x_\star)+\beta V(y_\star)\bigr),
\]
which tends to zero.  Thus the limit is independent of the reference
point.  Monotonicity, preservation of constants, translation invariance and
positive homogeneity pass to the limit from $\bS_{r,t}$.

Since $\bS_{s,t}$ preserves $\mathscr C_V$, the evolution property gives
\begin{align*}
 \mathcal E_s(\bS_{s,t}\varphi)
 &=\lim_{r\to-\infty}
   \bS_{r,s}(\bS_{s,t}\varphi)(x_\star)\\
 &=\lim_{r\to-\infty}\bS_{r,t}\varphi(x_\star)
 =\mathcal E_t(\varphi),
\end{align*}
which proves \eqref{eq:controlled-entrance-relation}.

To obtain the rate, set $g:=\bS_{s,t}\varphi$.  For $r<s$, the same
weighted estimate gives
\begin{align*}
 |\bS_{r,t}\varphi(x_\star)-g(x)|
 &=|\bS_{r,s}g(x_\star)-g(x)|\\
 &\leq[g]_{\beta,V}
 \left(2+\beta V(x)
       +\beta\mathsf M_{r,s}V(x_\star)\right).
\end{align*}
Let $r\to-\infty$, use
\eqref{eq:controlled-uniform-reference-moment}, and then apply
\eqref{eq:controlled-diffusion-exponential-contraction} to $g$.  This is
\eqref{eq:controlled-large-time-rate}.

It remains to prove the uniqueness.  Let
\((\widetilde{\mathcal E}_t)_{t\in\R}\) be another entrance family
satisfying the stated conditions and the entrance relation.  For $s<t$, put
the same $g$ as above.  By the quotient identity, for every $\varepsilon>0$
there is $c\in\R$ such that
\[
 \norm{g-c}_{\beta,V}\leq[g]_{\beta,V}+\varepsilon.
\]
Translation invariance, the uniform bound
\eqref{eq:controlled-tempered-entrance-family}, and the corresponding bound
for $\mathcal E_s$, which follows from
\eqref{eq:controlled-uniform-reference-moment}, imply
\[
 |\widetilde{\mathcal E}_t(\varphi)-\mathcal E_t(\varphi)|
 =|\widetilde{\mathcal E}_s(g)-\mathcal E_s(g)|
 \leq C\bigl([g]_{\beta,V}+\varepsilon\bigr),
\]
where $C$ is independent of $s$.  First let $\varepsilon\downarrow0$ and
then $s\to-\infty$ {in the term} \([g]_{\beta,V}\).  The exponential contraction of $[g]_{\beta,V}$ yields
\[
 \widetilde{\mathcal E}_t(\varphi)=\mathcal E_t(\varphi).
\]
Hence the entrance family \((\mathcal E_t)_{t\in\R}\) is unique in the
stated tempered class.
\end{proof}

\begin{remark}
\label{rem:entrance-risk-functionals}
Theorem~\ref{thm:controlled-large-time} has a direct long-run risk
interpretation.  The family \((\mathcal E_t)_{t\in\R}\) is the unique
tempered entrance monetary risk functional generated by optimal robust
control, and \eqref{eq:controlled-entrance-relation} is its compatibility
with the controlled risk evolution. 
In the uncontrolled case, the corresponding evolution is sublinear,
and under the standing hypotheses the same construction yields an
entrance family of coherent risk evaluations.
The same conclusion holds for any fixed-policy state-space evolution
satisfying the evolution property and the same terminal-perturbation
domination, Lyapunov and local oscillation hypotheses.
In the optimally controlled problem, subadditivity may fail, but
cash translation remains exact.
 The weighted
quotient seminorm is therefore natural from this viewpoint: it factors out
deterministic cash shifts, while
\eqref{eq:controlled-large-time-rate} gives exponential loss of the initial
state for the economically nontrivial, state-dependent part of the risk
profile.  Autonomous coefficients lead to an invariant long-run risk
functional, while periodic coefficients lead to a periodic family of such
functionals.
\end{remark}

\begin{corollary}
\label{cor:controlled-periodic-autonomous}
{
Under the assumptions and notation of
Theorem~\ref{thm:controlled-large-time}, fix
\(\varphi\in\mathscr C_V(\R^d)\).  Suppose first that the controlled
coefficients are \(T\)-periodic in calendar time.  Then
}
\begin{equation}\label{eq:controlled-periodic-evolution}
 \bS_{s+T,t+T}=\bS_{s,t},
 \qquad
 \mathcal E_{t+T}=\mathcal E_t,
\end{equation}
and, for every $n\geq1$,
\begin{equation}\label{eq:controlled-periodic-convergence}
 \left|\bS_{s,s+nT}\varphi(x)-\mathcal E_s(\varphi)\right|
 \leq C_0e^{-\lambda nT}[\varphi]_{\beta,V}
 \left(2+\beta V(x)+\beta\overline V_\star\right).
\end{equation}

If the coefficients are autonomous, then
$\bS_{s,t}=\mathsf T_{t-s}$ and the entrance family is independent of
calendar time.  Writing its common value as $\mathcal E$, one has
\begin{equation}\label{eq:controlled-invariant-expectation}
 \mathcal E(\mathsf T_t\varphi)=\mathcal E(\varphi)
\end{equation}
and
\begin{equation}\label{eq:controlled-autonomous-convergence}
 |\mathsf T_t\varphi(x)-\mathcal E(\varphi)|
 \leq C_0e^{-\lambda t}[\varphi]_{\beta,V}
 \left(2+\beta V(x)+\beta\overline V_\star\right).
\end{equation}
\end{corollary}

\begin{proof} 
No periodicity condition on the individual feedback rules is required:
the admissible class \(\cA\) is defined in elapsed time and is the same
for every initial calendar time.
For any fixed \(\alpha\in\cA\) and \(P\in\cP_{\cQ}\), periodicity
of the coefficients and pathwise uniqueness imply
\[
 X_{s+T+u}^{s+T,x,\alpha}
 =
 X_{s+u}^{s,x,\alpha},
 \qquad u\geq0,\quad P\text{-a.s.}
\]
Indeed, the two processes solve the same elapsed-time equation,
driven by the same canonical process and using the same feedback rule.
Taking first the supremum over \(P\) and then the infimum over
\(\alpha\) proves
\(\bS_{s+T,t+T}=\bS_{s,t}\).  Shift
both endpoints in \eqref{eq:controlled-entrance-expectation} by $T$ to
obtain $\mathcal E_{t+T}=\mathcal E_t$.  Formula
\eqref{eq:controlled-periodic-convergence} is
\eqref{eq:controlled-large-time-rate} with $t=s+nT$.  In the autonomous
case, the same shift argument makes $\mathcal E_t$ constant in $t$;
\eqref{eq:controlled-entrance-relation} becomes invariance, and
\eqref{eq:controlled-autonomous-convergence} follows from
\eqref{eq:controlled-large-time-rate}.
\end{proof}

\section{Controlled G-diffusion with running cost and eigenpair}
\subsection{Running costs and the additive Bellman evolution}

The zero-running-cost controlled evolution now has deterministic and stopped
DPPs, positive-time regularity, a local Harnack inequality, local
oscillation contraction, exponential weighted Harris contraction and a
unique tempered entrance expectation with exponential attraction. 
In this section, we consider the problem with a running cost.
Although the resulting Bellman evolution is no
longer constant preserving or positively homogeneous, differences of two
values with the same running cost retain a weighted span contraction.  In the
autonomous case this yields an additive eigenvalue, a corrector and convergence
of the long-time average cost.

Let
\[
 f:\R\times\R^d\times A\longrightarrow\R
\]
be a running cost.  Nature still enters the cost through the law of the state;
the direct dependence of $f$ on the covariance scenario is omitted in this
first formulation.

{\renewcommand{\theassumption}{R}
\begin{assumption}
\label{ass:running-cost}
The function $f$ is jointly continuous and, for some $c_f<\infty$,
\begin{equation}\label{eq:running-cost-growth}
 |f(q,x,a)|\leq c_f(1+V(x)),
 \qquad (q,x,a)\in\R\times\R^d\times A.
\end{equation}
On every compact calendar-time interval and bounded state set, its modulus of
continuity can be chosen uniformly in $a\in A$.
\end{assumption}
}

Throughout this section, Assumptions~\ref{ass:uniform-parabolic},
\ref{ass:controlled-dissipativity} and \ref{ass:running-cost} are in
force, and terminal functions belong to \(\mathscr C_V(\R^d)\), with
\(V(x)=1+\abs{x}^2\).
By truncation and the uniform moment estimate in
Lemma~\ref{lem:controlled-flow}, the measurability conclusion of
Lemma~\ref{lem:controlled-canonical-state-map} continues to hold for the
present running-cost class.

For $s\leq t$, define the accumulated cost
\begin{equation}\label{eq:accumulated-running-cost}
 \begin{aligned}
 C_{s,t}^{s,x,\alpha}
 &:=\int_s^t f(r,X_r^{s,x,\alpha},\alpha_r^{[s]})\,dr\\
 &=\int_0^{t-s}f(s+u,X_{s+u}^{s,x,\alpha},\alpha_u)\,du
 \end{aligned}
\end{equation}
and the running-cost value operator
\begin{equation}\label{eq:running-cost-operator}
 \bU_{s,t}^{f}\varphi(x)
 :=\inf_{\alpha\in\cA}\sup_{P\in\cP_{\cQ}}
 E^P\!\left[
  C_{s,t}^{s,x,\alpha}
  +\varphi(X_t^{s,x,\alpha})
 \right].
\end{equation}
Assumption~\ref{ass:running-cost}, the moment estimates and the all-policy
Lyapunov bound make \eqref{eq:running-cost-operator} finite on
$\mathscr C_V(\R^d)$.  Notice that $\bU_{s,t}^{0}=\bS_{s,t}$.

\begin{proposition}
\label{prop:running-cost-DPP}
Under Assumptions~\ref{ass:uniform-parabolic},
\ref{ass:controlled-dissipativity} and \ref{ass:running-cost},
let \(V(x)=1+\abs{x}^2\).  Then, for \(s\leq r\leq t\) and
\(\varphi\in\mathscr C_V(\R^d)\),
\begin{align}
 \bU_{s,t}^{f}\varphi(x)
 =\inf_{\alpha\in\cA}\sup_{P\in\cP_{\cQ}}
 E^P\!\left[
  \int_s^r f(v,X_v^{s,x,\alpha},\alpha_v^{[s]})\,dv
  +\bU_{r,t}^{f}\varphi(X_r^{s,x,\alpha})
 \right].
 \label{eq:running-cost-deterministic-DPP}
\end{align}
Equivalently,
\begin{equation}\label{eq:running-cost-evolution}
 \bU_{s,t}^{f}=\bU_{s,r}^{f}\bU_{r,t}^{f}.
\end{equation}
If $\tau$ is a finite elapsed-time canonical stopping time with
$0\leq\tau\leq t-s$,
then
\begin{align}
 \bU_{s,t}^{f}\varphi(x)
 =\inf_{\alpha\in\cA}\sup_{P\in\cP_{\cQ}}
 E^P\!\left[
  \int_s^{s+\tau}f(v,X_v^{s,x,\alpha},\alpha_v^{[s]})\,dv
  +\bU_{s+\tau,t}^{f}\varphi(X_{s+\tau}^{s,x,\alpha})
 \right].
 \label{eq:running-cost-stopped-DPP}
\end{align}
The same formula holds for the feedback-compatible stopping rules of
Proposition~\ref{prop:compatible-stopped-DPP}.
\end{proposition}

\begin{proof}
The state and time estimates in Lemma~\ref{lem:controlled-flow},
together with localisation and uniform moment bounds of order \(p>2\),
give joint continuity of the fixed-policy running-cost values, with
moduli uniform in the feedback rule on compact time-state sets.
Taking the infimum over feedbacks preserves these common local moduli
and gives the same local continuity for the value.
Consequently, the countable-neighbourhood selection and pasting
arguments used for the zero-cost DPP apply to the running-cost
problem as well.

The accumulated cost splits pathwise at the calendar time $r$:
\begin{align*}
 C_{s,t}^{s,x,\alpha}
 ={}\int_s^r f(v,X_v^{s,x,\alpha},\alpha_v^{[s]})\,dv
 +\int_r^t f(v,X_v^{s,x,\alpha},\alpha_v^{[s]})\,dv.
\end{align*}
The controlled flow identity identifies the second integral with the cost of
the shifted continuation problem.  The two inequalities in
\eqref{eq:running-cost-deterministic-DPP} now follow from exactly the law and
feedback pasting arguments used in
Theorem~\ref{thm:controlled-deterministic-DPP}.  The cost accumulated before
the pasting time is measurable with respect to the past and is therefore
unchanged by the continuation.  Truncation, the uniform moment estimates and
\eqref{eq:running-cost-growth} extend the result from bounded costs and
terminal data to $\mathscr C_V$.

Splitting the integral at the calendar time $s+\tau$ and using the stopped
controlled flow proves
\eqref{eq:running-cost-stopped-DPP}; feedback-compatible stopping rules are
handled as in Proposition~\ref{prop:compatible-stopped-DPP}.\end{proof}

\begin{remark}
It is easy to see that for fixed $s\leq t$, $\bU_{s,t}^{f}$ is monotone and translation covariant:
\begin{equation}\label{eq:running-cost-translation}
 \bU_{s,t}^{f}(\varphi+c)=\bU_{s,t}^{f}\varphi+c.
\end{equation}
For bounded terminal data it is nonexpansive in the supremum norm.  It is not,
in general, constant preserving or positively homogeneous.  Instead,
{\begin{equation*}
 \bU_{s,t}^{cf}(c\varphi)=c\bU_{s,t}^{f}\varphi,
 \qquad c\geq0,
\end{equation*}}
and adding a constant $k$ to the running cost gives
{\begin{equation*}
 \bU_{s,t}^{f+k}\varphi
 =\bU_{s,t}^{f}\varphi+k(t-s).
\end{equation*}}
\end{remark}
\begin{remark}
\label{rem:running-cost-DPP-weight}
The specific quadratic form of \(V\) is not essential to the
dynamic-programming argument in
Proposition~\ref{prop:running-cost-DPP}.
Splitting the accumulated cost, applying the controlled flow
identity, and pasting controls and probability laws do not
use this explicit form.

The choice \(V(x)=1+\abs{x}^2\) is used to justify the extension
from bounded data to the stated weighted class.
The running cost and terminal functions then have at most
quadratic growth, and the controlled moment estimates of order
\(p>2\) provide the uniform integrability needed to pass to the
limit in the truncation argument.

The same argument applies to a more general continuous weight,
provided that the corresponding finite-horizon weighted growth
bounds, uniform-integrability estimates, and continuity of the
fixed-policy costs and value functions are available.
A first-moment bound for the weight alone does not, in general,
suffice for the uniform truncation argument.
\end{remark}

For \((q,x,p,X)\in\R\times\R^d\times\R^d\times\mathbb S^d\),
define the running-cost Bellman--Isaacs operator
\begin{equation}\label{eq:running-cost-Hamiltonian}
 H_f(q,x,p,X)
 :=\inf_{a\in A}\sup_{Q\in\cQ}
 \left\{
  f(q,x,a)+b_{a,Q}(q,x)\cdot p
  +\frac12\tr({c_{a,Q}}(q,x)X)
 \right\}.
\end{equation}
When \(f=0\), this reduces to the operator \(H_0\) defined in
\eqref{eq:lower-Isaacs-operator}.

For \(L>0\), define the finite-strip weighted continuous class
\begin{equation}\label{eq:weighted-viscosity-strip-class}
 \mathscr C_{L,V}
 :=\left\{u\in C([0,L]\times\R^d):
       \sup_{\substack{0\leq\ell\leq L\\x\in\R^d}}
       \frac{\abs{u(\ell,x)}}{1+V(x)}<\infty\right\}.
\end{equation}
The growth bound is uniform on each finite strip.  Continuity here is
joint continuity in \((\ell,x)\); it does not require continuity in the
weighted norm as a function of \(\ell\).

\begin{proposition}
\label{prop:running-cost-PDE}
Suppose Assumptions~\ref{ass:uniform-parabolic},
\ref{ass:controlled-dissipativity} and \ref{ass:running-cost} hold.
Fix \(t_*\in\R\) and \(\varphi\in\mathscr C_V(\R^d)\), and put
\begin{equation}\label{eq:running-cost-remaining-value}
 u_\varphi^{t_*}(\ell,x)
 :=\bU_{t_*-\ell,t_*}^{f}\varphi(x),
 \qquad \ell\geq0.
\end{equation}
For every \(L>0\), \(u_\varphi^{t_*}\) is the unique viscosity solution
in \(\mathscr C_{L,V}\) of
\begin{equation}\label{eq:running-cost-HJB}
 \begin{cases}
 \partial_\ell u-H_f(t_*-\ell,x,Du,D^2u)=0,
       &(\ell,x)\in(0,L)\times\R^d,\\
 u(0,x)=\varphi(x),&x\in\R^d.
 \end{cases}
\end{equation}
The initial condition is attained locally uniformly in \(x\) as
\(\ell\downarrow0\).

\end{proposition}

\begin{proof}
\emph{The value equation and growth.}
The cost bound and the all-policy Lyapunov estimate give
\[
 \sup_{0\leq\ell\leq L}\abs{u_\varphi^{t_*}(\ell,x)}
 \leq C_{\varphi,f,L}(1+V(x)).
\]
The localisation and uniform-integrability argument of
Lemma~\ref{lem:controlled-weighted-extension}, using a moment of order
\(p>2\), applies also to the accumulated running cost.  It gives joint
continuity on each finite strip and local uniform convergence to
\(\varphi\) at \(\ell=0\).  For the latter, the expected absolute
running cost over an interval of length \(\ell\) tends to zero
uniformly for initial states in a compact set; the terminal term is
handled by local uniform continuity and the same moment bound.
Thus \(u_\varphi^{t_*}\in\mathscr C_{L,V}\).
The stopped DPP and the local test-function argument of
Proposition~\ref{prop:controlled-value-PDE} prove
\eqref{eq:running-cost-HJB}.  The running cost contributes the additional
term \(f(q,x,a)\); since it has no covariance argument, the continuous
approximate covariance selector and its realisation are unchanged.

\smallskip
\noindent\emph{Comparison in the weighted class.}
Fix \(L>0\).  On a bounded cylinder, the structural comparison estimate
in Proposition~\ref{prop:controlled-value-PDE} remains valid for \(H_f\),
with the additional spatial error
\(\sup_{a\in A}|f(q,x,a)-f(q,y,a)|\), which tends to zero locally
uniformly as \(|x-y|\to0\).  Hence bounded-domain parabolic comparison
\cite[Theorem~8.2]{CIL1992} applies.

To control spatial infinity for quadratic-growth solutions, use the
quartic barrier
\[
 \chi_0(x):=(1+|x|^2)^2.
\]
The linear-growth drift and bounded effective covariance imply
\[
 \sup_{a\in A,\,Q\in\cQ}\mathcal L_{a,Q}^{q}\chi_0(x)
 \leq C_\chi\chi_0(x),
 \qquad q\in\R,\quad x\in\R^d.
\]
Indeed, \(D\chi_0=4(1+|x|^2)x\) and
\(D^2\chi_0=8xx^{\mathsf T}+4(1+|x|^2)I_d\).
For \(\chi(\ell,x):=e^{(C_\chi+1)\ell}\chi_0(x)\), it follows that
\begin{equation}\label{eq:weighted-viscosity-comparison-barrier}
 \partial_\ell\chi
 -\sup_{a\in A,\,Q\in\cQ}
     \mathcal L_{a,Q}^{t_*-\ell}[\chi(\ell,\cdot)](x)
 \geq\chi.
\end{equation}
{Let \(U,Z\in\mathscr C_{L,V}\) be a viscosity subsolution and
supersolution, respectively, with \(U(0,\cdot)\leq Z(0,\cdot)\).}
For \(\varepsilon>0\), the function
\(Z_\varepsilon:=Z+\varepsilon\chi\) is a supersolution.  This follows
by subtracting \(\varepsilon\chi\) from a lower touching test and using
\begin{align*}
 &H_f(q,x,p+\varepsilon D\chi,X+\varepsilon D^2\chi)
       -H_f(q,x,p,X)\\
 &\qquad\leq\varepsilon
       \sup_{a\in A,\,Q\in\cQ}\mathcal L_{a,Q}^{q}\chi,
\end{align*}
together with \eqref{eq:weighted-viscosity-comparison-barrier}; the
common running cost cancels in this operator increment.
The uniform quadratic-growth bounds on \(U,Z\) imply
\(U\leq Z_\varepsilon\) on \([0,L]\times\partial B_R\) for every
sufficiently large \(R\), because \(\chi\) grows quartically.
Bounded-domain comparison gives \(U\leq Z_\varepsilon\) on
\([0,L)\times B_R\).  The radius can be arbitrarily large; letting
\(\varepsilon\downarrow0\) and using continuity at \(L\) proves
\(U\leq Z\) on the whole strip.  This proves uniqueness in
\(\mathscr C_{L,V}\), also when \(f=0\).
\end{proof}

\begin{corollary}
\label{cor:running-cost-difference-Pucci}
{
Under the assumptions of Proposition~\ref{prop:running-cost-PDE},
fix \(t_*\in\R\) and \(\varphi,\psi\in\mathscr C_V(\R^d)\).
Set \(u=u_\varphi^{t_*}\), \(v=u_\psi^{t_*}\), and \(w=u-v\).
Then, on every bounded open spatial set \(O\), \(w\) belongs to the
homogeneous parabolic Pucci class
}
\begin{align}
 \mathcal M^-_{\underline\lambda/2,\overline\lambda/2}(D^2w)
 -B_O|Dw|
 &\leq\partial_\ell w,
 \label{eq:running-cost-difference-Pucci-lower}\\
 \partial_\ell w
 &\leq
 \mathcal M^+_{\underline\lambda/2,\overline\lambda/2}(D^2w)
 +B_O|Dw|,
 \label{eq:running-cost-difference-Pucci-upper}
\end{align}
in the viscosity sense on \((0,\infty)\times O\).
In particular, the local Harnack-chain oscillation estimate of
Proposition~\ref{prop:controlled-local-oscillation} holds for \(w\), with
the same constants and independently of \(f\).
\end{corollary}

\begin{proof}
The same-point increment bound for \(H_f\) is the Pucci bound displayed
in the proof of Proposition~\ref{prop:controlled-difference-LOC}, since
\(f(q,x,a)\) cancels when the two generators are compared at the same
action and covariance.  The viscosity theorem of sums, with the local
uniform continuity of \(f\) removing the spatial doubling error, gives
\eqref{eq:running-cost-difference-Pucci-lower}--
\eqref{eq:running-cost-difference-Pucci-upper}.
Thus both \(w-m\) and \(M-w\) satisfy the nonnegative Harnack estimate
on every cylinder on which \(m\leq w\leq M\).  The Harnack-chain proof
of Proposition~\ref{prop:controlled-local-oscillation} applies with the
same constants.  No convexity or concavity of \(H_f\) is required.
\end{proof}

\begin{lemma}
\label{lem:running-cost-terminal-comparison}
For $\varphi,\psi\in\mathscr C_V(\R^d)$,
\begin{equation}\label{eq:running-cost-terminal-comparison}
 |\bU_{s,t}^{f}\varphi(x)-\bU_{s,t}^{f}\psi(x)|
 \leq\sup_{\alpha\in\cA}\sup_{P\in\cP_{\cQ}}
 E^P\!\left[
  |\varphi-\psi|(X_t^{s,x,\alpha})
 \right].
\end{equation}
\end{lemma}

\begin{proof}
For fixed $\alpha$, the difference of the two suprema over $P$ is bounded in
absolute value by the supremum of the absolute difference of the terminal
costs; the common accumulated running cost cancels.  The absolute difference
of two infima over $\alpha$ is bounded by the supremum over $\alpha$ of the
corresponding absolute differences.  This proves
\eqref{eq:running-cost-terminal-comparison}. In particular, if
\(\abs{\varphi-\psi}\leq c+dV\), then
\[
 \abs{\bU_{s,t}^{f}\varphi(x)-\bU_{s,t}^{f}\psi(x)}
 \leq c+d\,\mathsf M_{s,t}V(x),
\]
so the running-cost evolution satisfies
Assumption~\ref{ass:abstract-terminal-domination} with the same dominating function
\(m_{s,t}=\mathsf M_{s,t}V\) as the zero-cost evolution.
\end{proof}

\begin{theorem}
\label{thm:running-cost-contraction}
Under Assumptions~\ref{ass:uniform-parabolic},
\ref{ass:controlled-dissipativity} and \ref{ass:running-cost}, the constants
$a,\beta,q,C_0,\lambda$ in
Theorem~\ref{thm:controlled-diffusion-Harris} can be chosen so that
\begin{align}
 [\bU_{s,s+a}^{f}\varphi-\bU_{s,s+a}^{f}\psi]_{\beta,V}
 &\leq q[\varphi-\psi]_{\beta,V},
 \label{eq:running-cost-block-contraction}\\
 [\bU_{s,t}^{f}\varphi-\bU_{s,t}^{f}\psi]_{\beta,V}
 &\leq C_0e^{-\lambda(t-s)}[\varphi-\psi]_{\beta,V}
 \label{eq:running-cost-exponential-contraction}
\end{align}
for $s\leq t$ and $\varphi,\psi\in\mathscr C_V(\R^d)$.  The constants do
not depend on $f$ within the class of Assumption~\ref{ass:running-cost}.
\end{theorem}

\begin{proof}
We apply Theorem~\ref{thm:controlled-weighted-Harris} to
\(\mathsf U_{s,t}=\bU_{s,t}^{f}\) on \(\mathscr C_V(\R^d)\).
Proposition~\ref{prop:running-cost-DPP} gives the evolution property,
while \eqref{eq:running-cost-translation} gives translation covariance. Lemma~\ref{lem:running-cost-terminal-comparison}
verifies Assumption~\ref{ass:abstract-terminal-domination} with
\(m_{s,t}=\mathsf M_{s,t}V\), and
Proposition~\ref{prop:controlled-Lyapunov-verification} gives the same
Lyapunov constants as in the zero-cost problem.

It remains only to verify the local oscillation hypothesis for differences.
{For terminal data \(\varphi,\psi\),
Corollary~\ref{cor:running-cost-difference-Pucci} places the difference of
the two remaining-time values in the same homogeneous parabolic Pucci class
used in}
Proposition~\ref{prop:controlled-difference-LOC}.  Therefore the identical
finite Harnack-chain construction on
\(\mathfrak C_{a,C_R,O_R}\), followed by the same calendar-time reflection,
gives \eqref{eq:abstract-controlled-LOC} with the same
\(\vartheta_R\) and \(\mathfrak Q_{s,R}\).  The local moment bound
\eqref{eq:controlled-local-moment-bound} is unchanged because it depends only
on \(\mathsf M_{s,t}V\), not on the running cost.

Thus all hypotheses of the abstract theorem hold with the same
\(a,\beta,q,C_0,\lambda\) as in
Theorem~\ref{thm:controlled-diffusion-Harris}.  Equations
\eqref{eq:running-cost-block-contraction} and
\eqref{eq:running-cost-exponential-contraction} follow immediately.  The
argument also makes clear why the contraction constants are independent of
\(f\) within Assumption~\ref{ass:running-cost}: the common running cost
cancels in terminal perturbations and in the difference equation.
\end{proof}

\subsection{Autonomous average cost and the additive eigenvalue}

{
Assume in this subsection that $b,h,\sigma$ and $f$ are independent of calendar
time.
}  Write
\[
 \bU_t^f:=\bU_{0,t}^f.
\]
This is an additively homogeneous semigroup:
\[
 \bU_{t+r}^f=\bU_t^f\bU_r^f,
 \qquad
 \bU_t^f(\varphi+c)=\bU_t^f\varphi+c.
\]

\begin{theorem}
\label{thm:running-cost-additive-eigenpair}
Under Assumptions~\ref{ass:uniform-parabolic},
\ref{ass:controlled-dissipativity} and \ref{ass:running-cost}, there are
$\rho\in\R$ and $h\in\mathscr C_V(\R^d)$ such that
\begin{equation}\label{eq:running-cost-additive-eigenrelation}
 \bU_t^f h=h+\rho t,
 \qquad t\geq0.
\end{equation}
The number $\rho$ is unique and $h$ is unique up to an additive constant.
It is the viscosity additive eigenpair for
\begin{equation}\label{eq:running-cost-ergodic-HJB}
 \rho
 =\inf_{a\in A}\sup_{Q\in\cQ}
 \left\{
  f(x,a)+b_{a,Q}(x)\cdot Dh
  +\frac12\tr({c_{a,Q}}(x)D^2h)
 \right\}.
\end{equation}
Moreover,
for every $\varphi\in\mathscr C_V(\R^d)$,
\begin{equation}\label{eq:running-cost-average-limit}
 \frac1t\bU_t^f\varphi(x)\longrightarrow\rho
 \qquad\text{as }t\to\infty,
\end{equation}
locally uniformly in $x$, and
{the relative-value class converges exponentially:}
\begin{equation}\label{eq:running-cost-quotient-convergence}
 \inf_{c\in\R}
 \norm{\bU_t^f\varphi-\rho t-h-c}_{\beta,V}
 \leq C_0e^{-\lambda t}[\varphi-h]_{\beta,V}.
\end{equation}
\end{theorem}

\begin{proof}
The weighted continuous space \(\mathscr C_V(\R^d)\) is Banach,
and the constant functions form a closed subspace.  Hence
\(\mathscr C_V(\R^d)/\R\) is complete; the quotient identity
\eqref{eq:controlled-quotient-identity} identifies its norm with
\([\cdot]_{\beta,V}\).
Translation covariance allows \(\bU_a^f\) to act on this quotient
space, and \eqref{eq:running-cost-block-contraction} says that
the induced map is a strict contraction.
Banach's fixed-point theorem gives a unique equivalence class
\([h]\) such that
\[
 [\bU_a^f h]=[h].
\]
Thus $\bU_a^f h=h+c_a$ for a constant $c_a$; put $\rho=c_a/a$.

For $t\geq0$, commutativity of the autonomous semigroup and translation
covariance give
\[
 \bU_a^f(\bU_t^f h)
 =\bU_t^f(\bU_a^f h)
 =\bU_t^f h+c_a.
\]
Hence $[\bU_t^f h]$ is also a fixed point of the quotient map and therefore
equals $[h]$.  There is a scalar $c(t)$ such that
$\bU_t^f h=h+c(t)$.  The semigroup property implies
$c(t+r)=c(t)+c(r)$.  Continuity in $t$, which follows from the moment
estimates and Assumption~\ref{ass:running-cost}, gives $c(t)=\rho t$ and
proves \eqref{eq:running-cost-additive-eigenrelation}.
Since \(u(\ell,x)=h(x)+\rho\ell\) is the value solution of
\eqref{eq:running-cost-HJB}, spatial touching tests give the viscosity
equation \eqref{eq:running-cost-ergodic-HJB}.
Conversely, if \((\widetilde\rho,\widetilde h)\in
\R\times\mathscr C_V(\R^d)\) satisfies that stationary equation in
the viscosity sense, then
{\(\widetilde h(x)+\widetilde\rho\ell\) belongs to
\(\mathscr C_{L,V}\) for every \(L>0\) and solves}
\eqref{eq:running-cost-HJB} with initial value \(\widetilde h\).
The weighted comparison result in
Proposition~\ref{prop:running-cost-PDE} therefore implies
\(\bU_t^f\widetilde h=\widetilde h+\widetilde\rho t\).
Thus the following uniqueness argument covers viscosity eigenpairs in
this growth class as well as semigroup eigenpairs.

If $(\widetilde\rho,\widetilde h)$ is another eigenpair, then
\eqref{eq:running-cost-exponential-contraction} gives
\[
 [h-\widetilde h]_{\beta,V}
 \leq C_0e^{-\lambda t}
 [h-\widetilde h]_{\beta,V}.
\]
For sufficiently large $t$ this forces $h-\widetilde h$ to be constant.
The two eigenrelations then give $\rho=\widetilde\rho$.

Apply \eqref{eq:running-cost-exponential-contraction} to $\varphi$ and $h$,
and use \eqref{eq:running-cost-additive-eigenrelation}.  This proves the
quotient estimate \eqref{eq:running-cost-quotient-convergence}.

Finally, since $|\varphi-h|\leq C(1+V)$,
Lemma~\ref{lem:running-cost-terminal-comparison} and
\eqref{eq:all-policy-Lyapunov} show, locally uniformly in $x$, that
\[
 |\bU_t^f\varphi(x)-\bU_t^f h(x)|
 \leq C\bigl(1+\mathsf M_{0,t}V(x)\bigr)
 \leq C_x
\]
with $C_x$ locally bounded and independent of $t$.  Divide by $t$ and use
$\bU_t^f h=h+\rho t$ to obtain
\eqref{eq:running-cost-average-limit}.
\end{proof}

\begin{corollary}[Relative-value convergence]
\label{cor:running-cost-relative-value}
{
Under the assumptions of Theorem~\ref{thm:running-cost-additive-eigenpair},
for every reference point $x_\star$ and
$\varphi\in\mathscr C_V(\R^d)$,
\begin{align}
 &\left|
  \bigl(\bU_t^f\varphi(x)-\bU_t^f\varphi(x_\star)\bigr)
  -\bigl(h(x)-h(x_\star)\bigr)
 \right|\notag\\
 &\qquad\leq
 C_0e^{-\lambda t}[\varphi-h]_{\beta,V}
 \left(2+\beta V(x)+\beta V(x_\star)\right).
 \label{eq:running-cost-relative-value-convergence}
\end{align}}
\end{corollary}

\begin{proof}
{This follows by evaluating the difference at $x$ and $x_\star$ in the
contraction estimate used to prove
\eqref{eq:running-cost-quotient-convergence}.}
\end{proof}

\begin{remark}
The word ``ergodic'' is often attached to
\eqref{eq:running-cost-ergodic-HJB} in the control literature.  Here it means
an additive eigenvalue or long-time average-cost problem.  It does not assert
ergodicity of an invariant sublinear expectation in the capacity sense.
The outer infimum over feedback controls means that neither
\(\bU_t^f\) nor its normalized long-time functional is generally sublinear,
so the \((0,0)\) invariant-capacity ergodicity theory developed by Feng and
Zhao~\cite{FengZhao2021} does not apply directly here.
\end{remark}

\subsection{The genuinely nonstationary additive problem}

For arbitrary calendar-time dependence, the autonomous eigenrelation
\(\bU_t^fh=h+\rho t\) is generally too rigid.  The correct replacement is a
complete orbit in the quotient by constants, together with an additive
scalar cocycle.

The equations derived in this section are \emph{infinite-horizon corrector
equations}, in contrast with the finite-horizon value equations
\eqref{eq:controlled-remaining-HJB},
\eqref{eq:controlled-conjugate-HJB} and
\eqref{eq:running-cost-HJB}.
In those earlier equations, a terminal calendar time \(t_*\) and a
terminal payoff are prescribed; in remaining-time coordinates, the
payoff is imposed as the initial condition at \(\ell=0\).
Although these equations may be considered for every \(\ell>0\), each
value still has a finite remaining horizon \(\ell\) and terminal
date~\(t_*\).

Here the normalized corrector is instead obtained by sending the
terminal horizon to infinity, as in
\eqref{eq:nonstationary-corrector-limit}.
The resulting equation is not supplemented by an independently
prescribed payoff at a finite terminal date.
Its solution is selected by the additive evolution relation, the
normalization at \(x_\star\), and the uniform weighted bound on the
normalized spatial profile.
In the general time-inhomogeneous case, subtracting a scalar
primitive of the additive cocycle gives an all-time viscosity
solution of \eqref{eq:nonstationary-corrector-HJB}.
In the periodic case, temporal periodicity replaces the prescribed
terminal datum, and the constant average cost \(\rho\) is determined
together with the periodic corrector through
\eqref{eq:periodic-ergodic-HJB}.

Thus the finite-horizon and infinite-horizon problems involve the
same local Bellman--Isaacs operator, but different global conditions
select their solutions.
The autonomous equation \eqref{eq:running-cost-ergodic-HJB} is the
time-independent counterpart of these infinite-horizon corrector
equations.  Here infinite horizon refers to the corrector and
long-run average-cost problem, not to convergence of the undiscounted
total running cost.

\begin{theorem}
\label{thm:nonstationary-additive-cocycle}
Fix \(x_\star\in\R^d\).
Under Assumptions~\ref{ass:uniform-parabolic},
\ref{ass:controlled-dissipativity} and \ref{ass:running-cost}, there exist
a family
\[
 h_t\in\mathscr C_V(\R^d),\qquad t\in\R,
\]
and real numbers \(\Lambda(s,t)\), \(s\leq t\), such that
\begin{align}
 h_t(x_\star)&=0,
 \qquad
 \sup_{t\in\R}\norm{h_t}_{\beta,V}<\infty,
 \label{eq:nonstationary-corrector-normalisation}\\
 \bU_{s,t}^fh_t
 &=h_s+\Lambda(s,t),
 \label{eq:nonstationary-additive-relation}\\
 \Lambda(s,t)
 &=\Lambda(s,r)+\Lambda(r,t),
 \qquad s\leq r\leq t.
 \label{eq:nonstationary-cocycle}
\end{align}
The pair \(\bigl((h_t)_{t\in\R},\Lambda\bigr)\) is unique among all
pairs satisfying \eqref{eq:nonstationary-corrector-normalisation}
and \eqref{eq:nonstationary-additive-relation}.

Moreover, for every \(\varphi\in\mathscr C_V(\R^d)\),
\begin{align}
 &[\bU_{s,t}^f\varphi-h_s-\Lambda(s,t)]_{\beta,V}
 \leq C_0e^{-\lambda(t-s)}
 [\varphi-h_t]_{\beta,V},
 \label{eq:nonstationary-relative-contraction}\\
 &\left|
  \bigl(\bU_{s,t}^f\varphi(x)
        -\bU_{s,t}^f\varphi(x_\star)\bigr)
  -h_s(x)
 \right|\notag\\
 &\qquad\leq
 C_0e^{-\lambda(t-s)}
 [\varphi-h_t]_{\beta,V}
 \bigl(2+\beta V(x)+\beta V(x_\star)\bigr).
 \label{eq:nonstationary-relative-value}
\end{align}
\end{theorem}

\begin{proof}
First observe that the one-block forcing is uniformly bounded in the
weighted quotient norm:
\begin{equation}\label{eq:nonstationary-one-block-forcing}
 K_f:=\sup_{s\in\R}[\bU_{s,s+a}^f0]_{\beta,V}<\infty.
\end{equation}
Indeed, Assumption~\ref{ass:running-cost} and the short-time all-policy
Lyapunov estimate give
\[
 |\bU_{s,s+a}^f0(x)|
 \leq c_f\int_s^{s+a}
 \bigl(1+\mathsf M_{s,r}V(x)\bigr)\,dr
 \leq C_{f,a}(1+V(x)),
\]
uniformly in \(s\).  This implies
\eqref{eq:nonstationary-one-block-forcing}.

{
For \(T\geq s\), set
\[
 H_{s,T}:=\bU_{s,T}^f0,\qquad
 \widehat H_{s,T}(x):=H_{s,T}(x)-H_{s,T}(x_\star).
\]
In particular, \(H_{s,s}=\widehat H_{s,s}=0\) by the identity
property of the evolution.
}
If \(T=s+na\), the evolution property and
\eqref{eq:running-cost-exponential-contraction} yield
\begin{align*}
 &[H_{s,s+(n+1)a}-H_{s,s+na}]_{\beta,V}\\
 &\qquad\leq
 C_0e^{-\lambda na}
 [\bU_{s+na,s+(n+1)a}^f0]_{\beta,V}
 \leq C_0K_fe^{-\lambda na}.
\end{align*}
Since
\[
 \sum_{n=0}^{\infty} C_0K_f e^{-\lambda na}<\infty,
\]
the successive differences are summable in the weighted seminorm.  Moreover,
\[
 \widehat H_{s,s+(n+1)a}(x_\star)
 =\widehat H_{s,s+na}(x_\star)=0.
\]
Hence their difference also vanishes at \(x_\star\).  For any
\(g\in\mathscr C_V(\R^d)\) with \(g(x_\star)=0\),
\[
 |g(x)|
 \leq [g]_{\beta,V}
 \bigl(2+\beta V(x)+\beta V(x_\star)\bigr),
\]
and therefore
\[
 \norm{g}_{\beta,V}
 \leq \bigl(2+\beta V(x_\star)\bigr)[g]_{\beta,V}.
\]
Applying this with
\[
 g=\widehat H_{s,s+(n+1)a}-\widehat H_{s,s+na},
\]
and observing that subtraction of the normalizing constants does not change
the weighted seminorm, we obtain
\[
 \norm{\widehat H_{s,s+(n+1)a}-\widehat H_{s,s+na}}_{\beta,V}
 \leq
 \bigl(2+\beta V(x_\star)\bigr)C_0K_f e^{-\lambda na}.
\]
The right-hand side is summable in \(n\).  Thus
\(\bigl(\widehat H_{s,s+na}\bigr)_{n\geq1}\) is a Cauchy sequence in
\(\mathscr C_V(\R^d)\), with the estimate uniform in \(s\), and hence
converges in the weighted norm.  Denote its limit by \(h_s\).
If \(s+na\leq T<s+(n+1)a\), the short-time estimate gives a uniform bound for
\([\bU_{s+na,T}^f0]_{\beta,V}\), and contraction over the preceding \(n\)
blocks shows that \(\widehat H_{s,T}\) has the same limit.  Thus
\begin{equation}\label{eq:nonstationary-corrector-limit}
 h_s
 =\lim_{T\to\infty}
 \left(
   \bU_{s,T}^f0-(\bU_{s,T}^f0)(x_\star)
 \right)
 \qquad\text{in }\mathscr C_V(\R^d).
\end{equation}
The geometric-series bound proves
\eqref{eq:nonstationary-corrector-normalisation}.
 
{
Fix \(s<t\).  Apply \(\bU_{s,t}^f\) to the normalized approximations
defining \(h_t\).  The contraction estimate
\eqref{eq:running-cost-exponential-contraction} justifies passing to
the limit in the quotient modulo constants.  Translation covariance,
the evolution property and \eqref{eq:nonstationary-corrector-limit}
then show that \(\bU_{s,t}^fh_t\) and \(h_s\) represent the same
quotient element.  Equivalently,
}
\(
\bU_{s,t}^fh_t-h_s
\)
is constant in \(x\). Denote this constant by \(\Lambda(s,t)\). Then
\[
\bU_{s,t}^fh_t=h_s+\Lambda(s,t),
\]
which proves \eqref{eq:nonstationary-additive-relation}.
Set \(\Lambda(s,s):=0\), so the same relation holds also on the
diagonal.
Applying the
evolution property to \(h_t\) and using translation covariance gives
\[
 h_s+\Lambda(s,t)
 =\bU_{s,r}^f(h_r+\Lambda(r,t))
 =h_s+\Lambda(s,r)+\Lambda(r,t),
\]
which proves \eqref{eq:nonstationary-cocycle}.

To prove uniqueness, let
\((\widetilde h_t,\widetilde\Lambda)\) be another family satisfying
the same normalization, uniform bound and additive relation.
By \eqref{eq:running-cost-exponential-contraction},
\[
 [h_s-\widetilde h_s]_{\beta,V}
 \leq
 C_0e^{-\lambda(t-s)}
 [h_t-\widetilde h_t]_{\beta,V},
 \qquad s\leq t.
\]
By the uniform bounds on \(h_t\) and \(\widetilde h_t\), the seminorm
on the right-hand side is uniformly bounded in \(t\). Hence, for fixed
\(s\), letting \(t\to\infty\) gives
\[
 [h_s-\widetilde h_s]_{\beta,V}=0.
\]
Thus \(h_s-\widetilde h_s\) is constant in \(x\), and the normalization
at \(x_\star\) shows that \(h_s=\widetilde h_s\).
Finally, evaluating \eqref{eq:nonstationary-additive-relation} at
\(x_\star\) gives
\[
 \Lambda(s,t)=\bU_{s,t}^fh_t(x_\star),
\]
and therefore the cocycle is unique as well.

Now compare terminal data \(\varphi\) and \(h_t\) in
\eqref{eq:running-cost-exponential-contraction}, and use
\eqref{eq:nonstationary-additive-relation}.  This gives
\eqref{eq:nonstationary-relative-contraction}; evaluating the difference at
\(x\) and \(x_\star\) gives
\eqref{eq:nonstationary-relative-value}.
\end{proof}

\begin{corollary}
\label{cor:nonstationary-cocycle-tracking}
{
Under the assumptions and notation of
Theorem~\ref{thm:nonstationary-additive-cocycle}, for every fixed
\(\varphi\in\mathscr C_V(\R^d)\) there is \(C_\varphi<\infty\) such that
\begin{equation}\label{eq:nonstationary-cocycle-tracking}
 \left|
  \bU_{s,t}^f\varphi(x)-\Lambda(s,t)
 \right|
 \leq C_\varphi(1+V(x)),
 \qquad s\leq t.
\end{equation}
Consequently,
\begin{equation}\label{eq:nonstationary-average-tracking}
 \frac{\bU_{s,t}^f\varphi(x)}{t-s}
 -\frac{\Lambda(s,t)}{t-s}
 \longrightarrow0
 \qquad\text{as }t-s\to\infty,
\end{equation}
locally uniformly in \(x\).}
\end{corollary}

\begin{proof}
{
The uniform bound on \(h_t\), the terminal comparison lemma and the
all-policy Lyapunov estimate give
\[
 |\bU_{s,t}^f\varphi(x)-\bU_{s,t}^fh_t(x)|
 \leq C_\varphi\bigl(1+\mathsf M_{s,t}V(x)\bigr)
 \leq C_\varphi(1+V(x)).
\]
Combining this with
\(\bU_{s,t}^fh_t=h_s+\Lambda(s,t)\) and the uniform weighted bound on
\(h_s\) proves \eqref{eq:nonstationary-cocycle-tracking}.  Dividing by
\(t-s\) proves \eqref{eq:nonstationary-average-tracking}.}
\end{proof}

\begin{corollary}
\label{cor:nonstationary-corrector-HJB}
Under the assumptions and notation of
Theorem~\ref{thm:nonstationary-additive-cocycle}, choose
\(A:\R\to\R\), normalized by \(A(0)=0\), such that
\[
 \Lambda(s,t)=A(t)-A(s),\qquad s\leq t,
\]
and define
\begin{equation}\label{eq:nonstationary-gauged-corrector}
 W(s,x):=h_s(x)-A(s).
\end{equation}
Then \(W\) is jointly continuous and, for every compact interval
\(I\subset\R\),
\[
 \sup_{s\in I}\sup_{x\in\R^d}
 \frac{|W(s,x)|}{1+V(x)}<\infty,
\]
and 
by translation covariance and
\eqref{eq:nonstationary-additive-relation},
\begin{equation}\label{eq:nonstationary-gauged-evolution}
 \bU_{s,t}^f W(t,\cdot)=W(s,\cdot),
 \qquad s\leq t.
\end{equation}
Moreover, \(W\) is a viscosity solution on \(\R\times\R^d\) of the Hamilton-Jacobi-Bellman equation
\begin{equation}\label{eq:nonstationary-corrector-HJB}
 -\partial_s W-H_f(s,x,D_xW,D_x^2W)=0.
\end{equation}
Here \(s\) is calendar time, and \(H_f\) is as in
\eqref{eq:running-cost-Hamiltonian}.  In particular,
\[
 A(s)=-W(s,x_\star),
\]
so \(A\) is continuous and \((s,x)\mapsto h_s(x)\) is jointly
continuous.  No differentiability of \(A\) is required.
\end{corollary}

\begin{proof}
The normalized primitive \(A\) may be written explicitly as
\[
 A(t):=
 \begin{cases}
  \Lambda(0,t),&t\geq0,\\
  -\Lambda(t,0),&t<0.
 \end{cases}
\]
The cocycle identity gives \(A(0)=0\) and
\(\Lambda(s,t)=A(t)-A(s)\) for all \(s\leq t\).
Fix \(b\in\R\).  The additive relation yields
{\begin{equation*}
 W(s,x)=(\bU_{s,b}^fh_b)(x)-A(b),\qquad s\leq b.
\end{equation*}}
{For every \(L>0\), Proposition~\ref{prop:running-cost-PDE}, applied to
\(h_b\in\mathscr C_V(\R^d)\), shows that
\[
 u_b(\ell,x):=(\bU_{b-\ell,b}^fh_b)(x),\qquad 0\leq\ell\leq L,
\]
belongs to \(\mathscr C_{L,V}\) and solves
\eqref{eq:running-cost-HJB}.  Since
\(W(s,x)=u_b(b-s,x)-A(b)\), the time change \(\ell=b-s\)
transforms viscosity test functions and gives
\eqref{eq:nonstationary-corrector-HJB} on
\((b-L,b)\times\R^d\).}  The subtraction of \(A(b)\) does not alter
the equation because the Hamiltonian is independent of the
undifferentiated solution.  {The same representation gives joint
continuity and the weighted growth bound on \([b-L,b]\times\R^d\).
Since \(b\) and \(L\) are arbitrary, these conclusions hold on every}
finite calendar-time strip.  Finally, \(h_s(x_\star)=0\) gives
\(A(s)=-W(s,x_\star)\), hence continuity of \(A\) and joint
continuity of \(h_s(x)=W(s,x)+A(s)\).
\end{proof}

\begin{remark}
The Hamilton--Jacobi--Bellman equation in
Corollary~\ref{cor:nonstationary-corrector-HJB} is an all-time
corrector equation, rather than a finite-horizon terminal-value problem.  The function \(W\) need not be uniformly weighted-bounded in
calendar time, since its scalar component \(A(s)\) may grow.  The uniform
bound applies instead to the normalized spatial profile
\[
 W(s,x)-W(s,x_\star)=h_s(x),
 \qquad
 \sup_{s\in\R}\norm{h_s}_{\beta,V}<\infty.
\]
\end{remark}

The periodic case is closely related to the ergodic stochastic control
problem studied by Ma--Zhao~\cite{MaZhao2026}.  For controlled diffusions
with time-periodic coefficients, they obtain a constant long-run average
cost together with a periodic corrector solving the associated ergodic
dynamic-programming equation.  Corollary~\ref{cor:periodic-running-cost-eigenpair}
has the analogous structure in the present covariance-uncertain setting:
the nonstationary additive cocycle reduces, after a periodic scalar gauge
change, to a constant average cost \(\rho\) and a periodic corrector family.
The weighted contraction obtained here additionally yields exponential
convergence of relative values.

\begin{corollary}
\label{cor:periodic-running-cost-eigenpair}
Suppose Assumptions~\ref{ass:uniform-parabolic},
\ref{ass:controlled-dissipativity} and \ref{ass:running-cost} hold, and
that the controlled coefficients and \(f\) are \(T\)-periodic in
calendar time.  Then there exist a unique
\(\rho\in\R\) and a \(T\)-periodic corrector family
\(\widetilde h_t\in\mathscr C_V(\R^d)\), unique up to the addition of one
common constant, such that
\begin{equation}\label{eq:periodic-additive-eigenrelation}
 \bU_{s,t}^f\widetilde h_t
 =\widetilde h_s+\rho(t-s),
 \qquad s\leq t.
\end{equation}
Writing \(\widetilde h(s,x):=\widetilde h_s(x)\), the corrector is
jointly continuous and satisfies
\[
 \sup_{s\in[0,T]}\norm{\widetilde h_s}_{\beta,V}<\infty.
\]

The pair \((\rho,\widetilde h)\) solves the periodic
infinite-horizon corrector problem, expressed by the ergodic
Bellman--Isaacs equation
\begin{equation}\label{eq:periodic-ergodic-HJB}
 \begin{cases}
 -\partial_s\widetilde h
 -H_f(s,x,D_x\widetilde h,D_x^2\widetilde h)+\rho=0,
       &(s,x)\in\R\times\R^d,\\
 \widetilde h(s+T,x)=\widetilde h(s,x),
       &(s,x)\in\R\times\R^d,
 \end{cases}
\end{equation}
in the viscosity sense.  Conversely, every continuous \(T\)-periodic
viscosity pair \((r,k)\) satisfying
\[
 \sup_{s\in[0,T]}\norm{k(s,\cdot)}_{\beta,V}<\infty
\]
satisfies the corresponding additive evolution relation, with
\(r=\rho\).  Thus the periodic viscosity formulation and the additive
evolution formulation are equivalent in this class.
\end{corollary}

\begin{proof}
{Periodicity of the evolution and uniqueness of the normalized
corrector--cocycle pair imply
\[
 h_{t+T}=h_t,
 \qquad
 \Lambda(s+T,t+T)=\Lambda(s,t).
\]}
Cocycle property of \(\Lambda\), together with its periodicity, then shows that \(\Lambda(s,s+T)\) is independent of \(s\).
Set
\[
 \rho:=\frac{\Lambda(0,T)}{T}.
\]
Use the continuous function \(A\), normalized by \(A(0)=0\),
from Corollary~\ref{cor:nonstationary-corrector-HJB}, and put
\[
 p(t):=A(t)-\rho t.
\]
{
Then \(p\) is continuous and \(T\)-periodic, so it is bounded on
\([0,T]\) and hence on \(\R\).
}
Define \(\widetilde h_t:=h_t-p(t)\).
Joint continuity of \(h\) and continuity of \(A\) show that
\(\widetilde h\) is jointly continuous.  The uniform weighted bound
on \(h_t\) and boundedness of \(p\) give the stated weighted bound
on \(\widetilde h\).
Translation covariance gives
\[
 \bU_{s,t}^f\widetilde h_t
 =h_s+\Lambda(s,t)-p(t)
 =\widetilde h_s+\rho(t-s),
\]
which proves \eqref{eq:periodic-additive-eigenrelation}.  Uniqueness follows
from the contraction on the quotient and evaluation over one period.

For the viscosity equation, fix \(b\in\R\) and use
\eqref{eq:periodic-additive-eigenrelation} to write
\[
 (\bU_{b-\ell,b}^f\widetilde h_b)(x)
 =\widetilde h(b-\ell,x)+\rho\ell.
\]
By Proposition~\ref{prop:running-cost-PDE}, the left-hand side solves
\eqref{eq:running-cost-HJB}.  Changing variables from remaining time
\(\ell\) to calendar time \(s=b-\ell\), and subtracting the smooth
term \(\rho\ell\) in the test-function argument, gives
\eqref{eq:periodic-ergodic-HJB}.  Equivalently,
\(W(s,x)=\widetilde h(s,x)-\rho s\) and
\eqref{eq:nonstationary-corrector-HJB} give the same equation.
Neither argument differentiates the periodic scalar correction \(p\).

Conversely, let \((r,k)\) be a continuous \(T\)-periodic viscosity
pair satisfying \eqref{eq:periodic-ergodic-HJB} with \(r\) in place
of \(\rho\), and assume
\[
 \sup_{s\in[0,T]}\norm{k(s,\cdot)}_{\beta,V}<\infty.
\]
By the finite-horizon comparison argument of
Proposition~\ref{prop:running-cost-PDE},
\[
 \bU_{s,t}^fk(t,\cdot)
 =k(s,\cdot)+r(t-s),
 \qquad s\leq t.
\]
Comparing this with
\eqref{eq:periodic-additive-eigenrelation} and applying
\eqref{eq:running-cost-exponential-contraction} over \(nT\) shows that
\(k(s,\cdot)-\widetilde h_s=c(s)\) is spatially constant.
The two additive relations then give
\[
 c(t)-c(s)=(r-\rho)(t-s).
\]
Taking \(t=s+T\) and using periodicity yields \(r=\rho\), and hence
\(c\) is constant.  Thus the corrector is unique up to one common
additive constant.
\end{proof}

\begin{corollary}
\label{cor:periodic-running-cost-convergence}
Under the assumptions and notation of
Corollary~\ref{cor:periodic-running-cost-eigenpair}, for every
\(\varphi\in\mathscr C_V(\R^d)\),
\begin{equation}\label{eq:periodic-average-cost-limit}
 \frac1{t-s}\bU_{s,t}^f\varphi(x)\longrightarrow\rho
 \qquad\text{as }t-s\to\infty,
\end{equation}
locally uniformly in \(x\), uniformly in the initial phase
\(s\) modulo \(T\).
Moreover, relative values converge exponentially to the periodic
corrector:
\begin{align}
 &\left|
  \bigl(\bU_{s,t}^f\varphi(x)
        -\bU_{s,t}^f\varphi(x_\star)\bigr)
  -\bigl(\widetilde h_s(x)-\widetilde h_s(x_\star)\bigr)
 \right|\notag\\
 &\qquad\leq
 C_0e^{-\lambda(t-s)}
 [\varphi-\widetilde h_t]_{\beta,V}
 \bigl(2+\beta V(x)+\beta V(x_\star)\bigr)
 \label{eq:periodic-relative-value-convergence}
\end{align}
for all \(s\leq t\) and \(x,x_\star\in\R^d\).
\end{corollary}
\begin{proof}
{Let \(A\) be the scalar function from
Corollary~\ref{cor:nonstationary-corrector-HJB} and set
\[
 p(t):=A(t)-\rho t.
\]
By the proof of Corollary~\ref{cor:periodic-running-cost-eigenpair},
\(p\) is bounded and \(T\)-periodic, and
\[
 \Lambda(s,t)=\rho(t-s)+p(t)-p(s).
\]
Hence the uniform tracking estimate
\eqref{eq:nonstationary-cocycle-tracking} yields}
\[
 \left|\frac{\bU_{s,t}^f\varphi(x)}{t-s}-\rho\right|
 \leq\frac{C_\varphi(1+V(x))+2\norm{p}_\infty}{t-s},
 \qquad s<t.
\]
This proves \eqref{eq:periodic-average-cost-limit}, including its
uniformity in the initial phase.
Finally,
\eqref{eq:periodic-relative-value-convergence} is
\eqref{eq:nonstationary-relative-value} after the scalar gauge change from
\(h_t\) to \(\widetilde h_t\).
\end{proof}

\begin{remark}
\label{rem:average-cost-policy-scope}
The average-cost convergence results concern the asymptotics of the
optimized finite-horizon values.  The construction of a single
admissible feedback attaining the limiting average cost over an
infinite horizon is not asserted here.
\end{remark}


\subsection{The fixed Brownian-law case with progressive controls}
\label{sec:Brownian-progressive-controls}

Let \(A\) be a nonempty compact metric space, and let \(P_0\) be
Wiener measure on the canonical space \(C_0([0,\infty);\R^m)\),
with coordinate process \(B\) and its usual augmented filtration.
Write \(\mathcal U_{\mathrm{prog}}\) for all \(A\)-valued
progressively measurable controls for this filtration, indexed by
elapsed time.  For \(s\in\R\), \(x\in\R^d\) and
\(\alpha\in\mathcal U_{\mathrm{prog}}\), consider
\[
 dX_{s+u}^{s,x,\alpha}
 =b_0(s+u,X_{s+u}^{s,x,\alpha},\alpha_u)\,du
  +\sigma(s+u,X_{s+u}^{s,x,\alpha},\alpha_u)\,dB_u,
 \qquad X_s^{s,x,\alpha}=x.
\]
Here
\(b_0:\R\times\R^d\times A\to\R^d\),
\(\sigma:\R\times\R^d\times A\to\R^{d\times m}\), and
\(f:\R\times\R^d\times A\to\R\).
In the model of this paper with \(\cQ=\{I_m\}\), the ordinary drift
is \(b_0=b+\sum_{i=1}^m h_{ii}\).

Set \(V(x):=1+|x|^2\).  For \(s\leq t\) and
\(\varphi\in\mathscr C_V(\R^d)\), define
\begin{align*}
 \bU_{s,t}^{f,\mathrm{prog}}\varphi(x)
 :=\inf_{\alpha\in\mathcal U_{\mathrm{prog}}}
 E^{P_0}\!\left[
  \int_0^{t-s}f(s+u,X_{s+u}^{s,x,\alpha},\alpha_u)\,du
  +\varphi(X_t^{s,x,\alpha})
 \right],
 \qquad
 \bS_{s,t}^{\mathrm{prog}}:=\bU_{s,t}^{0,\mathrm{prog}}.
\end{align*}
Write \(\bU_{s,t}^{f,\mathrm{el}}\) for the same value with the
infimum restricted to the elementary canonical feedback class
\(\cA\).  Define the ordinary Bellman Hamiltonian
\[
 H_f^0(q,x,p,X)
 :=\inf_{a\in A}\left\{
 f(q,x,a)+b_0(q,x,a)\cdot p
 +\frac12\tr\!\left(\sigma(q,x,a)\sigma(q,x,a)^{\mathsf T}X\right)
 \right\}.
\]

\begin{theorem}
\label{thm:Brownian-progressive-summary}
Under the preceding setup, suppose that
Assumptions~\ref{ass:uniform-parabolic},
\ref{ass:controlled-dissipativity} and \ref{ass:running-cost}
hold with
\(
 \cQ=\{I_m\},\ b=b_0,\ h_{ij}=0.
\)
Then the following conclusions hold.
\begin{enumerate}[label=\textup{(\roman*)}]
\item 
For every \(s\leq t\) and \(\varphi\in\mathscr C_V(\R^d)\),
\begin{equation}\label{eq:Brownian-progressive-value-equivalence}
 \bU_{s,t}^{f,\mathrm{prog}}\varphi
 =\bU_{s,t}^{f,\mathrm{el}}\varphi.
\end{equation}
The deterministic and stopped dynamic programming principles hold
for the progressive-control problem (with the same compatibility
requirement for control-dependent state-path stopping rules).
For fixed \(t_*\in\R\), the function
\(u(\ell,x):=\bU_{t_*-\ell,t_*}^{f,\mathrm{prog}}\varphi(x)\)
is, on each finite strip, the unique viscosity solution in the
class \(\mathscr C_{L_*,V}\), \(L_*>0\), of
\[
 \partial_\ell u-H_f^0(t_*-\ell,x,D_xu,D_x^2u)=0,
 \qquad u(0,x)=\varphi(x).
\]

\item 
There are \(\beta>0\), \(\lambda>0\) and \(C\geq1\), independent
of \(f\) within the stated class, such that
\begin{equation}\label{eq:Brownian-progressive-contraction}
 [\bU_{s,t}^{f,\mathrm{prog}}\varphi
  -\bU_{s,t}^{f,\mathrm{prog}}\psi]_{\beta,V}
 \leq Ce^{-\lambda(t-s)}[\varphi-\psi]_{\beta,V},
 \qquad s\leq t.
\end{equation}
For zero running cost, the pullback limit
\[
 \mathcal E_t(\varphi)
 :=\lim_{r\to-\infty}\bS_{r,t}^{\mathrm{prog}}\varphi(x)
\]
exists independently of \(x\) and defines the unique entrance
family in the uniformly weighted-bounded class of
Theorem~\ref{thm:controlled-large-time}.  It satisfies
\[
 \mathcal E_s(\bS_{s,t}^{\mathrm{prog}}\varphi)
 =\mathcal E_t(\varphi),
 \qquad
 |\bS_{s,t}^{\mathrm{prog}}\varphi(x)-\mathcal E_t(\varphi)|
 \leq Ce^{-\lambda(t-s)}[\varphi]_{\beta,V}(1+V(x)),
\]
after increasing \(C\) if necessary.

\item
For every fixed \(x_\star\in\R^d\), there is a unique
corrector--cocycle pair \(((h_t)_{t\in\R},\Lambda)\) satisfying
\[
 h_t\in\mathscr C_V(\R^d),\qquad
 h_t(x_\star)=0,\qquad
 \sup_{t\in\R}\|h_t\|_{\beta,V}<\infty,
\]
and
\[
 \bU_{s,t}^{f,\mathrm{prog}}h_t=h_s+\Lambda(s,t),
 \qquad
 \Lambda(s,t)=\Lambda(s,r)+\Lambda(r,t),
 \quad s\leq r\leq t.
\]
Moreover,
\[
 [\bU_{s,t}^{f,\mathrm{prog}}\varphi-h_s-\Lambda(s,t)]_{\beta,V}
 \leq Ce^{-\lambda(t-s)}[\varphi-h_t]_{\beta,V}.
\]
There is a continuous function \(A_0:\R\to\R\), with
\(A_0(0)=0\) and \(\Lambda(s,t)=A_0(t)-A_0(s)\), such that
\(W(s,x):=h_s(x)-A_0(s)\) is jointly continuous, weighted-bounded
on compact time intervals, and satisfies
\[
 -\partial_sW-H_f^0(s,x,D_xW,D_x^2W)=0
 \quad\text{on }\R\times\R^d
\]
in the viscosity sense.

\item 
If \(b_0,\sigma,f\) are autonomous, then \(h_t=h\) and
\(\Lambda(s,t)=\rho(t-s)\), where \(\rho\) is unique and
\[
 \rho=H_f^0(x,Dh,D^2h)
\]
in the viscosity sense; the corrector is unique modulo constants.
If instead the coefficients and \(f\) are \(T\)-periodic, there
are a unique \(\rho\) and a jointly continuous \(T\)-periodic
corrector \(\widetilde h\), uniformly weighted-bounded over one
period and unique up to one common additive constant, such that
\[
 \bU_{s,t}^{f,\mathrm{prog}}\widetilde h_t
 =\widetilde h_s+\rho(t-s),
 \qquad
 -\partial_s\widetilde h
 -H_f^0(s,x,D_x\widetilde h,D_x^2\widetilde h)+\rho=0
\]
with the equation understood in the viscosity sense.
In either case,
\[
 \frac{\bU_{s,t}^{f,\mathrm{prog}}\varphi(x)}{t-s}
 \longrightarrow\rho
 \qquad\text{as }t-s\to\infty,
\]
locally uniformly in \(x\), and uniformly in the initial phase
modulo \(T\) in the periodic case.  The corresponding relative
values converge exponentially as in part \textup{(iii)}.
\end{enumerate}
\end{theorem}

\begin{proof}
We first justify enlargement of the admissible class.  Fix
\(s<t\), \(x\), \(\varphi\), and
\(\alpha\in\mathcal U_{\mathrm{prog}}\), and put \(H=t-s\).
There are adapted piecewise-constant controls \(\alpha^n\), on
finite deterministic partitions of \([0,H]\), such that
\[
 E^{P_0}\int_0^H
 \bigl(1\wedge d_A(\alpha_u^n,\alpha_u)\bigr)\,du
 \longrightarrow0.
\]
This is the usual approximation by elementary predictable processes;
see also \cite[Section~4]{ElKarouiTan2015} for piecewise-constant
approximation in control problems.  The passage from progressive to
predictable versions is up to \(du\times P_0\)-null sets, which
leaves both integrals in the state equation unchanged.  At each
of the finitely many grid times, choose a version of the decision
variable measurable with respect to the unaugmented canonical
sigma-field.  Such versions exist because the usual Brownian augmentation does not
change deterministic-time sigma-fields modulo \(P_0\)-null sets.  Extending the control by
one fixed action beyond \(H\) gives a member of \(\cA\).

Let \(X^n\) and \(X\) be the corresponding states driven by the
same Brownian motion.  The uniform state-Lipschitz bounds, the
Burkholder--Davis--Gundy inequality and Gronwall's lemma give
\begin{align*}
 E^{P_0}\sup_{0\leq u\leq H}|X^n_{s+u}-X_{s+u}|^2
 &\leq C_H E^{P_0}\int_0^H
 \Bigl(
 |b_0(s+u,X_{s+u},\alpha_u^n)
    -b_0(s+u,X_{s+u},\alpha_u)|^2\\
 &
 +\|\sigma(s+u,X_{s+u},\alpha_u^n)
    -\sigma(s+u,X_{s+u},\alpha_u)\|_{\mathrm F}^2
 \Bigr)\,du
 \longrightarrow0.
\end{align*}
The convergence follows from continuity in the action, localisation
on compact time-state sets, and the linear-growth bounds.
The state moments of every order \(p\geq2\) are bounded uniformly
in \(n\) on this fixed horizon.  Choosing \(p>2\) gives uniform
integrability of the quadratic-growth terminal and running costs.
Consequently, if \(J_{s,t}^{f,\alpha}\varphi(x)\) denotes the
expected cost under the fixed control \(\alpha\), then
\[
 J_{s,t}^{f,\alpha^n}\varphi(x)
 \longrightarrow J_{s,t}^{f,\alpha}\varphi(x).
\]
Since \(\cA\subset\mathcal U_{\mathrm{prog}}\),
\(\bU_{s,t}^{f,\mathrm{prog}}\varphi
 \leq\bU_{s,t}^{f,\mathrm{el}}\varphi\).
Conversely,
\[
 \bU_{s,t}^{f,\mathrm{el}}\varphi(x)
 \leq\lim_{n\to\infty}J_{s,t}^{f,\alpha^n}\varphi(x)
 =J_{s,t}^{f,\alpha}\varphi(x).
\]
Taking the infimum over \(\alpha\) proves
\eqref{eq:Brownian-progressive-value-equivalence}.

The stated coefficient conditions give the special case
\(\cQ=\{I_m\}\), \(b=b_0\), \(h_{ij}=0\), of
Assumptions~\ref{ass:uniform-parabolic},
\ref{ass:controlled-dissipativity} and \ref{ass:running-cost}.
The viscosity characterization follows from
Proposition~\ref{prop:running-cost-PDE} and value equivalence.
The deterministic and stopped DPPs for the progressive formulation
follow by the same conditioning and countable continuation-pasting
arguments as in Proposition~\ref{prop:running-cost-DPP}:
progressive controls are stable under these operations, the
fixed-law state equation is well posed, and the fixed-control
continuation values have the common local continuity moduli.
State-path stopping rules are stopping times for the continuous
adapted state under the usual augmented Brownian filtration;
the same nonanticipativity property gives pasting compatibility.

Part \textup{(ii)} now follows from
Theorems~\ref{thm:running-cost-contraction} and
\ref{thm:controlled-large-time}; part \textup{(iii)} follows from
Theorem~\ref{thm:nonstationary-additive-cocycle} and
Corollary~\ref{cor:nonstationary-corrector-HJB}; and part
\textup{(iv)} follows from
Theorem~\ref{thm:running-cost-additive-eigenpair} and
Corollaries~\ref{cor:periodic-running-cost-eigenpair} and
\ref{cor:periodic-running-cost-convergence}.
The equality of values holds separately on every finite horizon,
so this transfer needs no approximation estimate uniform over
infinite horizons.
\end{proof}

\noindent
The theorem concerns the optimized value evolutions.  It does not
assert the existence of a single admissible feedback attaining the
infinite-horizon average cost.  For general nonstationary coefficients,
a limiting scalar average is not asserted; the corrector and cocycle
remain the appropriate conclusions.

\subsection{Why the fixed-law approximation does not extend to arbitrary
canonical progressive controls}
\phantomsection
\addcontentsline{toc}{section}{Appendix}

In the singleton-law Brownian case considered in
Subsection~\ref{sec:Brownian-progressive-controls}, approximation
under the fixed probability measure suffices to identify the
progressive and elementary control values.  Under covariance
uncertainty, however, convergence under each fixed law need not be
uniform over the admissible laws. The following results show that this is a genuine obstruction:
elementary controls need not be dense in the unrestricted canonical
progressive class for the robust norm, and enlarging the control
class can change the robust value.

\begin{proposition}
\label{prop:progressive-robust-nondensity}
Let the noise be one-dimensional, let \(A=[0,1]\) and
\(\cQ=[1,2]\), and let \(\cA\) be the elementary canonical
control class.  Define the canonical progressively measurable control
\[
 \alpha^*_u(\omega)
 :=\mathbf 1_{\{\widehat Q_u(\omega)\leq3/2\}}.
\]
For every \(T>0\) and \(p\geq1\),
\[
 \inf_{\beta\in\cA}\sup_{P\in\cP_{\cQ}}
 E^P\!\left[\int_0^T
       |\alpha^*_u-\beta_u|^p\,du\right]
 =\frac{T}{2^p}>0.
\]
Thus elementary canonical controls are not dense in the class of all
\(A\)-valued canonical progressively measurable controls for this
robust integral distance.
\end{proposition}

\begin{proof}
Fix \(\beta\in\cA\).  Applying the covariance-response realisation
construction of Lemma~\ref{lem:elementary-covariance-realisation}
(for example, with auxiliary state equation \(dY_u=-Y_u\,du+dB_u\))
to the continuous response \(\widehat q(q,x,a)=1+a\) gives
\(P^\beta\in\cP_{\cQ}\) such that
\[
 \widehat Q_u=1+\beta_u
 \quad du\times P^\beta\text{-a.e. on }[0,T].
\]
Here the response is independent of the state: on each holding
interval the already chosen action determines the constant covariance,
and these pieces are concatenated.  Consequently,
\[
 \alpha^*_u=\mathbf 1_{\{\beta_u\leq1/2\}},\qquad
 |\alpha^*_u-\beta_u|\geq\frac12
 \quad du\times P^\beta\text{-a.e.}
\]
This proves the lower bound.  The constant elementary control
\(\beta_u\equiv1/2\) gives
\(|\alpha^*_u-\beta_u|=1/2\) on every path and at every time,
which proves the reverse inequality.
\end{proof}

\begin{example}
\label{ex:progressive-robust-value-gap}
Keep \(A=[0,1]\) and \(\cQ=[1,2]\), fix \(\kappa>0\), and
consider the one-dimensional controlled equation
\begin{align*}
 dX_u^\gamma
 ={}&\bigl[-\kappa X_u^\gamma
        +4\gamma_u(1-\gamma_u)-3(2\gamma_u-1)\bigr]du\\
 &+2(2\gamma_u-1)\,d\langle B\rangle_u+dB_u,
 \qquad X_0^\gamma=0.
\end{align*}
The coefficients are continuous in the action, globally Lipschitz in
the state uniformly in the action, and independent of calendar time.
Under \(P\in\cP_{\cQ}\), write \(Q_u^P=\widehat Q_u\).  Set
\[
 g(a,q):=4a(1-a)+(2a-1)(2q-3).
\]
The effective drift is \(-\kappa x+g(a,q)\), and the effective
covariance is \(q\in[1,2]\).  Since \(|g(a,q)|\leq5/4\),
\[
 2x\bigl(-\kappa x+g(a,q)\bigr)+q
 \leq-\kappa x^2+2+\frac{25}{16\kappa}.
\]
Thus the uniformly parabolic and all-policy Lyapunov assumptions hold.
In particular, this example does not rely on degeneracy or lack of
stability.

For each fixed admissible law, the equation has a unique strong
solution for every bounded canonical progressive control, by the
explicit variation-of-constants formula
\[
 X_T^\gamma
 =\int_0^T e^{-\kappa(T-u)}g(\gamma_u,Q_u^P)\,du
  +\int_0^T e^{-\kappa(T-u)}\,dB_u.
\]
Let \(\mathcal U_{\rm prog}\) denote all \([0,1]\)-valued
controls progressive for the raw canonical filtration, and define,
for terminal cost \(\varphi(x)=x\) and zero running cost,
\[
 V_{\rm el}(T):=\inf_{\beta\in\cA}\sup_{P\in\cP_{\cQ}}E^P[X_T^\beta],
 \qquad
 V_{\rm prog}(T):=\inf_{\gamma\in\mathcal U_{\rm prog}}
                \sup_{P\in\cP_{\cQ}}E^P[X_T^\gamma].
\]
Then
\[
 V_{\rm el}(T)=\frac{1-e^{-\kappa T}}{\kappa},
 \qquad V_{\rm prog}(T)=0.
\]
Indeed, for every elementary \(\beta\), the law \(P^\beta\) in
the preceding proof satisfies \(Q_u^{P^\beta}=1+\beta_u\), and
\[
 g(a,1+a)=4a(1-a)+(2a-1)^2=1.
\]
Hence \(\sup_P E^P[X_T^\beta]\geq(1-e^{-\kappa T})/\kappa\).
For the constant control \(\beta\equiv0\), one has
\(g(0,q)=3-2q\leq1\), and the constant-covariance law \(q=1\)
attains this bound.  This proves the elementary-value formula.

For every progressive \(\gamma\), the constant-covariance law
\(q=3/2\) gives \(g(\gamma_u,3/2)=4\gamma_u(1-\gamma_u)\geq0\),
so \(\sup_P E^P[X_T^\gamma]\geq0\).  The control \(\alpha^*\)
from the preceding proposition satisfies, under every admissible law,
\[
 g(\alpha^*_u,Q_u^P)=-|2Q_u^P-3|\leq0.
\]
Thus its worst-case expected terminal cost is zero, attained by the
constant-covariance law \(q=3/2\).  This proves the progressive-value
formula.

The above calculations show that there is also a direct obstruction
to uniform state approximation.  Under \(P^\beta\), drive
\(X^\beta\) and \(X^{\alpha^*}\) by the same canonical noise.
Their stochastic terms cancel, and
\[
 X_T^\beta-X_T^{\alpha^*}
 =\int_0^T e^{-\kappa(T-u)}
       \bigl(1+|2\beta_u-1|\bigr)\,du
 \geq\frac{1-e^{-\kappa T}}{\kappa},
 \qquad P^\beta\text{-a.s.}
\]
Consequently,
\[
 \inf_{\beta\in\cA}\sup_{P\in\cP_{\cQ}}
 E^P\!\left[|X_T^\beta-X_T^{\alpha^*}|^2\right]
 \geq
 \left(\frac{1-e^{-\kappa T}}{\kappa}\right)^2>0.
\]
Thus elementary controls cannot approximate the state generated by
\(\alpha^*\) uniformly over the admissible covariance laws, even at
a fixed terminal time.

The same example also separates average costs.  For running cost
\(f(x,a)=x\) and zero terminal cost, put
\(c_\kappa(T):=(1-e^{-\kappa T})/\kappa\).
The preceding pointwise-in-time bounds, attained by the same controls
and covariance laws, yield
\begin{align*}
 \inf_{\beta\in\cA}\sup_{P\in\cP_{\cQ}}
 E^P\!\left[\int_0^T X_u^\beta\,du\right]
 &=\frac{T-c_\kappa(T)}{\kappa},\\
 \inf_{\gamma\in\mathcal U_{\rm prog}}\sup_{P\in\cP_{\cQ}}
 E^P\!\left[\int_0^T X_u^\gamma\,du\right]
 &=0.
\end{align*}
Their limiting optimized average costs are therefore \(1/\kappa\)
and \(0\), respectively.
\end{example}

\begin{remark}
The preceding statements concern enlargement to all canonical
progressively measurable controls.  They do not exclude approximation
within a smaller class defined as a completion of elementary controls
in a robust norm.  For example, the admissible controls in
\cite{HuJi2017} belong to \(M_G^2\), a completion of simple processes,
rather than being arbitrary progressively measurable controls.
Thus the obstruction above shows that the fixed-law Brownian
approximation cannot in general be transferred to the full canonical
progressive class under covariance uncertainty: even the optimized
values can differ.  It does not exclude extensions based on additional
regularity, approximation or information-structure assumptions.

Replacing progressive measurability by predictability does not remove
the obstruction.  On the continuous canonical path space equipped
with its raw filtration, the two notions coincide; see
\cite[Proposition~3.4 of the arXiv version]{ClaisseTalayTan2016}.
In particular, the covariance-dependent control in the example above
is already predictable.  The relevant restriction is therefore the
sample-and-hold structure, which fixes the action on each holding
interval, rather than predictability alone.
\end{remark}

\section*{Acknowledgements}
\phantomsection
\addcontentsline{toc}{section}{Acknowledgements}
The author 
acknowledges financial support from the 
EPSRC
Established Career Fellowship (EP/S005293/2).
\section*{AI declarations}
During the preparation of this work, the author used ChatGPT (OpenAI) for auditing the manuscript and suggesting improvements. The research, ideas, 
mathematical results and arguments originated with the author.  All
AI-assisted material was independently reviewed, verified and revised
by the author, who takes full responsibility for the content of the
manuscript.


\begin{thebibliography}{99}
\footnotesize
\setlength{\itemsep}{1pt}
\setlength{\parskip}{0pt}

\bibitem{AchdouLauriereLions2021}
{Y.~Achdou, M.~Lauri\`ere and P.-L.~Lions,
Optimal control of conditioned processes with feedback controls,
\emph{J. Math. Pures Appl.} \textbf{148} (2021), 308--341.}

\bibitem{ArtznerDelbaenEberHeath1999}
P.~Artzner, F.~Delbaen, J.-M.~Eber and D.~Heath,
Coherent measures of risk,
\emph{Math. Finance} \textbf{9} (1999), 203--228.

\bibitem{BanerjeeGarofalo2015}
A.~Banerjee and N.~Garofalo,
Boundary behavior of nonnegative solutions of fully nonlinear parabolic
equations,
\emph{Manuscripta Math.} \textbf{146} (2015), 201--222.
Preprint version cited: arXiv:1212.5985v1 (2012).

\bibitem{BertsekasShreve1996}
D.~P. Bertsekas and S.~E. Shreve,
\emph{Stochastic Optimal Control: The Discrete-Time Case},
Athena Scientific, Belmont, MA, 1996.

\bibitem{CheriditoSonerTouziVictoir2007}
{P.~Cheridito, H.~M. Soner, N.~Touzi and N.~Victoir,
Second-order backward stochastic differential equations and fully nonlinear
parabolic PDEs,
\emph{Comm. Pure Appl. Math.} \textbf{60} (2007), no.~7, 1081--1110.}

\bibitem{ClaisseTalayTan2016}
J.~Claisse, D.~Talay and X.~Tan,
A pseudo-Markov property for controlled diffusion processes,
\emph{SIAM J. Control Optim.} \textbf{54} (2016), no.~2,
1017--1029.
Preprint version cited: arXiv:1501.03939v1 (2015).

\bibitem{CIL1992}
M.~G. Crandall, H. Ishii and P.-L. Lions,
User's guide to viscosity solutions of second order partial
 differential equations,
\emph{Bull. Amer. Math. Soc. (N.S.)} \textbf{27} (1992), 1--67.

\bibitem{CossoPham2019}
{A.~Cosso and H.~Pham,
Zero-sum stochastic differential games of generalized McKean--Vlasov type,
\emph{J. Math. Pures Appl.} \textbf{129} (2019), 180--212.}

\bibitem{DaPratoRockner2007}
{G.~Da Prato and M.~R\"ockner,
A note on evolution systems of measures for time-dependent stochastic
differential equations,
in \emph{Seminar on Stochastic Analysis, Random Fields and Applications V},
Progr. Probab., vol.~59, Birkh\"auser, Basel, 2007, 115--122.}

\bibitem{Delbaen2006}
{F.~Delbaen,
The structure of $m$-stable sets and in particular of the set of risk-neutral
measures,
in \emph{In Memoriam Paul-Andr\'e Meyer: S\'eminaire de Probabilit\'es XXXIX},
Lecture Notes in Math., vol.~1874, Springer, Berlin, 2006, 215--258.}

\bibitem{DenisHuPeng2011}
L.~Denis, M.~Hu and S.~Peng,
Function spaces and capacity related to a sublinear expectation:
application to $G$-Brownian motion paths,
\emph{Potential Anal.} \textbf{34} (2011), 139--161.

\bibitem{ElKarouiPengQuenez1997}
{N.~El Karoui, S.~Peng and M.~C. Quenez,
Backward stochastic differential equations in finance,
\emph{Math. Finance} \textbf{7} (1997), no.~1, 1--71.}

\bibitem{ElKarouiTan2015}
N.~El Karoui and X.~Tan,
Capacities, measurable selection and dynamic programming. Part II:
application in stochastic control problems,
{arXiv:1310.3364v2, 2015.}

\bibitem{EpsteinSchneider2003}
{L.~G. Epstein and M.~Schneider,
Recursive multiple-priors,
\emph{J. Econom. Theory} \textbf{113} (2003), no.~1, 1--31.}

\bibitem{FengQuZhao2026}
{C.~Feng, B.~Qu and H.~Zhao,
Entrance measures and dynamics for time-inhomogeneous McKean--Vlasov stochastic
differential equations,
\emph{Trans. Amer. Math. Soc.}, advance online publication, 2026.
\href{https://doi.org/10.1090/tran/9809}{doi:10.1090/tran/9809}.}

\bibitem{FengZhao2020}
{C.~Feng and H.~Zhao,
Random periodic processes, periodic measures and ergodicity,
\emph{J. Differential Equations} \textbf{269} (2020), no.~9, 7382--7428.}

\bibitem{FengZhao2021}
C.~Feng and H.~Zhao,
Ergodicity of sublinear Markovian semigroups,
\emph{SIAM J. Math. Anal.} \textbf{53} (2021), 5646--5681.

\bibitem{FengZhaoZhong2023}
{C.~Feng, H.~Zhao and J.~Zhong,
Existence of geometric ergodic periodic measures of stochastic differential
equations,
\emph{J. Differential Equations} \textbf{359} (2023), 67--106.}

\bibitem{GilboaSchmeidler1989}
{I.~Gilboa and D.~Schmeidler,
Maxmin expected utility with non-unique prior,
\emph{J. Math. Econom.} \textbf{18} (1989), no.~2, 141--153.}

\bibitem{HairerMattingly2011}
M.~Hairer and J.~C. Mattingly,
Yet another look at Harris' ergodic theorem for Markov chains,
in \emph{Seminar on Stochastic Analysis, Random Fields and Applications VI},
Progr. Probab. 63, Birkh\"auser/Springer, Basel, 2011, 109--117.

\bibitem{HuJi2017}
M.~Hu and S.~Ji,
Dynamic programming principle for stochastic recursive optimal control
problem driven by a \(G\)-Brownian motion,
\emph{Stochastic Process. Appl.} \textbf{127} (2017), no.~1, 107--134.

\bibitem{HuWangErgodic2018}
M.~Hu and F.~Wang,
Ergodic BSDEs driven by \(G\)-Brownian motion and applications,
\emph{Stoch. Dyn.} \textbf{18} (2018), no.~6, 1850050.

\bibitem{HuWangInfinite2018}
M.~Hu and F.~Wang,
Stochastic optimal control problem with infinite horizon driven by
\(G\)-Brownian motion,
\emph{ESAIM Control Optim. Calc. Var.} \textbf{24} (2018), no.~2, 873--899.

\bibitem{Kechris1995}
A.~S. Kechris,
\emph{Classical Descriptive Set Theory},
Graduate Texts in Mathematics, vol.~156,
Springer-Verlag, New York, 1995.

\bibitem{Krylov1987}
N.~V. Krylov,
\emph{Nonlinear Elliptic and Parabolic Equations of the Second Order},
Reidel, Dordrecht, 1987.

\bibitem{Krylov2023LdPlusOne}
N.~V. Krylov,
On diffusion processes with drift in \(L_{d+1}\),
\emph{Potential Anal.} \textbf{59} (2023), 1013--1037.
Preprint version cited: arXiv:2102.11465v1 (2021).

\bibitem{KrylovSafonov1979}
N.~V. Krylov and M.~V. Safonov,
An estimate for the probability of a diffusion process hitting a set of
positive measure,
\emph{Soviet Math. Dokl.} \textbf{20} (1979), 253--255.

\bibitem{KunzeLorenziLunardi2010}
{M.~Kunze, L.~Lorenzi and A.~Lunardi,
Nonautonomous Kolmogorov parabolic equations with unbounded coefficients,
\emph{Trans. Amer. Math. Soc.} \textbf{362} (2010), no.~1, 169--198.}

\bibitem{Lieberman1996}
G.~M. Lieberman,
\emph{Second Order Parabolic Differential Equations},
World Scientific, Singapore, 1996.

\bibitem{Lorenzi2013}
{L.~Lorenzi,
Nonautonomous Kolmogorov equations in the whole space: a survey on recent
results,
\emph{Discrete Contin. Dyn. Syst. Ser. S} \textbf{6} (2013), no.~3, 731--760.}

\bibitem{LuZhao2026}
J.~Lu and H.~Zhao,
{
Random periodic solutions of SDEs driven by $\alpha$-stable processes:
existence, uniqueness and numerical approximations,
\emph{Communications in Nonlinear Science and Numerical Simulation},
\textbf{152} (2026), Part E, 109451.
}

\bibitem{MaZhao2026}
C.~Ma and H.~Zhao,
Ergodic stochastic optimal control problems,
\emph{arXiv preprint arXiv:2608.15243}, 2026.

\bibitem{Nutz2012}
M.~Nutz,
Pathwise construction of stochastic integrals,
\emph{Electron. Commun. Probab.} \textbf{17} (2012), no.~24, 1--7.

\bibitem{NutzVanHandel2013}
M.~Nutz and R.~van Handel,
Constructing sublinear expectations on path space,
\emph{Stochastic Process. Appl.} \textbf{123} (2013), 3100--3121.

\bibitem{Peng2007}
S.~Peng,
$G$-expectation, $G$-Brownian motion and related stochastic calculus of
It\^o type,
in \emph{Stochastic Analysis and Applications}, Abel Symp. 2,
Springer, Berlin, 2007, 541--567.

\bibitem{Peng2019}
{S.~Peng,
\emph{Nonlinear Expectations and Stochastic Calculus under Uncertainty:
With Robust CLT and $G$-Brownian Motion},
Probability Theory and Stochastic Modelling, vol.~95,
Springer, 2019.}

\bibitem{PossamaiTanZhou2018}
{D.~Possama\"i, X.~Tan and C.~Zhou,
Stochastic control for a class of nonlinear kernels and applications,
\emph{Ann. Probab.} \textbf{46} (2018), no.~1, 551--603.}

\bibitem{Riedel2004}
F.~Riedel,
Dynamic coherent risk measures,
\emph{Stochastic Process. Appl.} \textbf{112} (2004), 185--200.

\bibitem{ShenStannatObermayer2013}
Y.~Shen, W.~Stannat and K.~Obermayer,
Risk-sensitive Markov control processes,
\emph{SIAM J. Control Optim.} \textbf{51} (2013), no.~5, 3652--3672.

\bibitem{Sirbu2014Elementary}
M.~S{\^\i}rbu,
Stochastic Perron's method and elementary strategies for zero-sum
differential games,
\emph{SIAM J. Control Optim.} \textbf{52} (2014), no.~3,
1693--1711.

\bibitem{Sirbu2014Robust}
M.~S{\^\i}rbu,
A note on the strong formulation of stochastic control problems
with model uncertainty,
\emph{Electron. Commun. Probab.} \textbf{19} (2014), no.~81,
1--10.

\bibitem{Sirbu2015Markov}
M.~S{\^\i}rbu,
Asymptotic Perron's method and simple Markov strategies in
stochastic games and control,
\emph{SIAM J. Control Optim.} \textbf{53} (2015), no.~4,
1713--1733.

\bibitem{SonerTouziZhang2011}
{H.~M. Soner, N.~Touzi and J.~Zhang,
Martingale representation theorem for the $G$-expectation,
\emph{Stochastic Process. Appl.} \textbf{121} (2011), no.~2, 265--287.}

\bibitem{SonerTouziZhang2012}
{H.~M. Soner, N.~Touzi and J.~Zhang,
Wellposedness of second order backward SDEs,
\emph{Probab. Theory Related Fields} \textbf{153} (2012), no.~1--2, 149--190.}

\bibitem{StroockVaradhan1979}
D.~W. Stroock and S.~R.~S. Varadhan,
\emph{Multidimensional Diffusion Processes},
Springer, Berlin, 1979.

\bibitem{Wang1992}
L.~Wang,
On the regularity theory of fully nonlinear parabolic equations. I,
\emph{Comm. Pure Appl. Math.} \textbf{45} (1992), 27--76.

\bibitem{ZhaoLearningAmbiguity2026}
H.~Zhao,
Learning covariance ambiguity in controlled \(G\)-diffusions:
offline certification and online adaptation,
companion manuscript, 2026.

\bibitem{ZhaoZhao2026}
D.~Zhao and H.~Zhao,
Weighted Harris--oscillation contraction and invariant expectations for
nonlinear Markov semigroups, \emph{arXiv:2610.39492}, 2026.

\end{thebibliography}
\end{document}